\documentclass[reqno]{amsart}
\usepackage{mathrsfs}
\usepackage{amsmath,amsfonts,amssymb,amsthm}
\usepackage{upgreek}
\usepackage[usenames,dvipsnames]{xcolor}
\usepackage{enumitem}
\usepackage{graphicx}
\usepackage{overpic}
\usepackage{cancel}
\usepackage[outdir=./]{epstopdf}
\usepackage[colorlinks=true]{hyperref} %para PDFLaTeX
\usepackage{url}
\usepackage[utf8]{inputenc}
\usepackage{booktabs,multirow}
\usepackage[font=footnotesize,width=\textwidth]{caption}
\hypersetup{linkcolor=Violet,citecolor=PineGreen}
\newtheorem{teor}{Theorem}[section]
\newtheorem{lema}[teor]{Lemma}
\newtheorem{prop}[teor]{Proposition}
\newtheorem{coro}[teor]{Corollary}
\theoremstyle{definition}
\newtheorem{defi}[teor]{Definition}

\newtheorem{nota}[teor]{Remark}
\newtheorem{notas}[teor]{Remarks}
\numberwithin{equation}{section}
\newcommand{\R}{\mathbb R}

\newcommand{\Z}{\mathbb{Z}}
\newcommand{\Q}{\mathbb{Q}}

\newcommand{\N}{\mathbb{N}}

\newcommand{\s}{\mathbb{S}}
\newcommand{\mA}{\mathcal{A}}
\newcommand{\mB}{\mathcal{B}}
\newcommand{\mC}{\mathcal{C}}
\newcommand{\mD}{\mathcal{D}}

\newcommand{\mK}{\mathcal{K}}
\newcommand{\mM}{\mathcal{M}}

\newcommand{\mR}{\mathcal{R}}

\newcommand{\mU}{\mathcal{U}}
\newcommand{\mV}{\mathcal{V}}
\newcommand{\ep}{\varepsilon}

\newcommand{\mI}{\mathcal{I}}
\newcommand{\mJ}{\mathcal{J}}
\newcommand{\mL}{\mathcal{L}}

\newcommand{\W}{\Omega}
\newcommand{\WR}{\W\times\R}
\newcommand{\RR}{\R\times\R}
\newcommand{\w}{\omega}
\newcommand{\lb}{\lambda}

\newcommand{\ma}{\mathfrak{a}}
\newcommand{\mb}{\mathfrak{b}}
\newcommand{\mc}{\mathfrak{c}}
\newcommand{\md}{\mathfrak{d}}
\newcommand{\mr}{\mathfrak{r}}

\newcommand{\ml}{\mathfrak{l}}
\newcommand{\mm}{\mathfrak{m}}
\newcommand{\mn}{\mathfrak{n}}
\newcommand{\muk}{\mathfrak{u}}
\newcommand{\mf}{\mathfrak{f}}
\newcommand{\mg}{\mathfrak{g}}
\newcommand{\mk}{\mathfrak{k}}
\newcommand{\mh}{\mathfrak{h}}
\newcommand{\ms}{\mathfrak{s}}

\newcommand{\upalfa}{$\upalpha$}
\newcommand{\upomeg}{$\upomega$}

\newcommand{\G}{\Gamma}

\newcommand{\wit}{\widetilde}

\newcommand{\ws}{\w{\cdot}s}
\newcommand{\wt}{\w{\cdot}t}
\newcommand{\bwt}{\bar\w{\cdot}t}
\newcommand{\bws}{\bar\w{\cdot}s}

\newcommand{\n}[1]{\left\|#1\right\|}

\newcommand{\lsm}{\left[\begin{smallmatrix}}
\newcommand{\rsm}{\end{smallmatrix}\right]}

\newcommand{\merg}{\mathfrak{M}_\mathrm{erg}(\W,\sigma)}

\begin{document}
\title[Concave-convex ODEs: from bifurcation to critical transitions]
{Concave-convex nonautonomous scalar ordinary differential equations: from bifurcation theory to critical transitions.}
\author[J. Due\~{n}as]{Jes\'{u}s Due\~{n}as}
\author[C. N\'{u}\~{n}ez]{Carmen N\'{u}\~{n}ez}
\author[R. Obaya]{Rafael Obaya}
\address{{\rm(J. Due\~{n}as, C.N\'{u}\~{n}ez, R. Obaya)}. Departamento de Matem\'{a}tica Aplicada, Universidad de Va\-lladolid, Paseo Prado de la Magdalena 3-5, 47011 Valladolid, Spain. The authors are members of IMUVA: Instituto de Investigaci\'{o}n en Matem\'{a}ticas, Universidad de Valladolid.}
\email[J.~Due\~{n}as]{jesus.duenas@uva.es}
\email[C.~N\'{u}\~{n}ez]{carmen.nunez@uva.es}
\email[R.~Obaya]{rafael.obaya@uva.es}
\thanks{All the authors were supported by Ministerio de Ciencia, Innovaci\'{o}n y Universidades (Spain)
under project PID2021-125446NB-I00 and by Universidad de Valladolid under project PIP-TCESC-2020.
J.~Due\~{n}as was also supported by Ministerio de Universidades (Spain) under programme FPU20/01627.}
\date{}
\begin{abstract}
A mathematical modeling process for phenomena with a single state variable that
attempts to be realistic must be given by a scalar nonautonomous differential equation
$x'=f(t,x)$ that is concave with respect to the state variable $x$ in some regions of its domain
and convex in the complementary zones. This article takes the first step towards developing a
theory to describe the corresponding dynamics: the case in
which $f$ is concave on the region $x\ge b(t)$ and convex on $x\le b(t)$, where $b$ is
a $C^1$ map, is considered. The different long-term dynamics that may appear are analyzed
while describing the bifurcation diagram for $x'=f(t,x)+\lb$. The results are used to
establish conditions on a concave-convex map $h$ and a nonnegative map $k$ ensuring the existence
of a value $\rho_0$ giving rise to the unique critical transition for the parametric family of
equations $x'=h(t,x)-\rho\,k(t,x)$, which is assumed to approach $x'=h(t,x)$ as time decreases,
but for which no conditions are assumed on the future dynamics. The developed theory is justified
by showing that concave-convex models fit correctly some laboratory experimental data,
and applied to describe a population dynamics model for which a large enough increase on the
peak of a temporary higher predation causes extinction.
\end{abstract}
\renewcommand{\subjclassname}{\textup{2020} Mathematics Subject Classification}

\maketitle
%%%%%%%%%%%%%%%%%%%%%%%%%%%%%%%%%%%%%%%%%%%%%%%%%%%%%%%%%%%%%%%%%%%%%%%%%%%%%%%%
%%%%%%%%%%%%%%%%%%%%%%%%%%%%%%%%%%%%%%%%%%%%%%%%%%%%%%%%%%%%%%%%%%%%%%%%%%%%%%%%
%%%%%%%%%%%%%%%%%%%%%%%%%%%%%%%%%%%%%%%%%%%%%%%%%%%%%%%%%%%%%%%%%%%%%%%%%%%%%%%%
%%%%%%%%%%%%%%%%%%%%%%%%%%%%%%%%%%%%%%%%%%%%%%%%%%%%%%%%%%%%%%%%%%%%%%%%%%%%%%%%
\section{Introduction}
This paper deals with scalar nonautonomous ordinary differential equations (ODEs
hereafter) $y'=f(t,y)$ of concave-convex type: the function $f\colon\RR\to\R$ and its
derivative $f_y$ with respect to the state variable $y$ are bounded and uniformly
continuous on $\R\times\mJ$ for each compact interval $\mJ\subset\R$;
and there exists a $C^1$ function $\kappa\colon\R\to\R$ with $\kappa$ and $\kappa'$
bounded and uniformly continuous such that, for each $t$ fixed, the map
$y\mapsto f(t,y)$ is concave when $y\ge\kappa(t)$ and convex when $y\le\kappa(t)$.
Under an additional coercive condition, we study the possibilities for the
long term behaviour of the solutions of the equation, and relate them to
local and global properties of the (nonautonomous) bifurcation diagram of
$y'=f(t,y)+\lb$. Among the several possible scenarios, a single or double saddle-node
bifurcation of hyperbolic solutions appears.
These bifurcation patterns bear significant influence on the description of relevant problems 
in applied sciences: see, for example, \cite{aajqw}, \cite{goldsztein}, \cite{haberman}, 
\cite{jung}, \cite{mitsui} and \cite{zaliapin}. We apply the conclusions of our study to mathematical
modeling and critical transitions of population dynamics problems.
Our methods combine tools from topological dynamics and ergodic theory.
The conclusions of our theoretical study, combined with numerical analysis,
can provide qualitative information in models in biology \cite{hill}, ecology \cite{scheffer,schefferbbbc},
circuit theory \cite{nagumoay,grasman} and neural networks \cite{wu,hoppen}, among other areas of interest.

The initial motivation for investigating this problem is twofold. The works of
Anagnostopoulou and J\"{a}ger~\cite{anagjager1},
Anagnostopoulou et al.~\cite{anapotras},
Due\~{n}as et al.~\cite{dno1,dno2,dno3,dno5},
Fuhrmann~\cite{fuhrmann},
Longo et al.~\cite{lnor,lno2,lno3},
N\'{u}\~{n}ez and Obaya~\cite{nuob6}, Nu\~{n}ez et al.~\cite{nuos4} and Remo et al.~\cite{remo1},
as well as some references therein, extensively investigate the
same questions in two relevant types of nonautonomous scalar ODEs: first, the case where $f$
is globally concave with respect to the state variable $x$ (i.e., with $x\mapsto f(t,x)$
concave for all $t\in\R$); and, second, the case where $f$ is globally d-concave
(which means $f_x$ is globally concave, always with respect to $x$). All these works
analyze the occurrence of nonautonomous bifurcation points of local saddle-node type.
This work extends the results obtained in the d-concave case,
since every regular enough d-concave and coercive function is concave-convex. Thus, the
field of applicability of this type of results is widened. The second motivation is to
take a first step in the description of the long-term behavior of more realistic systems,
which in general are given by functions that change several times from concavity to convexity
and vice versa. This paper analyzes this framework assuming the existence of only two zones,
separated by the graph of a $C^1$ function.

Our topological and ergodic methods require to include the nonautonomous ODE $y'=f(t,y)$ whose
long term behaviour we want to describe in a family of ODEs related by a continuous flow on
a compact metric space. So, we work with a compact metric space $\W$, with a continuous global
flow $\sigma\colon\R\times\W\to\W,\,(t,\w)\mapsto\wt:=\sigma(t,\w)$ on $\W$,
and with a continuous map $\mf\colon\WR\to\R$ such that the partial derivative $\mf_y$
with respect to the second variable exists and is continuous on $\WR$.
If there exists a point $\w_0\in\W$ with $f(t,y)=\mf(\w_0{\cdot}t,y)$, then our initial
equation is one of the family $\{y'=\mf(\wt,y)\,|\;\w\in\W\}$. This inclusion
is always possible under some boundedness and regularity conditions on $f$ which ensure
the existence and compactness of its hull $\W_f:=\text{closure}\{f_t\,|\;t\in\R\}\subset C(\WR,\R)$
in the compact-open topology, with $f_t(s,y):=f(t+s,y)$, and the continuity of the time
shift flow $\sigma_f(t,\w):=\w_t$. In this case, we define $\mf(\w,y)=\w(0,y)$, so that
$\mf(f{\cdot}t,y)=f(t,y)$. This {\em hull construction\/} is more detailed in
Section \ref{subsec:2process}. The resulting flow on $\W_f$ is transitive: the
$\sigma_f$-orbit of $f$ is obviously dense. But our initial analysis in this work does
not assume this condition on $(\W,\sigma)$, which makes it considerably more general.
The great advantage of the collective formulation is that the family $\{y'=\mf(\wt,x)\,|\;\w\in\W\}$
defines a (possibly local) continuous {\em skewproduct\/} flow on $\WR$, sending a point $(\w,y)$ in time $t$
to $(\wt,v(t,\w,y))$, where $v(t,\w,y)$ solves $y'=\mf(\wt,y)$ and satisfies $v(0,\w,y)=y$.

As indicated above, we assume the existence of a map $\mb\colon\W\to\R$ such that
$t\mapsto\mb(\wt)$ is $C^1$ on $\R$, the map $\mb'\colon\W\to\R$ given by
$\mb'(\w):=(d/dt)\,\mb(\wt)|_{t=0}$ is continuous, and, for each $\w\in\W$
fixed, the map $x\mapsto\mf(\w,x)$ is concave when $x\ge\mb(\w)$ and convex when
$x\le\mb(\w)$. The change of state variable $x(t)=y(t)-\mb(\wt)$ takes our initial family
to $\{x'=\mh(\wt,x)\,|\;\w\in\W\}$ for $\mh(\w,x):=\mf(\w,x+\mb(\w))-\mb'(\w)$.
We rewrite each equation $x'=\mh(\wt,x)$ as
\begin{equation}\label{eq:1conclin}
 x'=\mc(\wt)+\md(\wt)\,x+\mg(\wt,x)
\end{equation}
for $\mc(\w):=\mh(\w,0)=\mf(\w,\mb(\w))-\mb'(\w)$, $\md(\w):=\mh_x(\w,0)=\mf_y(\w,\mb(\w))$
and $\mg(\w,x):=\mh(\w,x)-\mc(\w)-\md(\w)\,x=\mf(\w,x+\mb(\w))-\mf(\w,\mb(\w))-\mf_y(\w,\mb(\w))\,x$.
In this way, for any $\w\in\W$,
$\mg(\w,0)=0$, $\mg_x(\w,0)=0$, and the map $x\mapsto\mg(\w,x)$ is convex on $(-\infty,0)$ and
concave on $(0,\infty)$.
%Let us represent by $\tau$ the skewproduct flow induced on $\WR$
%by the family of equations \eqref{eq:1conclin} (varying on $\W$) HACE FALTA???.
The classical coercivity condition $\lim_{x\to\pm\infty}\mh(\w,x)=\pm\infty$,
assumed throughout the paper, implies the boundedness of all the forward orbits
of the skewproduct flow induced by \eqref{eq:1conclin} on $\WR$,
as well as the existence of a global attractor.

In fact, the initial conditions on concavity and convexity that we assume for our analysis are quite
more general than those described in the previous paragraph: the function $x\mapsto\mg(\w,x)$
is not required to be convex on $(-\infty,0)$ and concave on $(0,\infty)$ for {\em all\/} $\w\in\W$,
but only for the elements of a subset with full measure for every ergodic measure on $\W$.
We also assume these conditions to be strict for the elements of a set with positive measure for every
ergodic measure. This idea follows the path set out in \cite{dno4,dno5} in the concave and d-concave cases:
since, in many cases, the ergodic measures are concentrated in topologically small sets
(with empty interior or even of the first Baire category), this {\em concave in measure\/} or
{\em d-concave in measure} assumptions are considerably weaker than the hypotheses assumed on previous works,
despite of what the conclusions are identical, considerably broadening the field of possible
applications. And the same can be said, of course, in the present {\em concave-convex in measure\/} case.

Scalar differential equations of type \eqref{eq:1conclin} have been studied in Campos et al.~\cite{cano}
assuming $\int_\W\md(\w)\,dm\le 0$ for every ergodic measure $m$ on $\W$. When all these integrals are strictly
negative, the dynamics is governed by the presence of an attractor $\mA\subset\WR$
that is homeomorphic to the set $\W$ and globally exponentially stable.
When some of these integrals vanish, the dynamics can be significatively more complex, admitting the existence
of a chaotic global attractor. In order to extend the conclusions previously obtained, in this paper
we analyze the case in which $\int_\W\md(\w)\,dm> 0$ for every ergodic measure $m$ on $\W$.
We will see that the dynamical possibilities under this condition are significatively different
from that of \cite{cano} even assuming the most restrictive conditions on the concave-convex map
giving rise to the ODE.

Much of this work is focused on the local and global properties of the bifurcation diagram
of the parametric family
\begin{equation}\label{eq:1param}
  x'=\mc(\wt)+\md(\wt)\,x+\mg(\wt,x)+\lb,\,\quad\w\in\W
\end{equation}
(also written as $x'=\mh(\wt,x)+\lb$), in the line of the works mentioned at the
beginning. In some of the analyzed cases, the patterns we observe are quite similar to those
associated with d-concave differential equations. However, the interplay between concave and
convex effects can yield other configurations which are not possible in the d-concave in
measure case,
partly described in this paper. In fact, this case also is compatible with the existence
of more than three uniformly separated copies of the base. While we are aware of how
such examples can be explicitly constructed, these are not included in the present
discussion as they are the subject of ongoing research. Another large part of the paper
is devoted to provide theorems showing the existence of size-induced critical transitions
in parametric perturbations of nonautonomous concave-convex differential equations, in
the line of the works of Alkhayuon and Ashwin \cite{alas}, Ashwin et al.~\cite{aspw},
Due\~{n}as et al.~\cite{dno3,dno4}, Kiers and Jones~\cite{kiersjones},
Longo et al.~\cite{lno3},
Slyman et al.~\cite{slymanetal} and Wieczorek et al.~\cite{sebastian1,wieczorekxj}.
A novel aspect of the approach to critical transitions in this paper is that it offers
rigorous results without relying on assumptions about future dynamics,
which were customary in our previous papers on critical transitions.
This allows the concept of critical transition to be defined solely based on the
behavior of locally pullback attractive solutions of the transition equation, rather than on the
complete classification of the global dynamics of the transition equation.
The last part of the paper concerns the applicability of the previously developed theory,
focusing in two different and complementary aspects: to show that the concave-convex
dynamics appears naturally in population dynamics, and to describe conditions ensuring the
occurrence of a unique critical transition in one of these models.

Let us take a closer look at the content of the document while briefly explaining
its structure.

Section \ref{sec:2} contains a summary of the basic concepts and results required in the subsequent analysis,
including a short description of the hull construction. From this point, we will work both with processes
(defined by single ODEs) and with skewproduct flows (defined by families of ODEs). Among the interesting
objects for processes (resp.~for skewproduct flows), we focus on hyperbolic solutions and on the set of
globally bounded solutions (resp.~on hyperbolic invariant subsets of $\WR$ which are homeomorphic to $\W$,
called {\em hyperbolic copies of the base}, and on the global attractor). The variation in the number and
structure of these objects will determine the occurrence of bifurcation points. A key theorem
on persistence of hyperbolic solutions and hyperbolic copies of the base completes this section.

Sections \ref{sec:3} and \ref{sec:4} present the core of the bifurcation results of the paper.
The first one contains an extension of part of the main results obtained in \cite{dno5} for the case
of a family of concave in measure ODEs to a family \eqref{eq:1param} which is of {\em concave-linear in measure\/}
type: in this case, the line $x=0$ divides $\WR$ into two ``half-planes”, the function $\mg$ being null
for $x<0$ and strictly concave in measure for $x>0$. Roughly speaking, we prove that the global bifurcation
diagram is determined by the existence of a unique bifurcation value $\lb_-$ of
nonautonomous saddle-node type: an
attractor-repeller pair of hyperbolic copies of the base exists for $\lb>\lb_-$, and these solutions
approach each other as $\lb$ decreases towards $\lb_-$, resulting at this point in the
lack of uniformly separated copies of the base. In addition, $\lb_-$ is the lower
value of the parameter such that each equation of the family has bounded solutions, and
none of the equations has bounded solutions if $\lb$ is smaller than a value $\lb_*\le\lb_-$.
These results are easy translated to the {\em linear-convex in measure\/} case: if $\mg$ is null
for $x>0$ and strictly convex in measure for $x<0$, then the global bifurcation diagram is
given by a horizontal inversion of the above described one, with $\lb^+$ as
bifurcation value.

A concave-convex map $\mg$ is equal to $\mg^-\!+\mg^+$ for a concave-linear $\mg^-$ and a linear-convex
$\mg^+$. Combining the results of Section \ref{sec:3} is
one the of the key points to prove one of the main results of Section \ref{sec:4},
which establishes conditions for a concave-convex in measure family \eqref{eq:1param}
ensuring that the global bifurcation diagram is determined
by the occurrence of two bifurcation values $\lb_-<\lb^+$: they are characterized by local
structures and behaviors of the copies of the base ``around them" that reproduce those described for
the concave-linear and linear-convex case. This short description shows the previously
mentioned resemblance with the most representative global bifurcation diagram
obtained in \cite{dno5} in the d-concave in measure case. But, in fact, in the situation of
this paper, the dynamics on $(\lb_-,\lb^+)$ can be highly more complex even in this
``simplest" case.
In addition, there are many more possibilities for the bifurcation diagram. The description of
some common properties for all the cases, as well as of certain additional conclusions
which can be reached when some of the main objects of the concave-linear and linear-convex
cases fall in a frequently appearing d-concavity band, are also included in
Section \ref{sec:4}. A dynamical classification in five different cases given
under the additional assumption of transitivity of the flow on $\W$, and examples of
all these cases, complete this part.

From Section \ref{sec:5}, we assume the second derivative $\mg_{xx}$ to exist and
be continuous on $\WR$. This section contains the results on critical transitions
outlined above. They refer to a single ODE
\begin{equation}\label{eq:1ct}
 x'=h(t,x)-\rho\,k(t,x)\,,
\end{equation}
where the hull extension $\mh$ of $h$ provides a concave-convex in measure family
$x'=\mh(\wt,x)$. We establish conditions on $\mh$ and $k$ ensuring the existence of a unique
critical value $\rho_0$. One of these hypotheses is $\lim_{t\to-\infty}k(t,x)=0$
uniformly on compact subsets of $\R$ $\W$, which roughly speaking means
that $x'=h(t,x)$ acts as ``past equation" of \eqref{eq:1ct}. Our conditions also include a
convenient hyperbolic structure on $x'=\mh(\wt,x)$, with three hyperbolic
copies of the base determining the global dynamics, which in particular ensures that
the upper and lower bounded solutions of \eqref{eq:1ct} are locally pullback attractive
and approach the upper and lower hyperbolic solutions of $x'=h(t,x)$ as
time decreases. The occurrence of tracking or tipping depends on the future behavior of
these two solutions: there is tracking if they remain uniformly separated and
tipping otherwise. And some additional conditions on $k$, related to its sign
and convexity properties on the area bounded by the two upper hyperbolic copies of
the base for $x'=\mh(t,x)$, guarantee tracking on an interval $[0,\rho_0]$ and
tipping on $(\rho_0,\infty)$. We insist once more in the lack of a priori hypotheses
concerning the future: in contrast to some previous studies,
we do not even assume the existence of a ``future equation".

The paper is completed with Section \ref{sec:6}. In the first part, developing ideas
from Meyer~\cite{meyer1}, the real data offered by Ris et al.~in~\cite{mosquitodata},
which measure the population growth rate in terms of the population size in an insect species,
are fitted using concave-convex splines. This is a representative fact of what we consider
a much more general property: despite the fact that many of the usually analyzed scalar
population models are given by a map which is d-concave (or d-concave in measure) at least on
the region on which the representative solutions lie, in order to design more realistic models
it is necessary to include the possibility of changes from convexity to concavity. In the
second part of Section \ref{sec:6}, the existence of a critical transition in a
population model subject to Allee effect and predation (see Courchamp et
al.~\cite{courchamp1}) is analyzed using the results from Section
\ref{sec:5}. Our study shows the existence of a maximum threshold for the pick size
in the temporary increase in predation, beyond which the population gets extinct.
%%%%%%%%%%%%%%%%%%%%%%%%%%%%%%%%%%%%%%%%%%%%%%%%%%%%%%%%%%%%%%%%%%%%%%%%%%%%%%%%
%%%%%%%%%%%%%%%%%%%%%%%%%%%%%%%%%%%%%%%%%%%%%%%%%%%%%%%%%%%%%%%%%%%%%%%%%%%%%%%%
%%%%%%%%%%%%%%%%%%%%%%%%%%%%%%%%%%%%%%%%%%%%%%%%%%%%%%%%%%%%%%%%%%%%%%%%%%%%%%%%
%%%%%%%%%%%%%%%%%%%%%%%%%%%%%%%%%%%%%%%%%%%%%%%%%%%%%%%%%%%%%%%%%%%%%%%%%%%%%%%%
\section{Basic concepts and preliminary results}\label{sec:2}
We will manage without recalling the definitions the basic concepts of
maximal solutions of ODEs,
flows, orbits, invariant and minimal sets, \upalfa-limit and \upomeg-limit sets,
invariant and ergodic measures. They can be found in \cite{cano}, \cite{dno1},
and references therein.
Now we summarize some less known concepts, also key to this paper.
%%%%%%%%%%%%%%%%%%%%%%%%%%%%%%%%%%%%%%%%%%%%%%%%%%%%%%%%%%%%%%%%%%%%%%%%%%%%%%%%
%%%%%%%%%%%%%%%%%%%%%%%%%%%%%%%%%%%%%%%%%%%%%%%%%%%%%%%%%%%%%%%%%%%%%%%%%%%%%%%%
\subsection{Skewproduct flows given by families of scalar nonautonomous ODEs}
\label{subsec:2skew}
Let $\sigma\colon\R\times\W\to\W$, $(t,\w)\mapsto\sigma(t,\w)=:\w{\cdot}t$ be a global
continuous flow on a compact metric space $\W$.
The set of continuous functions $\mh\colon\WR\to\R$ for which the
derivative $\mh_x$ with respect to the second variable exists and is continuous is
$C^{0,1}(\WR,\R)$, and $C^{0,2}(\WR,\R)$ is the subset of maps $\mh$ for which
$\mh_{xx}$ exists and is continuous. Each $\mh\in C^{0,1}(\WR,\R)$ provides
the family of scalar nonautonomous ODEs
\begin{equation}\label{eq:2fam}
x'=\mh(\w{\cdot}t,x)\,,\quad\w\in\W\,,
\end{equation}
with maximal solutions $(\alpha_{\w,x},\beta_{\w,x})\to\R,\,t\mapsto v(t,\w,x)$
satisfying $v(0,\w,x)=x$. So, $v(t+s,\w,x)=v(t,\ws,v(s,\w,x))$ whenever the right-hand term is defined, and hence,
if $\mV:=\bigcup_{(\w,x)\in\WR}((\alpha_{\w,x},\beta_{\w,x})\times\{(\w,x)\})$, then
\begin{equation*}
 \tau\colon\mV\subseteq\R\times\WR\to\WR\,,\;\;(t,\w,x)\mapsto (\wt,v(t,\w,x))
\end{equation*}
defines a (possibly local) continuous flow on $\WR$, of {\em skewproduct} type. 
We will refer to $\W$ as the {\em skewproduct base}.
In Section \ref{subsec:2process}, we will recall how a flow of this type arises from a
single ODE, by means of the hull construction.

A {\em $\tau$-equilibrium} is a map $\mb\colon\W\to\R$ whose graph is $\tau$-invariant;
i.e., with $v(t,\w,\mb(\w))=\mb(\wt)$ for all $\w\in\W$ and $t\in\R$.
A {\em $\tau$-copy of the base\/} or {\em $\tau$-copy of $\W$} is the compact graph
of a continuous equilibrium $\mb$. We represent it by $\{\mb\}$. The prefix $\tau$
will be often omitted. The map $\mb\colon\W\to\R$ is $C^1$ {\em along the base orbits} if, for any $\w\in\W$,
the map $t\mapsto \mb_\w(t)=\mb(\w{\cdot}t)$ is $C^1$ on $\R$. In this case, $\mb'(\w):=\mb'_\w(0)$.
Such a map is a {\em global upper} (resp.~{\em lower}) {\em solution} for \eqref{eq:2fam} if
$\mb'(\w)\ge\mh(\w,\mb(\w))$ (resp.~$\mb'(\w)\le\mh(\w,\mb(\w))$) for every $\w\in\W$, and it is
{\em strict} if the inequality strict for all $\w\in\W$.

Given a bounded $\tau$-invariant set $\mB\subset\W\times[-k,k]\subset\WR$ with $k>0$ projecting onto $\W$
(i.e., such that there exists at least a point $(\w,x)\in\mB$ for each $\w\in\W$),
the maps $\w\mapsto\inf\{x\in\R\,|\;(\w,x)\in\mB\}$ and $\w\mapsto\sup\{x\in\R\,|\;(\w,x)\in\mB\}$
define $\tau$-equilibria. We will refer to these maps as
the {\em lower\/} and {\em upper equilibria of $\mB$}. They are respectively
lower and upper semicontinuous if $\mB$ is compact, and hence they are
$m$-measurable for all $m\in\merg$.

Two compact subsets $\mK_1$ and $\mK_2$ of $\W\times\R$ are {\em ordered\/}, with
$\mK_1<\mK_2$, if $x_1<x_2$ whenever there exists $\w\in\W$ such that
$(\w,x_1)\in\mK_1$ and $(\w,x_2)\in\mK_2$. Two (ordered)
bounded equilibria $\mb_1,\mb_2$ with $\mb_1<\mb_2$ (i.e., with $\mb_1(\w)<\mb_2(\w)$
for all $\w\in\W$) are {\em uniformly separated\/} if
$\inf_{\w\in\W}(\mb_2(\w)-\mb_1(\w))>0$. The $\w$-section of a set $\mC\subseteq\WR$ is $(\mC)_\w:=\{x\in\R\mid\,(\w,x)\in\mC\}$.

If there exists a compact $\tau$-invariant set $\mA\subset\WR$ such that
$\lim_{t\to\infty} \text{dist}(\mC{\cdot}t,\mA)=0$ for every bounded set $\mC\subset\WR$,
where $\mC{\cdot}t=\{(\wt,v(t,\w,x))\,|\;(\w,x)\in\mC\}$ and
\[
\text{dist}(\mC_1,\mC_2)=\sup_{(\w_1,x_1)\in\mC_1}\left(\inf_{(\w_2,x_2)\in\mC_2}
\big(\mathrm{dist}_{\WR}((\w_1,x_1),(\w_2,x_2))\big)\right)\,,
\]
then $\mA$ is the (unique) {\em global attractor for $\tau$}.

A $\tau$-copy of the base $\{\mb\}$ is
{\em hyperbolic attractive\/} if it is uniformly exponentially stable as time increases;
i.e., if there exists $\rho>0$, $k\ge 1$ and $\gamma>0$ such that: if, for any $\w\in\W$,
$|\mb(\w)-x|<\rho$, then $v(t,\w,x)$ is defined for all $t\ge 0$, and in addition
$|\mb(\wt)-v(t,\w,x)|\le k\,e^{-\gamma\,t}\,|\mb(\w)-x|$ for $t\ge 0$. In this case,
$\rho$ is a {\em radius of uniform exponential stability} and $(k,\gamma)$ is a
{\em dichotomy constant pair} of $\mb$. Changing $t\ge 0$ by $t\le 0$ and $\gamma$ by $-\gamma$
provides the definition of {\em repulsive hyperbolic} $\tau$-copy of the base.

Let $\mK\subset\WR$ be $\tau$-invariant compact set projecting onto $\W$,
and let $\mathfrak{M}_\mathrm{inv}(\mK,\tau)$ and $\mathfrak{M}_\mathrm{erg}(\mK,\tau)$
be the (nonempty) sets of the $\tau$-invariant and $\tau$-ergodic measures on $\mK$.
A value $\gamma\in\R$ is a \emph{Lyapunov exponent of $\mK$} if there exists
$(\w,x)\in\mK$ such that $\gamma=\lim_{t\to\pm\infty}(1/t)\int_0^t\mh_x(\tau(r,\w,x))\,dr$.
As explained in \cite[Seccion 2.4]{dno4}, there exists $m^l,m^u\in\merg$, an $m^l$-measurable
equilibrium $\mb^l\colon\W\to\R$ and an $m^u$-measurable equilibrium $\mb^u\colon\W\to\R$,
both with graph contained in $\mK$,
such that the lower and upper Lyapunov exponents of $\mK$ are given by
$\int_\W \mh_x(\w,\mb^l(\w))\,dm^l$ and $\int_\W \mh_x(\w,\mb^u(\w))\,dm^u$, respectively.

There is a strong connection between the hyperbolicity of a copy of the base $\{\mb\}$
and the sign of all its Lyapunov exponents: $\{\mb\}$ is hyperbolic attractive (resp.~repulsive)
if and only if all its Lyapunov exponents are strictly negative (resp.~positive).
This property is a consequence of \cite[Theorem 2.8]{dno4}. We repeat here its statement,
since it will be fundamental in some key points of several proofs.
\begin{teor}\label{th:2copia}
Let $\mK\subset\WR$ be a $\tau$-invariant compact set projecting
onto $\W$. Assume that its upper and lower equilibria
coincide (at least) on a point of each minimal subset $\mM\subseteq\W$.
Then, all the Lyapunov exponents of $\mK$
are strictly negative (resp.~positive) if and only if $\mK$ is
an attractive (resp.~repulsive) hyperbolic $\tau$-copy of the base.
In addition, if either $\mK$ (and hence $\W$) is minimal or its upper and lower equilibria
coincide on a $\tau$-invariant subset $\W_0\subseteq\W$
with $m(\W_0)=1$ for all $m\in\merg$, then the condition on its upper and lower equilibria holds.
\end{teor}
We will often restrict the family \eqref{eq:2fam} to the elements of a $\sigma$-invariant
compact subset $\W_0\subseteq\W$.
Since $C^{0,1}(\WR,\R)$ is continuously embedded into $C^{0,1}(\W_0\times\R,\R)$, it makes perfect sense
to consider the restricted skewproduct flow on $\W_0\times\R$,
whose definition is the same. In many of these cases, $\W_0$ will be the closure $\W_{\bar\w}$
of the $\sigma$-orbit $\{\bwt\,|\;t\in\R\}$ of a point $\bar\w\in\W$. We recall that the
base flow is {\em transitive\/} if there exists $\bar\w\in\W$ with $\W_{\bar\w}=\W$, and {\em minimal\/}
if $\W_\w=\W$ for all $\w\in\W$.

As usual, an open ball in $\W$ is represented by $B_\W(\w,\delta)$.
%%%%%%%%%%%%%%%%%%%%%%%%%%%%%%%%%%%%%%%%%%%%%%%%%%%%%%%%%%%%
%%%%%%%%%%%%%%%%%%%%%%%%%%%%%%%%%%%%%%%%%%%%%%%%%%%%%%%%%%%%
\subsection{A particular process}\label{subsec:2process}
A map $h\colon\RR\to\R$ is {\em admissible\/}, or $h\in C^{0,0}(\RR,\R)$,
if the restriction of $h$ to $\R\times\mJ$ is bounded
and uniformly continuous for any compact set $\mJ\subset\R$. The map $h$
is {\em $C^1$-admissible\/} (resp.~{\em $C^2$-admissible\/}), or
$h\in C^{0,1}(\RR,\R)$ (resp.~$h\in C^{0,2}(\RR,\R)$), if its
derivative $h_x$ with respect to the second variable exists and is admissible
(resp.~if $h_x$ and $h_{xx}$ exist and are admissible).

Let us take $h\in C^{0,1}(\RR,\R)$ and let $x_h(t,s,x)$ be the maximal solution of
\begin{equation}\label{eq:2proceso}
 x'=h(t,x)
\end{equation}
with $x_h(s,s,x)=x$. By uniqueness of solutions, $x_h(t,s,x_h(s,r,x))=x_h(t,r,x)$. The map
$(t,s,x)\mapsto x_h(t,s,x)$ is often called a {\em process}, or an {\em admissible process}.

Two bounded (and hence globally defined) solutions $b_1(t)$ and
$b_2(t)$ of \eqref{eq:2proceso} are {\em uniformly separated\/} if
$\inf_{t\in\R}|b_1(t)-b_2(t)|>0$.

A bounded solution
$b(t)$ of \eqref{eq:2proceso} is {\em hyperbolic attractive\/}
(resp.~{\em hyperbolic repulsive\/}) if
there exist $k\ge 1$ and $\gamma>0$ such that
$\exp\Big(\int_s^t h_x(r,b(r))\,dr\Big)\le k\,e^{-\gamma(t-s)}$ whenever $t\ge s$
(resp.~$\exp\Big(\int_s^t h_x(r,b(r))\,dr\Big)\le k\,e^{\gamma (t-s)}$ whenever $t\le s$).
An attractive (resp.~repulsive) hyperbolic solution is uniformly
exponentially stable at $+\infty$ (resp.~at $-\infty$),
as deduced from the First Approximation Theorem (see Theorem~III.2.4 of \cite{hale});
i.e., there exist a {\em radius of uniform exponential stability\/} $\rho>0$ and a {\em dichotomy
constant pair\/} $(k,\gamma)$ with $k\ge 1$ and $\gamma>0$ such that,
if $|b(s)-x|\le\rho$, then $|b(t)-x_h(t,s,x)|
\le k\,e^{-\gamma(t-s)}|b(t)-x|$ for all $t\ge s$
(resp.~$|b(t)-x_h(t,s,x)|\le k\,e^{\gamma(t-s)}|b(s)-x|$
for all $t\le s$).

Let us take as starting point the family \eqref{eq:2fam}, with $\mh\in C^{0,1}(\WR,\R)$.
We fix $\bar\w\in\W$ and define $h_{\bar\w}(t,x):=\mh(\bwt,x)$. It is easy to check that
$h_{\bar\w}\in C^{0,1}(\RR,\R)$, and hence
\begin{equation*}
 x'=h_{\bar\w}(t,x) \qquad (\text{i.e., }x'=\mh(\bwt,x)\,)
\end{equation*}
induces an admissible process. In addition, if $x_{\bar\w}(t,s,x)$ is the solution with $x_{\bar\w}(s,s,x)=x$, then
$x_{\bar\w}(t,s,x):=v(t-s,\bws,x)$. The converse procedure is the already mentioned {\em hull construction},
which will play a fundamental role in Section \ref{sec:5}, and which we briefly describe.
The {\em hull\/} $\W_h$ of an admissible map $h\colon\RR\to\R$ is
the closure of the set $\{h{\cdot}t\,|\;t\in\R\}$ on the set $C(\RR,\R)$ provided with the compact-open
topology, where $h{\cdot}t(s,x):=h(t+s,x)$. According to \cite[Theorem IV.3]{selltopdyn} and
\cite[Theorem~I.3.1]{shenyi}, $\W_h$ is a compact metric space,
the time-shift map $\sigma_h\colon\R\times\W_h\to\W_h,(t,\w)\mapsto\w{\cdot}t$
defines a global continuous flow, and the map $\mh$ given by $\mh(\w,x)=\w(0,x)$
is continuous on $\W_h\times\R$. In addition, if $h$ is $C^1$-admissible
then $\W_h\subset C^{0,1}(\RR,\R)$,
and the continuous map $\mh_x(\w,x):=\w_x(0,x)$ is the derivative of
$\mh$ with respect to $x$; and, if $h$ is $C^2$-admissible
then $\W_h\subset C^{0,2}(\RR,\R)$,
and the continuous map $\mh_{xx}(\w,x):=\w_{xx}(0,x)$ is the
second derivative of $\mh$ with respect to $x$.
Note also that $(\W_h,\sigma_h)$ is a transitive flow, since
the $\sigma_h$-orbit of the point $h\in\W_h$ is dense.

For the reader's convenience, we complete this section with the next fundamental persistence result,
very often used in the proofs of our results. The first assertion, concerning processes,
is \cite[Theorem 2.3]{dno4}, and the interested reader can find there suitable references for its proof.
This proof can be easily adapted
to prove the second assertion, about skewproduct flows. For $h\in C^{0,1}(\RR,\R)$, we
denote $\n{h}_{1,\kappa}:=\sup_{(t,x)\in\R\times[-\kappa,\kappa]} |h(t,x)|+
\sup_{(t,x)\in\R\times[-\kappa,\kappa]} |h_x(t,x)|$; and, for $\mh\in C^{0,1}(\WR,\R)$, we denote
$\n{\mh}_{1,\kappa}:=\max_{(\w,x)\in\W\times[-\kappa,\kappa]} |\mh(\w,x)|+
\max_{(\w,x)\in\W\times[-\kappa,\kappa]} |\mh_x(\w,x)|$.
\begin{teor}\label{th:2pers}
Let $h\in C^{0,1}(\RR,\R)$ be fixed, let $b_h$ be an attractive (resp. repulsive) hyperbolic solution
of \eqref{eq:2proceso} with dichotomy constant pair $(k_0,\gamma_0)$, and take $\kappa>\sup_{t\in\R}|b_h(t)|$.
Then, for every $\gamma\in(0,\gamma_0)$ and $\ep>0$, there exists $\delta_\ep>0$ and $\rho_\ep>0$ such that,
if $g$ is $C^1$-admissible and $\n{h-g}_{1,\kappa}<\delta_\ep$, then
there exists an attractive (resp.~repulsive) hyperbolic solution $b_g$ of $x'=g(t,x)$
with radius of uniform exponential stability dichotomy $\rho_\ep$ and dichotomy constant pair $(k_0,\gamma)$
such that $\sup_{t\in\R}|b_h(t)-b_g(t)|<\ep$.

Let $\mh\in C^{0,1}(\WR,\R)$ be fixed, let $\{\mb_\mh\}$ be an attractive (resp. repulsive) hyperbolic
copy of the base for \eqref{eq:2fam} with dichotomy constant pair $(\gamma_0,k_0)$, and take $\kappa>\max_{\w\in\W}|\mb_\mh(\w)|$.
Then, for every $\gamma\in(0,\gamma_0)$ and $\ep>0$, there exists $\delta_\ep>0$ and $\rho_\ep>0$ such that,
if $\mg\in C^{0,1}(\WR,\R)$ and $\|\mh-\mg\|_{1,\kappa}<\delta_\ep$, then
there exists an attractive (resp. repulsive) hyperbolic copy of the base $\{\mb_\mg\}$
of $x'=\mg(\wt,x)$ with radius of uniform exponential stability $\rho_\ep$ and dichotomy constant pair $(k_0,\gamma)$
such that $\max_{\w\in\W}|\mb_\mh(\w)-\mb_\mg(\w)|<\ep$.
\end{teor}
%%%%%%%%%%%%%%%%%%%%%%%%%%%%%%%%%%%%%%%%%%%%%%%%%%%%%%%%%%%%%%%%%%%%%%%%%%%%%
%%%%%%%%%%%%%%%%%%%%%%%%%%%%%%%%%%%%%%%%%%%%%%%%%%%%%%%%%%%%%%%%%%%%%%%%%%%%%
%%%%%%%%%%%%%%%%%%%%%%%%%%%%%%%%%%%%%%%%%%%%%%%%%%%%%%%%%%%%%%%%%%%%%%%%%%%%%
%%%%%%%%%%%%%%%%%%%%%%%%%%%%%%%%%%%%%%%%%%%%%%%%%%%%%%%%%%%%%%%%%%%%%%%%%%%%%
\section{A global bifurcation diagram in the concave-linear case}\label{sec:3}
Let $(\W,\sigma)$ be a continuous flow on a compact metric space.
In this section, we describe the bifurcation diagram for a family of concave-linear
ODEs of the type
\begin{equation}\label{eq:3conclin}
 x'=\mc(\wt)+\md(\wt)\,x+\mg^-(\wt,x)+\lb\,,\quad\w\in\W\,,
\end{equation}
where $\lb\in\R$, and derive analogous conclusions for linear-convex equations. Apart from their own interest,
these descriptions will be fundamental auxiliary tools for the study of bifurcations for
families of concave-convex equations in the next section. We will work under (all or part of) the
next hypotheses:
\begin{enumerate}[leftmargin=25pt,label=\rm{\bf{l\arabic*}}$^-$]
\item\label{l1} $\mc,\,\md\in C(\W,\R)$ and $\mg^-\in C^{0,1}(\WR,\R)$,
\item\label{l2} $\int_\W \md(\w)\, dm>0$ for all $m\in\merg$,
\item\label{l3} $\lim_{x\to\infty} (\md(\w)\,x+\mg^-(\w,x))=-\infty$ uniformly on $\W$,
\item\label{l4} $\mg^-(\w,x)\le 0$ for all $(\w,x)\in\W\times\R$, $\mg^-(\w,x)=0$ for all $(\w,x)\in\W\times(-\infty,0]$,
    and $m(\{\w\in\W\,|\;\mg^-(\w,x)<0$ for all $x>0\})=1$ for all $m\in\merg$,
\item\label{l5} $m(\{\w\in\W\,|\;x\mapsto \mg^-(\w,x)$ is concave$\})=1$ for all $m\in\merg$,
\item\label{l6} $m(\{\w\in\W\,|\;x\mapsto \mg^-_x(\w,x)$ is strictly decreasing on $\mJ\})>0$
    for all compact interval $\mJ\subset[0,\infty)$ and $m\in\merg$.
\end{enumerate}
In this section, we adapt the results obtained in \cite[Section 3]{dno5}
to this new setting. This requires to understand the relation between the set of
conditions \ref{l1}-\ref{l6} and the conditions \ref{c1}-\ref{c4} there assumed,
which we repeat for the reader's convenience. They refer to a map $\mh\colon\WR\to\R$:
\begin{enumerate}[leftmargin=20pt,label=\rm{\bf{c\arabic*}}]
\item\label{c1} $\mh\in C^{0,1}(\WR,\R)$,
\item\label{c2} $\limsup_{x\to\pm\infty}\mh(\w,x)<0$ uniformly on $\W$,
\item\label{c3} $m(\{\w\in\W\,|\;x\mapsto\mh(\w,x) \text{ is concave}\})=1$ for all $m\in\merg$,
\item\label{c4} $m(\{\w\in\W\,|\; x\mapsto\mh_x(\w,x)$ is strictly decreasing on $\mJ\})>0$
for all compact interval $\mJ\subset\R$ and all $m\in\merg$,
\end{enumerate}
For $\mh(\w,x):=\mc(\w)+\md(\w)\,x+\mg^-(\w,x)$, conditions \ref{l1} and \ref{c1} are equivalent.
Roughly speaking, conditions \ref{l2}, \ref{l3} and \ref{l4} will substitute the coercivity condition \ref{c2}:
it is clear that \ref{l3} is equivalent to $\limsup_{x\to\infty}(\mc(\w)+\md(\w)\,x+\mg^-(\w,x)+\lb)<0$
uniformly on $\W$ for all $\lb\in\R$; hypothesis \ref{l4} means that the dynamics of
\eqref{eq:3conclin} for $x\le 0$ coincides with that of the linear family
\begin{equation}\label{eq:3linearone}
 x'=\mc(\wt)+\md(\wt)\,x+\lb\,,\quad\w\in\W\,;
\end{equation}
and also that the right-hand side of \eqref{eq:3linearone} is greater than or equal to that
of \eqref{eq:3conclin}, which will be a fundamental property to apply comparison
arguments; and \ref{l2} determines the global dynamics of this linear family.
Conditions \ref{l5} and \ref{l6} are properties on concavity (in full measure)
of $x\mapsto\mc(\w)+\md(\w)\,x+\mg^-(\w,x)$ and strict concavity (in positive measure) of the same maps
on the compact subsets of $[0,\infty)$; and together with Remark \ref{rm:3con0}, they play the role of
conditions \ref{c3} and \ref{c4}.
\begin{nota}\label{rm:3con0}
Observe that \ref{l1}, \ref{l4} and \ref{l5} ensure that
$m(\{\w\in\W\,|\; \mg^-_x(\w,x)<0$ for all $x>0\})=1$ for every
$m\in\mathfrak{M}_\mathrm{erg}(\W,\sigma)$.
\end{nota}
We represent by \eqref{eq:3conclin}$^\lb$, \eqref{eq:3conclin}$_\w$, and \eqref{eq:3conclin}$^\lb_\w$
the families and equation obtained from \eqref{eq:3conclin} by fixing $\lb$, $\w$, and both of them, respectively;
and we proceed in the same way with other families of equations.
We call $\tau_\lb^-$ and $\tau_\lb^l$ the (local and global) skewproduct flows induced by
\eqref{eq:3conclin}$^\lb$ and \eqref{eq:3linearone}$^\lb$, with $\tau_\lb^-(t,\w,x)=
(\wt,v_\lb^-(t,\w,x))$ and $\tau_\lb^l(t,\w,x)=(\wt,v_\lb^l(t,\w,x))$; by $(\alpha^-_{\w,x,\lb},\beta^-_{\w,x,\lb})$
the maximal interval of definition of $v_\lb^-(t,\w,x)$;
by $\mB_\lb^-$ and $\mB_\lb^l$ the subsets of $\WR$ of globally bounded orbits for $\tau_\lb^-$ and $\tau_\lb^l$,
which are invariant for the respective flows; and by $(\W_\lb^b)^-$ the (possibly empty) projection of
$\mB_\lb^-$ over $\W$. In addition, whenever $(\W_\lb^b)^-$ is nonempty, we represent by
$\mr_\lb^-\colon(\W_\lb^b)^-\to\R$ and $\ma_\lb^-\colon(\W_\lb^b)^-\to\R$ the maps such that
\begin{equation}\label{def:3Bl-}
 \mB_\lb^-=\{(\w,x)\,|\;\w\in(\W_\lb^b)^-\text{ and }\mr_\lb^-(\w)\le x\le\ma_\lb^-(\w)\}\,,
\end{equation}
which are equilibria for the restriction of $\tau$ to $(\W_\lb^b)^-\times\R$, and are
respectively lower and upper semicontinuous (see Section \ref{subsec:2skew}).

Let us first analyze the properties of $\mB_\lb^l$.
\begin{prop}\label{prop:3CL-DE}
Assume that $\mc$ and $\md$ are continuous, and that \ref{l2} holds. Then,
\begin{itemize}[leftmargin=20pt]
\item[\rm(i)] the set $\mB^l_\lb$ is a repulsive hyperbolic $\tau_\lb^l$-copy of the base
for all $\lb\in\R$. In addition, $\mB^l_\lb=\{\mb_\lb\}$, with
$\lim_{t\to-\infty}|v^l_\lb(t,\w,x)-\mb_\lb(\wt)|=0$ for all $(\w,x)$, and
$\lim_{t\to\infty}v^l_\lb(t,\w,x)=\pm\infty$ if and only if $\pm x>\pm\mb_\lb(\w)$
\item[\rm(ii)] The map $\R\to C(\W,\R),\,\lb\mapsto\mb_\lb$ is continuous
and strictly decreasing, and $\lim_{\lb\to\pm\infty}\mb_\lb(\w)=\mp\infty$
uniformly on $\W$.
\end{itemize}
\end{prop}
\begin{proof}
(i) We include a proof of this classical result for the reader's convenience.
Condition \ref{l2} ensures that all the Lyapunov exponents of $\mB^l_\lb$ are strictly positive:
the lower one is $\int_\W\md(\w)\,dm^l$ for an $m^l\in\merg$ (see Section \ref{subsec:2skew}).
This fact ensures the exponential dichotomy of the family $x'=\md(\wt)\,x,\;\w\in\W$ (see,
e.g.,~\cite[Proposition 1.12(iv)]{duenasphdthesis}).
In turn, the global exponential dichotomy ensures the exponential dichotomy of
the equation $x'=\md(\wt)\,x$ over $\R$ for any $\w\in\W$ (see, e.g., \cite[Theorem 1.60]{johnson1}).
So,
\begin{equation}\label{def:3exprmu}
 b^\w_\lb(t):=-\int_t^\infty e^{\int_s^t\md(\w{\cdot}l)\,dl}(\mc(\ws)+\lb)\,ds
\end{equation}
is the unique bounded solution of \eqref{eq:3linearone} for each $\w\in\W$ and $\lb\in\R$
(see, e.g., \cite[Lecture 8]{coppel1}). Since $\w\mapsto \mb_\lb(\w):=b^\w_\lb(0)$
(which satisfies $\mb_\lb(\wt)=b^\w_\lb(t)$) is continuous on $\W$,
$\mB^l_\lb$ is the compact $\tau_\lb^l$-invariant set
$\{(\w,\mb_\lb(\w))\,|\;\w\in\W\}$, and hence Theorem \ref{th:2copia} ensures that
$\mB^l_\lb$ is a repulsive hyperbolic copy of the base. The last assertions in (i)
follow from $v_\lb^l(t,\w,x)=e^{\int_0^t\md(\ws)\,ds}(x-\mb_\lb(\w))+\mb_\lb(\wt)$,
since $\lim_{t\to\infty}e^{\int_0^t\md(\ws)\,ds}=\infty$ for all $\w\in\W$
(see, e.g., \cite[Proposition 1.56]{johnson1}).
\smallskip

(ii) Note that $\mb_\lb(\w)=\mb_0(\w)-\lb\,\mb(\w)$ for
$\mb(\w):=\int_0^\infty e^{\int_s^0\md(\w{\cdot}l)\,dl}\,ds$,
which is strictly positive and continuous on $\W$, and hence positively bounded from below.
All the assertions in (ii) follow easily from here.
\end{proof}
\begin{coro}\label{coro:3CLsolouna}
Assume that \ref{l1} and \ref{l2} hold. Then,
for any $\w\in\W$ and $\lb\in\R$, there exists at most a bounded solution
$b^\w_\lb\colon\R\to(-\infty,0]$ of \eqref{eq:3conclin}$_\w^\lb$, which exists if and only if
$\mb_\lb(\wt)\le 0$ for all $t\in\R$; and, in this case,
$b^\w_\lb(t)=\mb_\lb(\wt)$. In addition, any bounded solution $\bar x\colon\R\to\R$
of \eqref{eq:3conclin}$_\w^\lb$ satisfies $\bar x(t)\ge\mb_\lb(\wt)$ for all $t\in\R$.
\end{coro}
\begin{proof}
The first assertions follow immediately from Proposition \ref{prop:3CL-DE} and the relation between
\eqref{eq:3conclin} and \eqref{eq:3linearone}. To check the last one,
we assume for contradiction that $\bar x(s)<\mb_\lb(\ws)$. A classical comparison argument shows that
$\bar x(t)=v^-_\lb(t-s,\ws,\bar x(s))\le v_\lb^l(t-s,\ws,\bar x(s))<v_\lb^l(t-s,\ws,\mb_\lb(\ws))=\mb_\lb(\wt)$
for $t>s$, so that Proposition \ref{prop:3CL-DE}(i) ensures that $\bar x$ it is unbounded from below as time increases.
\end{proof}
The next result proves that the bifurcation diagram for \eqref{eq:3conclin} is
similar to that described in \cite[Theorem 3.3]{dno5} in the so-called ``concave in measure case".
Recall that $\mr_\lb^-$ and $\ma_\lb^-$ are defined by \eqref{def:3Bl-}.
\begin{teor}\label{th:3CLbifur}
If \ref{l1}-\ref{l6} hold, then there exist real values $\lb_*,\,\lb_-$ with $\lb_*\le\lb_-$
such that
\begin{itemize}[leftmargin=22pt]
\item[\rm(i)] $\lb>\lb_-$ if and only if $(\W_\lb^b)^-=\W$ and $(\ma_\lb^-,\mr_\lb^-)$ is an attractor-repeller
    pair of copies of the base for \eqref{eq:3conclin}$^\lb$. In particular, for $\lb>\lb_-$,
    $\lim_{t\to\infty}(v_\lb^-(t,\w,x)-\ma^-_\lb(\wt))=0$ if and only if $x>\mr_\lb^-(\w)$, and
    $\lim_{t\to-\infty}(v_\lb^-(t,\w,x)-\mr^-_\lb(\wt))=0$ if and only if $x<\ma_\lb^-(\w)$.
    Besides,
    the maps $(\lb_-,\infty)\to C(\W,\R),\lb\mapsto-\mr_\lb^-,\,\ma_\lb^-$ are continuous and strictly increasing,
    with $\lim_{\lb\to\infty}\mr_\lb^-=-\infty$ and $\lim_{\lb\to\infty}\ma_\lb^-=\infty$ uniformly on $\W$.
\item[\rm(ii)] $(\W^b_{\lb_-})^-=\W$; $\mr_{\lb_-}^-=\lim_{\lb\to(\lb_-)^+}\mr_\lb^-$
    and $\ma_{\lb_-}^-=\lim_{\lb\to(\lb_-)^+}\ma_\lb^-$ pointwise on $\W$;
    and $\inf_{\w\in\W}(\ma_{\lb_-}^-(\w)-\mr_{\lb_-}^-(\w))=0$.
\item[\rm (iii)] $\ma_{\lb_1}^->\ma_{\lb_-}^-\ge \mr_{\lb_-}^->\mr_{\lb_2}^-$ if $\lb_1>\lb^+$ and $\lb_2>\lb^+$.
\item[\rm(iv)] $\lb<\lb_*$ if and only if $(\W_{\lb}^b)^-$ is empty.
\item[\rm(v)] If $\lb_*<\lb_-$ and $\lb\in[\lb_*,\lb_-)$, then $(\W_{\lb}^b)^-$ is a proper subset of $\W$.
\end{itemize}
In addition,
\begin{itemize}[leftmargin=22pt]
\item[\rm(vi)] if $\W$ is minimal, then $\lb_*=\lb_-$.
\item[\rm(vii)] For each $\w\in\W$, let $\lb_-^\w$ be the upper special value of the parameter associated
    to the restriction of the $\lb$-parametric family \eqref{eq:3conclin} to the closure $\W_\w$ of
    $\{\wt\,|\;t\in\R\}$. Then, $\lb_*=\inf_{\w\in\W}\lb_-^\w$ and $\lb_-=\sup_{\w\in\W}\lb_-^\w$.
\end{itemize}
\end{teor}
\begin{proof}
\textsc{Step 1.} We fix $\lb\in\R$ and $m\in\merg$, and take two bounded and $m$-measurable $\tau_\lb^-$-equilibria
$\eta_1,\eta_2\colon\W\to\R$ with $m(\W_0)=1$, where $\W_0:=
\{\w\in\W\,|\;\eta_1(\w)<\eta_2(\w)\}$. Our goal is to check that
$\int_\W(\md(\w)+\mg^-_x(\w,\eta_1(\w)))\,dm>0>\int_\W(\md(\w)+\mg^-_x(\w,\eta_2(\w)))\,dm$.
According to \cite[Proposition 3.2]{dno4}, it suffices to check that
$m(\{\w\in\W\,|\;\mg^-_x(\w,\eta_1(\w))>\mg^-_x(\w,\eta_2(\w))\})>0$.
\par
If $m(\{\w\in\W_0\,|\;\eta_2(\w)>0\})=1$, then \ref{l6} proves the claim. So, we
assume $m(\{\w\in\W_0\,|\;\eta_2(\w)\le 0\})>0$.
Lusin's Theorem provides a compact
subset $\Delta\subseteq\{\w\in\W_0\,|\;\eta_2(\w)\le 0\}$ with $m(\Delta)>0$ such that
$\eta_1|_\Delta$ and $\eta_2|_\Delta$ are continuous. Let us check that, for any
$\w_0\in\Delta$, there exists $s_{\w_0}>0$ such that $\eta_2(\w_0{\cdot}s_{\w_0})>0$:
otherwise $t\mapsto \eta_1(\w_0{\cdot}t),\,\eta_2(\w_0{\cdot}t)$ define two different
nonpositive bounded solutions of \eqref{eq:3conclin}$_{\w_0}^\lb$, which contradicts
Corollary \ref{coro:3CLsolouna}. Hence, there exists
$s_{\w_0}>0$ such that $\eta_1(\w_0{\cdot}s_{\w_0})<0<\eta_2(\w_0{\cdot}s_{\w_0})$. Since
$\w\mapsto v_\lb^-(s_{\w_0},\w,\eta_i(\w))=\eta_i(\ws_{\w_0})$ is continuous on $\Delta$
for $i=1,2$, there exists $\rho_{\w_0}>0$ such that
$\eta_1(\ws_{\w_0})<0<\eta_2(\ws_{\w_0})$ for all $\w\in\Delta\cap B_\W(\w_0,\rho_{\w_0})$.
Obviously, $\Delta=\bigcup_{\w\in\Delta}(\Delta\cap B_\W(\w,\rho_\w))$. The compactness of $\Delta$
and $m(\Delta)>0$ ensure the existence of $\w_1\in\Delta$ such that $m(\Delta\cap B_\W(\w_1,\rho_{\w_1}))>0$, and
hence the $\sigma$-invariance of $m$ yields $m(\{\w\in\W_0\,|\;\eta_1(\w)<0<\eta_2(\w)\})>0$.
Since, for $\w$ in this set, $\mg_x^-(\w,\eta_1(\w))=\mg_x^-(\w,0)=0$, the assertion follows from
Remark \ref{rm:3con0}.
\smallskip

{\sc Step 2.} The property established in the previous step allows us to repeat the proof of
\cite[Theorem~3.3]{dno4} in order to prove that, for any $\lb\in\R$, there are at most two disjoint and
ordered $\tau_\lb^-$-invariant compact sets projecting onto $\W$, in which case
they are given by two hyperbolic copies of the base and agree with the lower and
upper bounded solutions: $\{\mr_\lb^-\}$, repulsive, and $\{\ma_\lb^-\}$, attractive.
\smallskip\par
{\sc Step 3.} We fix $\lb$ and take $m_2$ such that $\mc(\w)+\md(\w)\,x+\mg^-(\w,x)+\lb\le-\delta$
for a $\delta>0$ and all $\w\in\W$ if $x\ge m_2$. It is easy to check that
$\limsup_{t\to(\beta^-_{\w,x,\lb})^-}v_\lb^-(t,\w,x)\le m_2$ for all $(\w,x)\in\WR$,
and that any globally bounded solution is upper-bounded by $m_2$. We call
$m_1:=\inf_{\w\in\W}\mb_\lb(\w)$. According to Corollary \ref{coro:3CLsolouna},
any globally bounded solution is lower-bounded by $m_1$. In addition, for all $(\w,x)\in\WR$,
$\liminf_{t\to(\alpha^-_{\w,x,\lb})^+}v_\lb^-(t,\w,x)\ge m_1$: a comparison argument yields
$v^-_\lb(t,\w,x)-m_1\ge v^l_\lb(t,\w,x)-\mb_\lb(\wt)$ for all $t\le 0$, and
$\lim_{t\to-\infty}(v^l_\lb(t,\w,x)-\mb_\lb(\wt))=0$, as Proposition \ref{prop:3CL-DE}(i) ensures.
These properties reproduce \cite[Proposition~3.5(i)\&(ii)]{dno4} and allow us to prove the
remaining assertions of that proposition, which describes the properties of the set $\mB^-_\lb$,
its upper and lower bounds, and the global dynamics in some cases.
\smallskip\par
{\sc Step 4.} The previous properties allow us to repeat step by step the proofs of
\cite[Proposition 3.2 and Theorem 3.3]{dno5}.
\end{proof}
\begin{teor}\label{th:3CLbifuruna}
Assume all the hypotheses of Theorem {\rm\ref{th:3CLbifur}}.
Let us fix $\bar\w\in\W$, and let $r^-_\lb$ and $a^-_\lb$ represent the lower and upper
bounded solutions of $x'=\mc(\bwt)+\md(\bwt)\,x+\mg^-(\bwt,x)+\lb$ in the case of existence. Then,
\begin{itemize}[leftmargin=22pt]
\item[\rm(i)] there exists $\lb^{\bar\w}_-\in\R$ such that $x'=\mc(\bwt)+\md(\bwt)\,x+\mg^-(\bwt,x)+\lb$ has:
\begin{itemize}[leftmargin=10pt]
\item[-] two uniformly separated hyperbolic solutions
$r_\lb^-$ and $a_\lb^-$ for $\lb>\lb^{\bar\w}_-$, with
$r_\lb^-<a_\lb^-$, $r_\lb^-$ repulsive and
$a_\lb^-$ attractive;
\item[-] bounded but neither hyperbolic solutions nor uniformly separated solutions
for $\lb=\lb^{\bar\w}_-$;
\item[-] and no bounded solutions for $\lb<\lb^{\bar\w}_-$.
\end{itemize}
\item[\rm(ii)]The maps $(\lb^{\bar\w}_-,\infty)\to C(\R,\R)\,,\;\lb\mapsto-r^-_\lb,\,a^-_\lb$
are continuous for the uniform topology and strictly increasing
on $(\lb^{\bar\w}_-,\infty)$; $a^-_{\lb^{\bar\w}_-}=\lim_{\lb\to(\lb^{\bar\w}_-)^+}a^-_\lb$ and
$r^-_{\lb^{\bar\w}_-}=\lim_{\lb\to(\lb^{\bar\w}_-)^+}r^-_\lb$ pointwise on $\R$;
and $\lim_{\lb\to\infty}a^-_\lb=
-\lim_{\lb\to\infty}r^-_\lb=\infty$ uniformly on $\R$.
\item[\rm(iii)] The value $\lb^{\bar\w}_-$ is common for all $\w\in\W$ such that $\W_\w=\W_{\bar\w}$.
\item[\rm(iv)] If $\mu^{\bar\w}_-:=\inf\{\lb\in\R\,|\;\mb_\lb(\bwt)\le 0$ for all $t\in\R\}$, then
$\mu_-^{\bar\w}>\lb^{\bar\w}_-$ and $r^-_\lb(t)=\mb_\lb(\bwt)$ for all $t\in\R$ if $\lb\ge\mu_-^{\bar\w}$.
\end{itemize}
In addition, let $x_\lb(t,s,x)$ be the maximal solution of
$x'=\mc(\bwt)+\md(\bwt)\,x+\mg^-(\bwt,x)+\lb$ with $x_\lb(s,s,x)=x$,
defined on $(\alpha_{s,x,\lb},\beta_{s,x,\lb})$. Then,
\begin{itemize}[leftmargin=22pt]
\item[\rm(v)] if $\lb\ge\lb^{\bar\w}_-$, then
\begin{itemize}[leftmargin=10pt]
\item[-] $\lim_{t\to(\alpha_{s,x,\lb})^+}x_\lb(t,s,x)=\infty$
if and only if $x>a^-_\lb(s)$,
\item[-] and $\lim_{t\to(\beta_{s,x,\lb})^-}x_\lb(t,s,x)=-\infty$
if and only if $x<r^-_\lb(s)$;
\end{itemize}
if $\lb>\lb^{\bar\w}_-$, then
\begin{itemize}[leftmargin=10pt]
\item[-] $\lim_{t\to-\infty}|x_\lb(t,s,x)-r^-_\lb(t)|=0$
if and only if $x<a^-_\lb(s)$,
\item[-] $\lim_{t\to\infty}|x_\lb(t,s,x)-a^-_\lb(t)|=0$ if and only if
$x>r^-_\lb(s)$;
\end{itemize}
and, if $\lb<\lb^{\bar\w}_-$ and $x\in\R$, then
\begin{itemize}[leftmargin=10pt]
\item[-] $\lim_{t\to(\alpha_{s,x,\lb})^+}x_\lb(t,s,x)=\infty$ or
$\lim_{t\to(\beta_{s,x,\lb})^-}x_\lb(t,s,x)=-\infty$ (perhaps both).
\end{itemize}
\end{itemize}
\end{teor}
\begin{proof}
To prove (i), (ii), (iii) and (v), we can reason as in the proofs of \cite[Theorems 3.4 and 3.5]{dno5}
using the information of the proof of Theorem \ref{th:3CLbifur}. Note that the notation $\lb_-^{\bar\w}$
is coherent with that established in Theorem \ref{th:3CLbifur}(vii).
Let us prove (iv). Corollary \ref{coro:3CLsolouna}
shows that $r^-_\lb(t)\ge \mb_\lb(\bwt)$. Since $\mb_\lb(\bwt)$ solves \eqref{eq:3conclin}$_{\bar\w}^\lb$
for $\lb\ge\mu_-^{\bar\w}$, then $\lb_-^{\bar\w}\le\mu_-^{\bar\w}$, and Corollary \ref{coro:3CLsolouna} also
shows that $r^-_\lb(t)=\mb_\lb(\bwt)$ if $\lb\ge\mu_-^{\bar\w}$. In addition, if $\lb\ge\mu_-^{\bar\w}$, then
$t\mapsto\mb_\lb(\bwt)$ is a repulsive hyperbolic solution of \eqref{eq:3conclin}$_{\bar\w}^\lb$,
since $\mg_x^-(\bar\w,\mb_\lb(\bar\w))\equiv 0$. Since (i) ensures the lack of
hyperbolic solutions for $\lb^{\bar\w}_-$, we conclude that $\lb_-^{\bar\w}<\mu_-^{\bar\w}$.
\end{proof}
\begin{nota}\label{rm:3convex}
Let us consider now the flow $\tau_\lb^+$ induced by the family of equations
\begin{equation}\label{eq:3linconv}
 x'=\mc(\wt)+\md(\wt)\,x+\mg^+(\wt,x)+\lb\,,\quad\w\in\W
\end{equation}
for $\lb\in\R$, and the conditions
\begin{enumerate}[leftmargin=25pt,label=\rm{\bf{l\arabic*}}$^+$]
\item\label{l1+} $\mc,\,\md\in C(\W,\R)$ and $\mg^+\in C^{0,1}(\WR,\R)$,
\item\label{l2+} $\int_\W \md(\w)\, dm>0$ for all $m\in\merg$,
\item\label{l3+} $\lim_{x\to-\infty} (\md(\w)\,x+\mg^+(\w,x))=\infty$ uniformly on $\W$,
\item\label{l4+} $\mg^+(\w,x)\ge 0$ for all $(\w,x)\in\W\times\R$, $\mg^+(\w,x)=0$ for all $(\w,x)\in\W\times[0,\infty)$,
    and $m(\{\w\in\W\,|\;\mg^+(\w,x)>0$ for all $x<0\})=1$ for all $m\in\merg$,
\item\label{l5+} $m(\{\w\in\W\,|\;x\mapsto \mg^+(\w,x)$ is convex$\})=1$ for all $m\in\merg$,
\item\label{l6+} $m(\{\w\in\W\,|\;x\mapsto \mg^+_x(\w,x)$ is strictly increasing on $\mJ\})>0$
    for all compact interval $\mJ\subset(-\infty,0]$ and $m\in\merg$.
\end{enumerate}
Then, Corollary \ref{coro:3CLsolouna} and Theorems \ref{th:3CLbifur} and \ref{th:3CLbifuruna}
have ``symmetric" formulations. In particular, there exist $\lb^+\le\lb^*$ such that:
\begin{itemize}[leftmargin=15pt]
\item[-] $\lb<\lb^+$ if and only if $\tau_\lb^+$ has an attractor-repeller pair of copies of the base
$(\ma^+_\lb,\mr^+_\lb)$, in which case $\ma^+_\lb$ is attractive, $\mr^+_\lb$ is repulsive, and $\mr_\lb^+>\ma_\lb^+$.
\item[-] $\lb>\lb^*$ if and only if none of the equations \eqref{eq:3linconv}$_\w^\lb$
has a bounded solution.
\item[-] If $\lb^*>\lb^+$ and $\lb\in(\lb^+,\lb^*]$,
then $(\W_{\lb}^b)^+$ (the projection onto $\W$ of the set of globally bounded orbits for \eqref{eq:3linconv}$^\lb$) is a proper subset of $\W$.
\end{itemize}
In addition, for a fixed $\bar\w\in\W$,
\begin{itemize}[leftmargin=15pt]
\item[-] any bounded solution of \eqref{eq:3linconv}$_{\bar\w}^\lb$ satisfies $\bar x(t)\le \mb_\lb(\bwt)$.
\item[-] If $\mu^+_{\bar\w}:=\sup\{\lb\in\R\,|\;\mb_\lb(\bwt)\ge 0$ for all $t\in\R\}$ and
    $\lb\le\mu^+_{\bar\w}$, then $\mb_\lb(\bwt)$ is the upper bounded solution of \eqref{eq:3linconv}$_{\bar\w}^\lb$,
    and it is hyperbolic repulsive.
\item[-] The absence of bounded solutions for \eqref{eq:3linconv}$_{\bar\w}^\lb$
    holds on a positive halfline $(\lb^+_{\bar\w},\infty)$,
    and the existence of an attractor repeller-pair of solutions on a negative halfline
    $(-\infty,\lb^+_{\bar\w})$, with $\lb^+_{\bar\w}>\mu^+_{\bar\w}$; and, if such a pair exists,
    the upper hyperbolic solution $r^+_\lb$ (resp.~the lower hyperbolic solution $a^+_\lb$) is repulsive (resp.~attractive),
    and they determine the global dynamics. In addition,  $a^+_\lb$ and $-r^+_\lb$ are strictly increasing on
    $(-\infty,\lb^+_{\bar\w})$, with (uniform on $\R$) limit $-\infty$ as $\lb\to-\infty$.
\item[-] \eqref{eq:3linconv}$_{\bar\w}^{\lb^+_{\bar\w}}$ has no uniformly
    separated bounded solutions.
\end{itemize}
These properties are just a summary of all the information provided by the linear-convex analogues of
Theorems \ref{th:3CLbifur} and \ref{th:3CLbifuruna}.
The simplest way to prove these analogues is to observe that the change of
variables $y(t)=-x(t)$ takes \eqref{eq:3linconv} to
\begin{equation*}
 y'=-\mc(\wt)+\md(\wt)\,y-\mg^+(\wt,-y)-\lb\,,\quad\w\in\W\,,
\end{equation*}
and that $\md$ and the map $\mg^-$ defined as $\mg^-(\w,x):=-\mg^+(\w,-x)$ satisfy all the
conditions \ref{l1}-\ref{l6}; and it is easy to reinterpret the conclusions of
Theorems \ref{th:3CLbifur} and \ref{th:3CLbifuruna} in the light of these changes.
\end{nota}
\begin{nota}\label{rm:3localSNB}
Assume that we are dealing with a parametric family of flows $\tau_\lb$ on $\WR$, and
that, for a value $\lb_0$ of the parameter, there exists a nonhyperbolic compact $\tau_{\lb_0}$-invariant set
projecting onto $\W$ whose upper and lower equilibria are not uniformly separated. If, for $\lb$ in a neighborhood of
$\lb_0$, the local dynamics is that of the saddle-node bifurcation of
\cite[Theorem 3.3]{dno5}, similar to that of Theorem \ref{th:3CLbifur},
then we will say that the parametric family {\em exhibits a local saddle-node bifurcation as $\lb\downarrow\lb_0$}.
And if the local dynamics is that of the saddle-node bifurcation of the analogous convex results
(as those described in Remark \ref{rm:3convex}), then we
we will say that it {\em exhibits a local saddle-node bifurcation as $\lb\uparrow\lb_0$}.
\end{nota}
%%%%%%%%%%%%%%%%%%%%%%%%%%%%%%%%%%%%%%%%%%%%%%%%%%%%%%%%%%%%%%%%%%
%%%%%%%%%%%%%%%%%%%%%%%%%%%%%%%%%%%%%%%%%%%%%%%%%%%%%%%%%%%%%%%%%%
\section{Concave-convex bifurcations}\label{sec:4}
As in the previous section, $(\W,\sigma)$ is a continuous flow on a compact metric space.
In this section, we describe bifurcation diagrams for a family of concave-convex
equations of the type
\begin{equation}\label{eq:4conconv}
 x'=\mc(\wt)+\md(\wt)\,x+\mg(\wt,x)+\lb\,,\quad\w\in\W\,,
\end{equation}
under the next hypotheses:
\begin{enumerate}[leftmargin=23pt,label=\rm{\bf{cc\arabic*}}]
\item\label{cc1} $\mc,\,\md\in C(\W,\R)$ and $\mg\in C^{0,1}(\WR,\R)$,
\item\label{cc2} $\int_\W \md(\w)\, dm>0$ for all $m\in\merg$,
\item\label{cc3} $\lim_{x\to\pm\infty} (\md(\w)\,x+\mg(\w,x))=\mp\infty$ uniformly on $\W$,
\item\label{cc4} $\mg(\w,0)=\mg_x(\w,0)=0$ for all $\w\in\W$, $\mg(\w,x)\le 0$ for all
    $(\w,x)\in\W\times[0,\infty)$, $\mg(\w,x)\ge 0$ for all $(\w,x)\in\W\times(-\infty,0]$, and
    $m(\{\w\in\W\,|\;\mg(\w,x)<0$ for all $x>0$ and $\mg(\w,x)>0$ for all $x<0\})=1$ for all $m\in\merg$,
\item\label{cc5} $m(\{\w\in\W\,|\;x\mapsto \mg(\w,x)$ is concave on $[0,\infty)$ and convex on $(-\infty,0]\})=1$ for all $m\in\merg$,
\item\label{cc6} $m(\{\w\in\W\,|\;x\mapsto \mg_x(\w,x)$ is strictly decreasing on $\mJ$ and strictly increasing on $-\mJ\})>0$
    for all compact interval $\mJ\subset[0,\infty)$ and $m\in\merg$.
\end{enumerate}
We denote by $\tau_\lb$ the (possibly local) skewproduct flow induced by
\eqref{eq:4conconv}$^\lb$ on $\WR$, with $\tau_\lb(t,\w,x)=
(\wt,v_\lb(t,\w,x))$.
\begin{nota}\label{rm:4existeA}
According to \cite[Proposition~5.5]{dno4},
conditions \ref{cc1} and \ref{cc3} ensure the existence of the global attractor $\mA_\lb$
of $\tau_\lb$ for each $\lb\in\R$, which is compact and composed by the graphs of the
globally bounded $\tau_\lb$-orbits, and which we write as
\[
 \mA_\lb=\bigcup_{\w\in\W}\big(\{\w\}\times[\ml_\lb(\w),\muk_\lb(\w)]\big)\,.
\]
The maps $\ml_\lb\colon\W\to\R$ and $\muk_\lb\colon\W\to\R$ are lower and upper semicontinuous,
respectively. The comparison results of \cite[Proposition~5.5(iv)]{dno4}
are also valid in this setting, and will be often used. In particular, they
allow us to reason as in the proof of \cite[Theorem 5.5(i)]{dno1} in order to prove
that, for all $\w\in\W$, $\lb\mapsto\ml_\lb(\w)$ is
strictly increasing and left-continuous, and $\lb\mapsto\muk_\lb(\w)$ is
strictly increasing and right-continuous on $\R$ for all $\w\in\W$.
\end{nota}
Let us define the maps
\[
 \mg^-(\w,x):=
 \left\{\begin{array}{ll}
 0 &\text{for }x<0\,,\\
 \mg(\w,x)&\text{for }x\ge 0\,,
 \end{array}\right.
 \qquad
 \mg^+(\w,x):=
 \left\{\begin{array}{ll}
 \mg(\w,x) &\text{for }x<0\,,\\
 0&\text{for }x\ge 0
 \end{array}\right.
\]
(i.e., $\mg^-:=\min(0,\mg)$ and $\mg^+:=\max(0,\mg)$), and consider the families
\begin{equation}\label{eq:4conclin}
 x'=\mc(\wt)+\md(\wt)\,x+\mg^-(\wt,x)+\lb\,,\quad\w\in\W\,,
\end{equation}
\begin{equation}\label{eq:4linconv}
 x'=\mc(\wt)+\md(\wt)\,x+\mg^+(\wt,x)+\lb\,,\quad\w\in\W\,,
\end{equation}
which respectively satisfy the concave-linear hypotheses
\ref{l1}-\ref{l6} and the linear-convex hypotheses \ref{l1+}-\ref{l6+}. For these
equations, we will use the
notation established in Section \ref{sec:3} for \eqref{eq:3conclin} and \eqref{eq:3linconv}.
The corresponding bifurcation diagrams, and the dynamic for each fixed value of $\lb$,
will be fundamental tools to study the concave-convex dynamics.
An example is the next result, which we will often use.
\begin{lema}\label{lema:4copy}
Assume that \ref{cc1}-\ref{cc6} hold. Then,
\begin{itemize}[leftmargin=22pt]
\item[\rm(i)] a map $x\colon(\alpha,\beta)\to[0,\infty)$ (resp.~$x\colon(\alpha,\beta)\to(-\infty,0]$)
is a solution of \eqref{eq:4conconv}$^\lb_\w$ if and only if it is a solution of \eqref{eq:4conclin}$^\lb_\w$
(resp.~of \eqref{eq:4linconv}$^\lb_\w$).
\item[\rm(ii)] Let $\mb\colon\W\to\R$ be a continuous map, and assume that $\mb\ge 0$ (resp~$\mb\le 0$). Then,
$\{\mb_\lb\}$ is an attractive or repulsive hyperbolic $\tau_\lb$-copy of the base if and
only if it is an attractive or repulsive hyperbolic $\tau_\lb^-$-copy of the base
(resp.~$\tau_\lb^+$-copy of the base).
\end{itemize}
\end{lema}
\begin{proof}
We reason for nonnegative maps. Note that $\mg=\mg^-$ and $\mg_x=\mg^-_x$ on $\W\times[0,\infty)$. This proves (i) and that, if $\mb\ge 0$, then $\{\mb\}$ is a $\tau_\lb$ copy of the base if and only if it is a $\tau_\lb^-$-copy of the base. In addition, in this case, $\mg_x(\w,\mb_\lb(\w))=\mg^-_x(\w,\mb_\lb(\w))$,
and hence the upper and lower Lyapunov exponents of $\{\mb_\lb\}$ are common for $\tau_\lb$ and $\tau^-_\lb$, which ensures that the hyperbolicity and its type are also common (see Theorem \ref{th:2copia}).
\end{proof}
Let us recover the significative values of the parameter for the auxiliary families.
The values $\lb_*\le\lb_-$ are associated to \eqref{eq:4conclin} by Theorem \ref{th:3CLbifur}.
Let $\lb^+\le\lb^*$ be the analogous
values for the family \eqref{eq:4linconv} (see Remark \ref{rm:3convex}).
In what follows, $(\ma_\lb^-,\mr_\lb^-)$ (with $\mr_\lb^-<\ma_\lb^-$) represents the attractor-repeller
pair of copies of the base of \eqref{eq:4conclin} for $\lb>\lb_-$, and $(\ma_\lb^+,\mr_\lb^+)$
(with $\mr_\lb^+>\ma_\lb^+$) is that of \eqref{eq:4linconv} for $\lb<\lb^+$.
Their respective limits $\mr_{\lb^-}^-\le\ma_{\lb-}^-$ as $\lb\downarrow\lb_-$ and
$\ma_{\lb^+}^+\le\mr_{\lb^+}^+$ as $\lb\uparrow\lb^+$ are equilibria for the corresponding flows and,
as in the case of the hyperbolic copies of the base, they bound the respective sets of bounded solutions.
Recall also that $\ma_\lb^-$ and $-\mr_\lb^-$ (resp.~$\ma_\lb^+$ and $-\mr_\lb^+$) are strictly increasing
on $[\lb_-,\infty)$ (resp.~$(-\infty,\lb^+]$): see Theorem \ref{th:3CLbifur}  and Remark \ref{rm:3convex}.
Our first result establishes some general properties under conditions \ref{cc1}-\ref{cc6}.
\begin{prop}\label{prop:4tapas}
Assume that \ref{cc1}-\ref{cc6} hold. Then,
\begin{itemize}[leftmargin=22pt]
\item[\rm(i)] $\{\ml_\lb\}$ is an attractive hyperbolic $\tau_\lb$-copy of the base
if $\lb<\lb^+$, with $\ml_{\lb}\le\ma^+_{\lb}$, and the (strictly increasing)
map $(-\infty,\lb^+)\to C(\W,\R),\,\lb\mapsto\ml_\lb$ is continuous.
In addition, $\ml_{\lb}=\ma^+_{\lb}$ if $\ma^+_{\lb}\le 0$ (which happens if $-\lb$ is large enough).
\item[\rm(ii)]
$\{\muk_\lb\}$ is an attractive hyperbolic $\tau_\lb$-copy of the base
if $\lb>\lb_-$, with $\muk_{\lb}\ge\ma^-_{\lb}$, and the (strictly increasing)
map $(\lb_-,\infty)\to C(\W,\R),\,\lb\mapsto\muk_\lb$ is continuous.
In addition, $\muk_{\lb}=\ma^-_{\lb}$ if $\ma^-_{\lb}\ge 0$ (which happens if $\lb$ is large enough).
\item[\rm(iii)] If $\lb_-<\lb_1<\lb_2$ (resp.~$\lb_1<\lb_2<\lb^+$), then
$\mb_{\lb_2}<\mb_{\lb_1}\le\mr^-_{\lb_1}<\ma^-_{\lb_1}\le\muk_{\lb_1}<\muk_{\lb_2}$
(resp.~$\ml_{\lb_1}<\ml_{\lb_2}\le\ma^+_{\lb_2}<\mr^+_{\lb_2}\le\mb_{\lb_2}<
\mb_{\lb_1}$). In addition, $\mb_{\lb_1}=\mr^-_{\lb_1}$ if $\mb_{\lb_1}\ge 0$
(which happens if $-\lb_1$ is large enough), and $\mb_{\lb_2}=\mr^+_{\lb_2}$
if $\mb_{\lb_2}\le 0$ (which happens if $\lb_2$ is large enough).
\end{itemize}
\end{prop}
\begin{proof}
(i) We fix $\lb<\lb^+$. The map $\ma^+_\lb$ is a bounded
continuous global upper solution for $\tau_\lb$, since $(\ma^+_\lb)'(\w)\ge \mc(\w)+\md(\w)\,\ma^+_\lb(\w)+
\mg(\wt,\ma^+_\lb(\w))+\lb$. According to \cite[Proposition 5.5(iv)]{dno4}, $\ma^+_\lb\ge\ml_\lb$.
Hence, $\mb_\lb^-(\w):=\lim_{s\to\infty}v_\lb(s,\w{\cdot}(-s), \ma_\lb^+(\w{\cdot}(-s)))$,
which satisfies $\ma_\lb^+\ge\mb_\lb^-\ge\ml_\lb$, is an upper semicontinuous $\tau_\lb$-equilibrium
(see, e.g., \cite[Section 2.2]{dno1} and references therein).
So, $\mK:=\{(\w,x)\mid\,\ml_\lb(\w)\le x\le \mb_\lb^-(\w)\}$ is
a compact $\tau_\lb$-invariant set. Its upper Lyapunov exponent is
$\gamma_\mK^+=\int_\W(\md(\w)+\mg_x(\w,\mb(\w)))\,dm$, where $m\in\merg$ and $\mb\colon\W\to\R$ is an
$m$-measurable $\tau_\lb$-equilibrium with graph in $\mK$. Note that
$\int_\W\mg_x(\w,\mb(\w))\,dm=\int_\W\mg^+_x(\w,\mb(\w))\,dm+\int_\W\mg^-_x(\w,\mb(\w))\,dm
\le \int_\W\mg^+_x(\w,\mb(\w))\,dm\le \int_\W\mg^+_x(\w,\ma^+_\lb(\w))\,dm$:
the first inequality is due to $m(\{\w\,|\;\mg_x^-(\w,{\cdot})\le 0\})=1$ (which follows from
\ref{cc4} and \ref{cc5}), and the second one is due the nondecreasing character in full measure of
$x\mapsto\mg^+_x(\w,x)$ (a new consequence of \ref{cc5}). So,
$\gamma_\mK^+\le \int_\W(\md(\w)+\mg^+_x(\w,\ma^+_\lb(\w)))\,dm$.
The attractive hyperbolicity of $\ma^+_\lb$ for $\tau_\lb^+$ (see Remark \ref{rm:3convex})
ensures that the last integral is strictly negative (see Theorem \ref{th:2copia}), and
hence $\gamma_\mK^+<0$. If $\mM\subseteq\W$ is a minimal subset, then
$\mK^\mM:=\{(\w,x)\in\mK\,|\;\w\in\mM\}$ is a compact invariant set for the
restriction of $\tau_\lb$ to $\mM\times\R$.
Since all its Lyapunov exponents (also Lyapunov exponents for $\mK$)
are negative, \cite[Theorem 3.4]{cano} ensures that $\mK^\mM$ is a $\tau_\lb$-copy of the
base. In these conditions, Theorem \ref{th:2copia} ensures that $\mK$ itself is an attractive
hyperbolic $\tau_\lb$-copy of the base. It follows easily that $\mK=\{\mb_\lb^-\}=\{\ml_\lb\}$,
which proves the attractive hyperbolicity of $\ml_\lb$.

To prove the continuity of $(-\infty,\lb^+)\to C(\W,\R),\,\lb\mapsto\ml_\lb$,
we use the persistence ensured by Theorem \ref{th:2pers} combined with the
monotonicity and left-continuity established in Remark \ref{rm:4existeA} and with
fact that $v_\lb(t,\w,x)$ is unbounded if $x<\ml_\lb(\w)$.

It remains the last assertion in (i). Recall that $\ma_\lb^+\le 0$ if $-\lb$ is large enough:
see Remark \ref{rm:3convex}. As seen above, $\ma^+_\lb\ge\ml_\lb$. Since
Lemma \ref{lema:4copy} shows that
any $\tau_\lb$-orbit strictly below $\{\ma^+_\lb\}$ is a $\tau^+_\lb$-orbit,
and hence it is unbounded, we conclude that $\ma^+_\lb=\ml_\lb$, as asserted.
\smallskip

(ii) The proof is similar. Now, $\mb_\lb^+(\w):=\lim_{s\to\infty}v_\lb(s,\w{\cdot}(-s),
\ma_\lb^-(\w{\cdot}(-s)))\ge\ma_\lb^-(\w)$ is a lower semicontinuous $\tau_\lb$-equilibrium below
$\muk_\lb$ if $\lb>\lb_-$, and we can repeat the previous argument for
$\mK:=\{(\w,x)\mid\,\mb_\lb^+(\w)\le x\le\muk_\lb(\w)\}$.
\smallskip

(iii) We work for $\lb_-<\lb_1<\lb_2$. The first inequality is proved by Proposition
\ref{prop:3CL-DE}(ii). Corollary \ref{coro:3CLsolouna} shows the second one.
The third inequality follows from Theorem \ref{th:3CLbifur}(i); the
forth one follows from (ii); and the last one
is explained in Remark \ref{rm:4existeA}. The case $\lb_1<\lb_2<\lb^+$ is similar.
The last assertions can be deduced from Corollary \ref{coro:3CLsolouna} and its
linear-convex analogue, as in the proof of Theorem \ref{th:3CLbifuruna}(iv).
\end{proof}

With the aim of partly describe the bifurcation diagram in a particulary interesting case, we define
\begin{equation}\label{def:4mu}
\begin{split}
 &\mu_-:=\inf\{\lb\in\R\,|\;\mb_\lb(\w)\le 0 \text{ for all $\w\in\W$}\}\,,\\
 &\mu^+:=\sup\{\lb\in\R\,|\;\mb_\lb(\w)\ge 0 \text{ for all $\w\in\W$}\}\,,
\end{split}
\end{equation}
and observe that $\mu_-=\sup_{\bar\w\in\W}\mu_-^{\bar\w}$ and $\mu^+=\inf_{\bar\w\in\W}\mu^+_{\bar\w}$,
with $\mu_-^{\bar\w}$ and $\mu^+_{\bar\w}$ defined in Theorem \ref{th:3CLbifuruna}(iv) and
Remark \ref{rm:3convex}. Definition \ref{def:4mu}
combined with \eqref{def:3exprmu} and with $\mb_\lb(\w)=b_\lb^\w(0)$ provides
\begin{equation}\label{eq:4mu-mu+}
\begin{split}
\mu_-&=-\inf_{\w\in\W}\frac{\int_0^\infty \exp\left(\int_s^0\md(\w{\cdot}l)\,dl\right)\,\mc(\ws)\, ds}{\int_0^\infty \exp\left(\int_s^0\md(\w{\cdot}l)\,dl\right)\, ds}\,,\\
\mu^+&=-\sup_{\w\in\W}\frac{\int_0^\infty \exp\left(\int_s^0\md(\w{\cdot}l)\,dl\right)\,\mc(\ws)\, ds}{\int_0^\infty \exp\left(\int_s^0\md(\w{\cdot}l)\,dl\right)\, ds}\,.\\
\end{split}
\end{equation}

Theorem \ref{th:4CCbifur} explains the evolution of the global dynamics in the case
$\lb_*\le\lb_-<\mu^+\le\mu_-<\lb^+\le\lb^*$. As we will see in Section \ref{subsec:4todos},
this particular order is not the unique possible one, but it always holds in the autonomous case.
The information of Remark \ref{rm:4existeA} and Proposition \ref{prop:4tapas}
contributes to the global depiction of this evolution.
The respulsivity of an invariant compact set also comes into play to understand it.
Given a compact $\tau_\lb$-invariant set $\mK\subset\WR$ projecting onto $\W$ and such that
$\mK_\w$ is an interval for all $\w$, we say that it is {\em repulsive\/}
if there exists a bounded open set $\mU$ containing $\mK$ such that
every $\tau_\lb$-orbit starting at $\mU\setminus\mK$ leaves $\mU$ as time increases.
\begin{teor} \label{th:4CCbifur}
Assume that \ref{cc1}-\ref{cc6} hold, and that $\lb_*\le\lb_-<\mu^+\le\mu_-<\lb^+\le\lb^*$ hold.
Then,
\begin{itemize}[leftmargin=25pt]
\item[\rm(i)] $\{\ma^-_\lb\}$ (resp. $\{\ma^+_\lb\}$) is an attractive hyperbolic
$\tau_\lb$-copy of the base for all $\lb>\lb_-$ (resp. $\lb<\lb^+$),
and $\{\mr^-_\lb\}$ (resp. $\{\mr^+_\lb\}$) is a repulsive hyperbolic
$\tau_\lb$-copy of the base for all $\lb\in(\lb_-,\mu^+]$ (resp. $\lb\in[\mu_-,\lb^+)$).
\item[\rm(ii)] $\ml_\lb=\ma^+_\lb$ for all $\lb\le\lb^+$ and $\muk_\lb=\ma^-_\lb$ for all $\lb\ge\lb_-$. In particular,
    $\mA_\lb=\bigcup_{\w\in\W}(\{\w\}\times[\ma_\lb^+(\w),\ma_\lb^-(\w)])$ for all $\lb\in[\lb_-,\lb^+]$.
\item[\rm(iii)] If the $\sigma$-orbit of a point $\bar\w\in\W$ is dense, then $\lim_{t\to\infty}(\muk_\lb(\bwt)-\ma^+_\lb(\bwt))=0$ for
    $\lb\in[\lb_*,\lb_-)$, and $\lim_{t\to\infty}(\ma^-_\lb(\bwt)-\ml_\lb(\bwt))=0$ for
    $\lb\in(\lb^+,\lb^*]$.
\item[\rm(iv)] $\mA_\lb$ is an attractive hyperbolic $\tau_\lb$-copy of the base if $\lb\not\in[\lb_*,\lb^*]$:
    $\mA_\lb=\{\ma_\lb^+\}$ for $\lb<\lb_*$ and $\mA_\lb=\{\ma_\lb^-\}$ for $\lb>\lb^*$.
\end{itemize}
Regarding the internal dynamics of the global attractor,
\begin{itemize}[leftmargin=25pt]
\item[\rm(v)] for every $\lb\in(\lb_-,\lb^+]$, there exists an upper semicontinuous
    $\tau_\lb$-equilibrium $\mm_\lb^-\colon\W\to\R$ such that
    $\lim_{t\to\infty}(v_\lb(t,\w,x)-\muk_\lb(\wt))=0$ if and only if $x>\mm_\lb^-(\w)$. In addition,
    $\lb\mapsto\mm_\lb^-(\w)$ is
    left-continuous and strictly decreasing on $(\lb_-,\lb^+]$ for all $\w\in\W$,
    $\mm_\lb^-=\mr_\lb^-$ whenever $\mr_\lb^-\ge 0$, and $\mm_\lb^-=\mr_\lb^+$ for all
    $\lb\in[\mu_-,\lb^+]$.
\item[\rm(vi)] for every $\lb\in[\lb_-,\lb^+)$, there exists a lower semicontinuous
    $\tau_\lb$-equilibrium $\mm_\lb^+\colon\W\to\R$ such that
%    $\lim_{t\to\infty}(v_\lb(t,\w,x)-\ma_\lb^+(\wt))=0$
    $\lim_{t\to\infty}(v_\lb(t,\w,x)-\ml_\lb(\wt))=0$ if and only if $x<\mm_\lb^+(\w)$. In addition,
    $\lb\mapsto\mm_\lb^+(\w)$ is
    right-continuous and strictly decreasing on $[\lb_-,\lb^+)$ for all $\w\in\W$,
    $\mm_\lb^+=\mr_\lb^+$ whenever $\mr_\lb^+\le 0$, and $\mm_\lb^+=\mr_\lb^-$ for all
    $\lb\in[\lb_-,\mu^+]$.
\item[\rm(vii)] $\mm_\lb^-=\mm_\lb^+=\mr_\lb^-$ for all $\lb\in(\lb_-,\mu^+)$,
    $\mm_\lb^-=\mm_\lb^+=\mr_\lb^+$ for $\lb\in(\mu_-,\lb^+)$,
    and the set $\{(\w,x)\mid\, \mm_\lb^+(\w)\le x\le\mm_\lb^-(\w)\}$ is a repulsive compact
    $\tau_\lb$-invariant subset of $\WR$ for all $\lb\in(\lb_-,\lb^+)$.
\item[\rm(viii)] For any $\lb\in(\lb_-,\lb^+)$ and every $m\in\merg$,
    \begin{equation*}
    \int_\W\mh_x(\w,\mm_\lb^-(\w))\, dm\ge 0\quad\text{and}\quad\int_\W\mh_x(\w,\mm_\lb^+(\w))\, dm\ge 0\,.
    \end{equation*}
\end{itemize}
\end{teor}
\begin{proof} We take $\lb_1>\lb_-$ and $\lb_2\in(\lb_-,\mu^+]$. Then, Theorem \ref{th:3CLbifur}(iii),
Corollary~\ref{coro:3CLsolouna}, and the definition \eqref{def:4mu} of $\mu^+$ respectively show the
three inequalities $\ma_{\lb_1}^->\mr_{\lb_2}^-\ge\mb_{\mu^+}\ge 0$, and
Lemma \ref{lema:4copy}(ii) proves the assertions about $\{\ma^-_{\lb_1}\}$ and $\{\mr^-_{\lb_2}\}$.
Those concerning $\{\ma_\lb^+\}$ on $(-\infty,\lb^+)$ and $\{\mr^+_\lb\}$ on $[\mu_-,\lb^+)$
are proved using the linear-convex analogues of Theorem \ref{th:3CLbifur}(iii) and
Corollary~\ref{coro:3CLsolouna}.
\smallskip
(ii) As before, we have $\ma_\lb^->\mr_{\mu_+}^-\ge 0$ if $\lb\ge\lb_-$. So, Proposition \ref{prop:4tapas}(ii) shows that
$\muk_\lb=\ma_\lb^-$. The second assertion in (ii) is similar, and both imply the last one.
\smallskip

(iv) We reason for $\lb<\lb_*$. By (ii), it suffices to check that $t\mapsto v_\lb(t,\w,x)$
is unbounded (as time decreases) if $x>\ma_\lb^+(\w)$. For contradiction, we
assume that $t\mapsto v_\lb(t,\bar\w,\bar x)$ is bounded for a point $(\bar\w,\bar x)$ with $\bar x>\ma_\lb^+(\bar\w)$.
The linear-convex version of Theorem~\ref{th:3CLbifur}(i) yields
$\lim_{t\to-\infty}(v_\lb^+(t,\bar\w,\bar x)-\mr_\lb^+(\bwt))=0$. Hence, since $\mr_\lb^+>\mr_{\mu^+}^+$, there exists
$s<0$ such that $v^+_\lb(t,\bar\w,\bar x)>\mr_{\mu_+}^+(\bwt)\ge 0$  for all $t\le s$, with the last inequality
proved by the linear-convex version
of Corollary \ref{coro:3CLsolouna}. Since $\mg^+\ge\mg$, a standard comparison
argument shows that $v_\lb(t,\bar\w,\bar x)\ge v_\lb^+(t,\bar\w,\bar x)$ for all $t<0$,
so $v_\lb(t,\bar\w,\bar x)>0$ for all $t\le s$. Lemma \ref{lema:4copy} ensures that $t\mapsto v_\lb(t,\bar\w,\bar x)$
solves \eqref{eq:4conclin}$^\lb_{\bar\w}$ in $(-\infty,s]$, and hence $v_\lb(t,\bar\w,\bar x)=v_\lb^-(t-s,\bws,v_\lb(s,\bar\w,\bar x))$
for $t\le s$. Finally, we take a sequence $(t_n)\downarrow-\infty$ for which $(\w_0,x_0):=
\lim_{n\to\infty}(\bwt_n,v^-_\lb(t_n-s,\bws,v_\lb(s,\bar\w,\bar x)))$
exists, and check that $v_\lb^-(t,\w_0,x_0)$ is bounded, which impossible for $\lb<\lb_*$: see
Theorem \ref{th:3CLbifur}(iv).
\smallskip

(v)-(vii)
For $\lb\in(\lb_-,\lb^+]$, we define $\mm_\lb^-(\w):=\inf\{x\in\R\,|\;\lim_{t\to\infty}(v_\lb(t,\w,x)-\ma_\lb^-(\wt))=0\}=
\inf\{x\in\R\,|\;\lim_{t\to\infty}(v_\lb(t,\w,x)-\muk_\lb(\wt))=0\}$ (see (ii)).
We will prove in Lemma \ref{lema:4m} that $\mm_\lb^-$ is an upper semicontinuous $\tau_\lb$-equilibrium,
and that $\lb\mapsto\mm_\lb^-(\w)$ is strictly decreasing and left-continuous on $(\lb_-,\lb^+]$ for all $\w\in\W$.
\par
If $\mr_\lb^-\ge 0$ (which is the case if $\lb\in(\lb_-,\mu^+]$),
then Lemma \ref{lema:4copy}(ii) ensures
that $(\ma_\lb^-,\mr^-_\lb)$ is a pair of hyperbolic $\tau_\lb$-copies of the base and that
$v_\lb(t,\w,x)=v_\lb^-(t,\w,x)$ whenever $x\ge \mr^-_\lb(\w)$.
Hence, $\mm_\lb^-\le\mr_\lb^-$: see Theorem \ref{th:3CLbifur}(i). In addition, if $x<\mr_\lb^-(\w)$,
then $v_\lb(t,\w,x)<v_\lb(t,\w,\mr_\lb^-(\w))=\mr_\lb^-(\wt)$, and the uniform separation
of $\mr_\lb^-$ and $\ma_\lb^-$ precludes $\lim_{t\to\infty}(v_\lb(t,\w,x)-\ma_\lb^-(\wt))=0$. Therefore,
$\mm_\lb^-=\mr_\lb^-$.
\par
Let us fix $\lb\in[\mu_-,\lb^+]$. The inequality $\mg\le\mg^+$ and a standard comparison argument
ensure $v_\lb(t,\w,x)\le\mr_\lb^+(\wt)$ for all $t\ge 0$
whenever $x<\mr_\lb^+(\w)$. Since $\mr_\lb^+$ and $\ma_\lb^-$ are
uniformly separated (as Proposition \ref{prop:4tapas} ensures),
we conclude that $\mr^+_\lb\le \mm_\lb^-$. Now we take $(\bar\w,\bar x)$ with $\bar x>\mr^+_\lb(\bar\w)$.
The goal is to check that $\lim_{t\to\infty}(v_\lb(t,\bar\w,\bar x)-\ma_\lb^-(\bwt))=0$, which yields $\mr^+_\lb=\mm_\lb^-$.
We will prove below the existence of $s\ge 0$ with $v_\lb(s,\bar\w,\bar x)>0$.
Then, since $\lb\ge\mu_-$, we have $v_\lb(s,\bar\w,\bar x)>\mb_\lb(\bws)=\mr_\lb^-(\bws)$: see
Proposition \ref{prop:4tapas}(iii).
So, $v_\lb(t,\bar\w,\bar x)=v_\lb(t-s,\bar\w{\cdot}s,v_\lb(s,\bar\w,\bar x))\ge
v_\lb^-(t-s,\bar\w{\cdot}s,v_\lb(s,\bar\w,\bar x))$ for $t\ge s$, and the result follows from
$\lim_{t\to\infty}(v_\lb^-(t-s,\bar\w{\cdot}s,v_\lb(s,\bar\w,\bar x))-\ma_\lb^-(\bws{\cdot}(t-s)))=0$
(see Theorem \ref{th:3CLbifur}(i)).
Hence, there is nothing left to check if $\bar x>0$. We assume $\bar x\le 0$ and, for contradiction that
$v_\lb(s,\bar\w,\bar x)\le 0$ for all $s\ge 0$.
Then, the linear-convex analogue of Lemma
\ref{lema:4copy}(i) ensures that $v_\lb(t,\bar\w,\bar x)=v_\lb^+(t,\bar\w,\bar x)$ for all
$t\ge 0$, so $v_\lb^+(t,\bar\w,\bar x)$ is bounded as time increases (see Remark \ref{rm:4existeA}), while
the linear-convex analogue of Theorem \ref{th:3CLbifur}(i) ensures that $v_\lb^+(t,\bar\w,\bar x)$
is also bounded as time decreases. But this contradicts the definition of $\mr_\lb^+$, which
completes the proof of (v).
That of (vi) is analogous, with $\mm_\lb^+(\w):=\sup\{x\in\R\,|\;\lim_{t\to\infty}(v_\lb(t,\w,x)-\ma_\lb^+(\wt))=0\}
=\sup\{x\in\R\,|\;\lim_{t\to\infty}(v_\lb(t,\w,x)-\ml_\lb(\wt))=0\}$
for $\lb\in[\lb_-,\lb^+)$.
\par
The first assertions in (vii) are immediate consequences of (v) and (vi).
The compactness of $\{(\w,x)\mid\, \mm_\lb^+(\w)\le x\le\mm_\lb^-(\w)\}$ follows from the
lower semicontinuity of $\mm_\lb^+$ and the upper semicontinuity of $\mm_\lb^-$,
and the $\tau_\lb$-invariance follows from the fact that they are $\tau_\lb$-equilibria.
The repulsive character of the compact is ensured by the definition given in the previous
paragraph for $\mm_\lb^-$ and $\mm_\lb^+$.
\smallskip\par
(viii) We work with $\mm_\lb^-$: the proof for $\mm_\lb^+$ is analogous.
Assume for contradiction that there exists $m\in\merg$ with
$\int_\W\mh_x(\w,\mm_\lb^-(\w))\, dm=-\rho<0$. Birkhoff's Ergodic Theorem provides
a $\sigma$-invariant set $\W_1\subseteq\W$ with $m(\W_1)=1$ such that, for each
$\w\in\W_1$, there exists $t_\w>0$ with $\int_0^t\mh_x(\ws,\mm_\lb^-(\ws))\,ds\le -(\rho/2)\,t$
for all $t\ge t_\w$. Let us fix $\w_0\in\W_1$ and $t_0=t_{\w_0}$.
We take $\delta>0$ such that $|\mh_x(\w,\mm_\lb^-(\w))-\mh_x(\w,x)|<\rho/4$ for all $\w\in\W$
and $|\mm_\lb^-(\w)-x|<\delta$, and $\ep\in(0,\delta)$ such that $|v_\lb(t,\w,x)-\mm_\lb^-(\wt)|<\delta$
for all $\w\in\W$, $t\in[0,t_0]$ and $|x-\mm_\lb^-(\w)|<\ep$.
Then,
$|\mh_x(\wt,\mm_\lb^-(\wt))-\mh_x(\wt,v_\lb(t,\w,x))|<\rho/4$ for all $\w\in\W$, $t\in[0,t_0]$ and $|x-\mm_\lb^-(\w)|<\ep$,
and hence, since $(v_\lb)_x(0,\w_0,x)=1$,
\begin{equation}\label{eq:4integral}
\begin{split}
 (v_\lb)_x(t_0,\w_0,x)&=\exp\int_0^{t_0}\mh_x(\w_0{\cdot}s,v_\lb(s,\w_0,x))\, ds\\
 &<\exp\left((\rho/4)\,t_0+\int_0^{t_0}\mh_x(\w_0{\cdot}s,\mm_\lb^-(\w_0{\cdot}s))\, ds\right)\le e^{-(\rho/4)\, t_0}
\end{split}
\end{equation}
if $|x-\mm_\lb^-(\w_0)|<\ep$.
On the other hand, the Mean Value Theorem ensures that, if $|x-\mm_\lb^-(\w_0)|<\ep$,
then there exists $y\in(\mm_\lb^-(\w_0),x)$ such that $|v_\lb(t_0,\w_0,x)-\mm_\lb^-(\w_0{\cdot}t_0)|=
(v_\lb)_x(t_0,\w_0,y)\,|x-\mm_\lb^-(\w_0)|<\ep\, e^{-(\rho/4)\, t_0}<\ep$.
So, the choice of $\ep$ ensures that $|v_\lb(t,\w_0,x)-\mm_\lb^-(\w_0{\cdot}t)|<\delta$ for all $t\in[0,2t_0]$ if
$|x-\mm_\lb^-(\w)|<\ep$. This means that \eqref{eq:4integral} holds for $2t_0$ instead of $t_0$,
providing $|v_\lb(2t_0,\w_0,x)-\mm_\lb^-(\w_0{\cdot}(2t_0))|<\ep\, e^{-(\rho/4)\, 2t_0}$ if $|x-\mm_\lb^-(\w_0)|<\ep$.
An induction procedure yields $|v_\lb(nt_0,\w_0,x)-\mm_\lb^-(\w_0{\cdot}(nt_0))|<\ep\, e^{-(\rho/4)\,nt_0}$ if $|x-\mm_\lb^-(\w_0)|<\ep$.
This is a contradiction with the fact that $\lim_{t\to\infty}(v_\lb(t,\w_0,x)-\ma_\lb^-(\w_0{\cdot}t))=0$ if $x>\mm_\lb^-(\w_0)$.
\end{proof}
The next lemma has been used in the proof of points (v) and (vi). The conditions and notation are those of Theorem \ref{th:4CCbifur}.
\begin{lema}\label{lema:4m}
In the conditions of Theorem {\rm\ref{th:4CCbifur}},
\begin{itemize}[leftmargin=22pt]
\item[\rm(i)] the map $\mm_\lb^-$ is an upper semicontinuous $\tau_\lb$-equilibrium for $\lb\in(\lb_-,\lb^+]$,
    and $\lb\mapsto\mm_\lb^-(\w)$ is strictly decreasing and left-continuous on $(\lb_-,\lb^+]$ for all $\w\in\W$.
\item[\rm(ii)]  The map $\mm_\lb^+$ is a lower semicontinuous $\tau_\lb$-equilibrium,
    and $\lb\mapsto\mm_\lb^-(\w)$ is strictly decreasing and left-continuous on $[\lb_-,\lb^+)$ for all $\w\in\W$.
\end{itemize}
\end{lema}
\begin{proof}
We will prove (i). Note that $\mm_\lb^-(\w)$ is bounded for each $\lb\in(\lb_-,\lb^+]$.
In fact, $\ma^-_\lb\ge\mm_\lb^-\ge\ml_\lb$. The first inequality is clear, and the second one
follows from $\ml_\lb\le\ma_\lb^-$ (see Proposition \ref{prop:4tapas}(iii) and Theorem \ref{th:3CLbifur}(ii))
and from the properties of the global attractor.
It is easy to check that $v_\lb(t,\w,x)>\mm_\lb^-(\wt)$ if and only if $x>\mm_\lb^-(\w)$, so that
$v_\lb(t,\w,\mm_\lb^-(\w))=\mm_\lb^-(\wt)$; i.e., $\mm_\lb^-$ is a $\tau_\lb$-equilibrium.
\par
Let us fix $\lb_0\in(\lb_-,\lb^+]$.
The uniform asymptotic stability of $\ma_{\lb_0}^-=\muk_{\lb_0}$ (see Theorem \ref{th:4CCbifur}(i)-(ii))
provides $\delta_0>0$ such that,
if $\w\in\W$ and $|x-\ma^-_{\lb_0}(\w)|<\delta_0$, then
$\lim_{t\to\infty}(v_{\lb_0}(t,\w,x)-\ma^-_{\lb_0}(\wt))=0$. These
values of $\lb_0$ and $\delta_0$ will be fixed for the rest of the proof.
\par
Now we will check the upper semicontinuity of $\mm_{\lb_0}^-$.
We fix $\w_0\in\W$, a sequence $(\w_n)$ with limit $\w_0$, and $x>\mm_{\lb_0}^-(\w_0)$.
The goal is to check that $x\ge\mm_{\lb_0}^-(\w_n)$ if $n $ is large enough, since
this ensures $\mm_{\lb_0}^-(\w_0)\ge\limsup_{n\to\infty}\mm_{\lb_0}^-(\w_n)$.
We have $\lim_{t\to\infty}(v_{\lb_0}(t,\w_0,x)-\ma^-_{\lb_0}(\w_0{\cdot}t))=0$,
so there exists $s>0$ such that $|v_{\lb_0}(s,\w_0,x)-\ma^-_{\lb_0}(\w_0{\cdot}s)|<\delta_0/3$.
Hence, $|v_{\lb_0}(s,\w_n,x)-\ma^-_{\lb_0}(\w_n{\cdot}s)|<\delta_0$ for large enough $n$,
which yields $\lim_{t\to\infty}(v_{\lb_0}(t,\w_n,x)-\ma^-_{\lb_0}(\w_n{\cdot}t))=0$. Therefore,
$x\ge\mm_{\lb_0}^-(\w_n)$, as asserted.
\par
To prove the strictly decreasing character with respect to $\lb$ it suffices
to check that $\mm_{\lb}>\mm_{\lb_0}$ if $\lb<\lb_0$ and $\lb_0-\lb$ is small enough.
Let us first check that $\mm_{\lb}\ge\mm_{\lb_0}$.
We fix $\lb<\lb_0$ close enough to $\lb_0$ to ensure that $|\ma^-_{\lb_0}(\w)-\ma^-_{\lb}(\w)|<\delta_0/2$
for all $\w\in\W$ (see Theorem \ref{th:3CLbifur}(i)).
We take $\w_0\in\W$ and $x>\mm_\lb^-(\w_0)$. Let us check that $x>\mm_{\lb_0}^-(\w_0)$. As before, we take
$s>0$ such that $|v_\lb(s,\w_0,x)-\ma_\lb^-(\w_0{\cdot}s)|<\delta_0/2$. Hence,
$|v_\lb(s,\w_0,x)-\ma^-_{\lb_0}(\w_0{\cdot}s)|\le |v_\lb(s,\w_0,x)-\ma^-_\lb(\w_0{\cdot}s)|+
|\ma^-_\lb(\w_0{\cdot}s)-\ma^-_{\lb_0}(\w_0{\cdot}s)|<\delta_0/2+\delta_0/2=\delta_0$.
Since $v_{\lb_0}(s,\w_0,x)>v_\lb(s,\w_0,x)$, $\ma^-_{\lb_0}(\w_0{\cdot}s)-v_{\lb_0}(s,\w_0,x)<\ma^-_{\lb_0}(\w_0{\cdot}s)-v_{\lb}(s,\w_0,x)<\delta_0$.
If $|\ma^-_{\lb_0}(\w_0{\cdot}s)-v_{\lb_0}(s,\w_0,x)|<\delta_0$, then the choice of $\delta_0$ yields
$\lim_{t\to\infty}(v_{\lb_0}(t,\w_0,x)-\ma^-_{\lb_0}(\w_0{\cdot}t))=0$, and hence
$x>\mm_{\lb_0}^-(\w_0)$. Otherwise, $v_{\lb_0}(s,\w_0,x)>\ma^-_{\lb_0}(\w_0{\cdot}s)$, which is only possible
if $x>\ma^-_{\lb_0}(\w_0)\ge\mm^-_{\lb_0}(\w_0)$. So, $\mm^-_{\lb}\ge\mm^-_{\lb_0}$, as asserted. Now
we assume for contradiction that $\mm^-_{\lb}(\w_0)=\mm^-_{\lb_0}(\w_0)$ for an $\w_0\in\W$. Since
$(\mm^-_{\lb})'(\w)<\mc(\w)+\md(\w)\,\mm^-_\lb(\w)+\mg(\w,\mm^-_\lb(\w))+\lb_0$ for all $\w\in\W$,
a standard comparison argument shows that $\mm^-_\lb(\w_0{\cdot}t)<\mm^-_{\lb_0}(\w_0{\cdot}t)$
for all $t>0$, which contradicts $\mm^-_{\lb}\ge\mm^-_{\lb_0}$.
\par
Finally, let us prove that $\lb\mapsto\mm^-_\lb(\w)$ is left-continuous at $\lb_0$ for all $\w\in\W$.
It is easy to check that the map $\tilde\mm^-(\w):=\lim_{\lb\to(\lb_0)^-}\mm^-_\lb(\w)$ is a
$\tau_{\lb_0}$-equilibrium, and clearly $\tilde\mm^-\ge \mm^-_{\lb_0}$.
We take $\lb\in(\lb_-,\lb_0)$,
so that $\tilde\mm^-<\mm^-_{\lb}\le\ma^-_\lb\le\ma^-_{\lb_0}-\rho$, where $\rho:=\inf_{\w\in\W}(\ma^-_{\lb_0}(\w)-
\ma^-_\lb(\w))>0$ (see Theorem \ref{th:3CLbifur}(i)). This precludes
$\lim_{t\to\infty}(v_{\lb_0}(t,\w,\tilde\mm^-(\w))-\ma^-_{\lb_0}(\wt))=
\lim_{t\to\infty}(\tilde \mm^-(\wt)-\ma^-_{\lb_0}(\wt))=0$ for any $\w\in\W$,
and hence $\tilde\mm^-(\w)\le \mm^-_{\lb_0}(\w)$.
\end{proof}
Under the conditions of Theorem \ref{th:4CCbifur}, the simplest global dynamics
for the family \eqref{eq:4conconv} occurs when $\lb_*=\lb_-$ and $\lb^*=\lb^+$ (as in the case of a minimal
base flow), and when the maps $\mm_\lb^+$ and $\mm_\lb^-$ of
points (v) and (vi) coincide for all $\lb\in(\lb_-,\lb^+)$. In this case, it follows from
Theorem \ref{th:4CCbifur} that $\{\mm_\lb\}$ is a $\tau_\lb$-copy of the base
for all $\lb\in(\lb_-,\lb^+)$, with $\{\mm_\lb\}$ repulsive,
and that the map $(\lb_-,\lb^+)\to C(\W,\R),\,\lb\mapsto\mm_\lb$ is continuous and strictly increasing.
The corresponding global
diagram, which is the unique possible one in the autonomous case (as we will explain in Section \ref{subsec:4todos}),
is depicted in Figure \ref{fig:SFina}.
\begin{figure}
     \centering
         \includegraphics[width=0.6\textwidth]{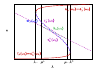}
         \caption{Depiction of the bifurcation diagram for \eqref{eq:4conconv} under the hypotheses of Theorem \ref{th:4CCbifur}
         in the simplest case, with $\mm_\lb^+=\mm_\lb^-$ for all $\lb\in(\lb_-,\lb^+)$. The two attractive hyperbolic copies of the base
         (in red) and the repulsive copy of the base (in blue)
         are depicted by their value at a fixed point $\w_0\in\W$.
         The copies of the base determining the behavior of the linear,
         concave-linear and linear-convex associated problems are also represented at
         $\w_0$. The large black points depict the fact that two possibilities
         arise at the bifurcation points: either $\mm_{\lb_-}(\w_0)=\muk_{\lb_-}(\w_0)$ or $\mm_{\lb_-}(\w_0)<\muk_{\lb_-}(\w_0)$,
         and either $\mm_{\lb^+}(\w_0)=\ml_{\lb^+}(\w_0)$ or $\mm_{\lb^+}(\w_0)>\ml_{\lb^+}(\w_0)$.}
        \label{fig:SFina}
\end{figure}

Assuming always \ref{cc1}-\ref{cc6}, but under additional conditions less restrictive than those of
Theorem \ref{th:4CCbifur}, some aspects of the evolution of the global dynamics can be
established. The next result provides a sample of this assertion.
The notion of local saddle-node bifurcation is explained
in Remark \ref{rm:3localSNB}.
\begin{prop}\label{prop:4signo}
Assume that \ref{cc1}-\ref{cc6}. Then,
\begin{itemize}[leftmargin=22pt]
\item[\rm(i)] if $\lb_-<\mu^+$, then \eqref{eq:4conconv}
has a local saddle-node bifurcation as $\lb\downarrow\lb_-$, due to the approach of $\muk_\lb=\ma^-_\lb$ and
$\mm_\lb=\mr^-_\lb$.
\item[\rm(ii)] If $\lb^+>\mu_-$, then \eqref{eq:4conconv}
has a local saddle-node bifurcation as $\lb\uparrow\lb^+$, due to the approach of $\mm_\lb=\ma_\lb^+$ and
$\ml_\lb=\ma_\lb^-$.
\end{itemize}
\end{prop}

\begin{proof}
Note that, if $\lb\in(\lb_-,\mu^+)$, then $\ma_\lb^->\mr_\lb^-\ge \mb_\lb\ge 0$ for $\lb\in(\lb_-,\mu^+)$:
see Proposition \ref{prop:4tapas}(iii) and \eqref{def:4mu}.
With this in mind, the first assertion follows easily from Lemma \ref{lema:4copy} and from
the bifurcation described in Theorem \ref{th:3CLbifur}. The proof of the second one is similar.
\end{proof}
%%%%%%%%%%%%%%%%%%%%%%%%%%%%%%%%%%%%%%%%%%%%%%%%%%%%%%%%%%%%%%%%%%
%%%%%%%%%%%%%%%%%%%%%%%%%%%%%%%%%%%%%%%%%%%%%%%%%%%%%%%%%%%%%%%%%%
\subsection{A d-concavity band with dynamical consequences}\label{subsec:5band}
In this section, the existence of a {\em d-concavity band\/} around $\W\times\{0\}$ for a regular enough
$\mg$ satisfying \ref{cc1}, \ref{cc4} and \ref{cc5}, as well as of a band of strict d-concavity, is shown.
Then, under hypotheses \ref{cc1}-\ref{cc6}, we establish conditions relating the maps $\mr^\pm_\lb$,
$\ma^\pm_\lb$, $\ml_\lb$ and $\muk_\lb$ to this d-concavity band in order to prove the existence of
local saddle-node bifurcations, of the type described in Remark \ref{rm:3localSNB}, using
the bifurcation results for d-concave flows proved in \cite{dno5}.  Recall that
$\tau_\lb^-$ and $\tau_\lb^+$ are the concave-linear and linear-convex skewproduct flows
associated to \eqref{eq:4conconv}$^\lb$, respectively defined from the solutions of
\eqref{eq:4conclin}$^\lb$ and \eqref{eq:4linconv}$^\lb$, respectively; and that $\lb_-$ and $\lb^+$ are
the bifurcation points of the corresponding bifurcation problems. The rest of the notation established
at the beginning of this section will be also used.
\begin{defi}\label{def:4band}
Assume that \ref{cc1}, \ref{cc4} and \ref{cc5} hold.
For each $\w\in\W$, we define the (possibly empty) subsets of $\R$
\[
\begin{split}
\mI_\alpha(\w)&:=\{x<0\mid\,\mg_{x}(\w,\cdot)\text{ is concave on }[x,0]\}\,,\\
\mI_\beta(\w)&:=\{x>0\mid\,\mg_{x}(\w,\cdot)\text{ is concave on }[0,x]\}\,,\\
\end{split}
\]
the maps $\alpha\colon\W\to[-\infty,0]$ and $\beta\colon\W\to[0,\infty]$
\[\alpha(\w):=\left\{\begin{array}{ll}
\inf\mI_\alpha(\w)\;&\text{if }\mI_\alpha(\w)\neq\emptyset\,,\\
0&\text{otherwise}\,,
\end{array}
\right.
\quad
\beta(\w):=\left\{\begin{array}{ll}
\sup\mI_\beta(\w)\;&\text{if }\mI_\beta(\w)\neq\emptyset\,,\\
0&\text{otherwise}\,,
\end{array}
\right.
\]
and the {\em d-concavity band} of $\mg$
\[
 \mD_\mg:=\{(\w,x)\mid\,\alpha(\w)\le x\le \beta(\w)\}\subseteq\WR\,.
\]
\end{defi}
\begin{defi}\label{def:4strictband}
Assume that \ref{cc1}, \ref{cc4} and \ref{cc5} hold.
For each $\w\in\W$, we define the (possibly empty) subsets of $\R$
\[
\begin{split}
\mI_\alpha^*(\w)&:=\{x<0\mid\,\mg_{x}(\w,\cdot)\text{ is strictly concave on }[x,0]\}\,,\\
\mI_\beta^*(\w)&:=\{x>0\mid\,\mg_{x}(\w,\cdot)\text{ is strictly concave on }[0,x]\}\,,\\
\end{split}
\]
the maps $\alpha^*\colon\W\to[-\infty,0]$ and $\beta^*\colon\W\to[0,\infty]$
\[\alpha^*(\w):=\left\{\begin{array}{ll}
\inf\mI_\alpha^*(\w)\;&\text{if }\mI_\alpha^*(\w)\neq\emptyset\,,\\
0&\text{otherwise}\,,
\end{array}
\right.
\quad
\beta^*(\w):=\left\{\begin{array}{ll}
\sup\mI_\beta^*(\w)\;&\text{if }\mI_\beta^*(\w)\neq\emptyset\,,\\
0&\text{otherwise}\,,
\end{array}
\right.
\]
and the {\em strict d-concavity band} of $\mg$
\[
\mD_\mg^*:=\{(\w,x)\mid\,\alpha^*(\w)\le x\le \beta^*(\w)\}\subseteq\WR\,.
\]
\end{defi}
It is obvious that $\alpha(\w)\le \alpha^*(\w)\le 0\le \beta^*(\w)\le \beta(\w)$ for all $\w\in\W$;
i.e., $\mD^*_\mg\subseteq\mD_\mg$. The next proposition shows that the nondegeneracy of these bands of
d-concavity and of strict d-concavity is a very general property.
In practice, it occurs in a wide range of functions beyond analytic functions.
\begin{prop}\label{prop:4band}
Assume that \ref{cc1}, \ref{cc4}, \ref{cc5} hold, and that $\mg\in C^{0,3}(\WR,\R)$.
Then,
\begin{itemize}[leftmargin=22pt]
\item[{\rm(i)}] the maps $\alpha$ and $\beta$ are respectively lower and upper semicontinuous,
and hence the set $\mD_\mg$ is closed.
\item[{\rm(ii)}] There exists a $\sigma$-invariant compact set $\W_0\subseteq\W$ with $m(\W_0)=1$ for all $m\in\merg$
such that, for all $\w\in\W_0$, $\mg_{xx}(\w,x)\le 0$ for all $x\ge 0$, 
$\mg_{xx}(\w,x)\ge 0$ for all $x\le 0$, $\mg_{xx}(\w,0)=0$ and $\mg_{xxx}(\w,0)\le 0$.
    If, in addition, $(\W,\sigma)$ is minimal, then $\W=\W_0$.
\item[{\rm(iii)}] If $\w\in\W_0$ and there exists $n_\w\in\N$ such that 
$x\mapsto\mg_\w(x):=\mg(\w,x)$ belongs to $C^{n_\w}(\R,\R)$ and satisfies $\mg_\w(0)=\mg'_\w(0)=\dots=\mg_\w^{(n_\w-1)}(0)=0\ne\mg_\w^{(n_\w)}(0)$, then
    $\alpha^*(\w)<0<\beta^*(\w)$. And this is the situation if $\w\in\W_0$ and $\mg_\w$ 
    is analytic and not identically $0$.
\end{itemize}
\end{prop}
\begin{proof}
(i) Let us work with $\alpha$.
For contradiction, we assume that $\w_0:=\lim_{n\to\infty}\w_n$ and that
$\alpha_0:=\lim_{n\to\infty}\alpha(\w_n)<\alpha(\w_0)$, which requires $\alpha_0<0$.
We take $x_0\in(\alpha_0,\alpha(\w_0))$, look for $n_0\in\N$ such that $\alpha(\w_n)\le x_0$ for all $n\ge n_0$,
and deduce from Definition \ref{def:4band} that $\mg_{xxx}(\w_n,x)\le 0$ for all $n\ge n_0$ and $x\in[x_0,0]$.
Hence, $\mg_{xxx}(\w_0,x)\le 0$ for all $x\in[x_0,0]$, which yields the contradiction $\alpha(\w)\le x_0$.
The proof is analogous for $\beta$, and the final assertion is clear.
\smallskip\par
(ii) We define $\mC:=\{\w\in\W\,|\;x\mapsto \mg(\w,x)$ is concave on $[0,\infty)$ and convex on $(-\infty,0]\}$
and $\W_0:=\cap_{s\in\Q}\,\mC{\cdot}s$, with $\mC{\cdot}s:=\sigma(s,C)$. The continuity of $\mg_{xx}$ shows that $\mC$ is closed,
and hence $\W_0$ is compact.
Given any $m\in\merg$, \ref{cc5} ensures that $m(\mC)=1$, so $m(\mC{\cdot}s)=1$ for all $s\in\R$, and therefore $m(\W_0)=1$.
We take $\w\in\W_0$ and $t\in\R$, write $t$ as the limit of a sequence $(s_n)$ in $\Q$, fix any $s\in\Q$,
and deduce from $\ws_n\in\mC_s$ for all $n$ that
$\wt=\lim_{n\to\infty}\w{\cdot}s_n$ is in the closed set $\mC_s$ and hence in $\W_0$, which is therefore $\sigma$-invariant.
By definition of $\mC\supseteq\W_0$, $\mg_{xx}(\w,x)\le 0$ for $x\ge 0$ and $\mg_{xx}(\w,x)\ge 0$ for $x\le 0$ if $\w\in\W_0$, and
hence $\mg_{xx}(\w,0)=0$ and $\mg_{xxx}(\w,0)\le 0$. The last assertion follows from the $\sigma$-invariance and compactness of $\W_0$.
\smallskip\par
(iii) We fix $\w\in\W_0$ and take $n_\w\in\N$ as in the statement. It follows from \ref{cc4} and (ii) that $n_\w\ge3$. According to
(ii), $\mg_\w''(x)=\mg_{xx}(\w,x)$ is nonincreasing at $x=0$, and hence $n_\w$ is odd and $\mg_\w^{(n_\w)}(0)<0$.
So, there exists $\ep>0$ such that $\mg_\w'''(x)=\mg_{xxx}(\w,x)<0$ for all $x\in[-\ep,0)\cup(0,\ep]$.
Therefore, $\mg'_\w(x)=\mg_\w(\w,x)$ is strictly concave on $[-\ep,\ep]$, and hence $\alpha^*(\w)<0<\beta^*(\w)$, as asserted.
The final assertion is a trivial consequence of the analyticity of $\mg$, its non constant character, and \ref{cc4},
according to which $\mg_\w(0)=\mg(\w,0)=0$.
\end{proof}
\begin{defi}\label{def:5includedinband}
Assume that \ref{cc1}, \ref{cc4} and \ref{cc5} hold.
We say that the graphs of two continuous maps $\mb_1,\mb_2\colon\W\to\R$ with $\mb_1(\w)<\mb_2(\w)$ for all $\w\in\W$
are {\em jointly included (in measure) in the strict d-concavity band $\mD^*_\mg$} if
\begin{itemize}[leftmargin=25pt]
\item[{\rm(1)}] $m(\{\w\in\W\mid\,[\mb_1(\w),\mb_2(\w)]\subseteq(\mD_\mg)_\w\})=1$ for all $m\in\merg$,
\item[{\rm(2)}] $m(\{\w\in\W\mid\,[\mb_1(\w),\mb_2(\w)]\subseteq(\mD^*_\mg)_\w\})>0$ for all $m\in\merg$.
\end{itemize}
\end{defi}
To fully understand the hypotheses of the next result, recall that
$\ml_{\lb_1}<\ml_{\lb_2}\le \ma_{\lb_2}^+<\mr_{\lb_2}^+\le \mr_{\lb_1}^-<\ma_{\lb_1}^-\le \muk_{\lb_1}<\muk_{\lb_2}$
if $\lb_-<\lb_1<\lb_2<\lb^+$:
see Proposition \ref{prop:4tapas}(iii).

\begin{teor}\label{th:4band}
Assume that \ref{cc1}-\ref{cc6} hold with $\mg\in C^{0,2}(\WR,\R)$,
and that $\lb_-<\lb^+$, and take $\lb_1<\lb_2$ in $(\lb_-,\lb^+)$.
\begin{itemize}[leftmargin=22pt]
\item[{\rm(i)}] If $\ml_{\lb_2}$ and $\mr_{\lb_1}^-$ are jointly included in
the strict d-concavity band, then \eqref{eq:4conconv} exhibits a local saddle-node
bifurcation as $\lb\uparrow\xi_+$, with $\xi_+\ge\lb^+$.
\item[{\rm(ii)}] If $\mr_{\lb_2}^+$ and $\muk_{\lb_1}$ are jointly included in
the strict d-concavity band, then \eqref{eq:4conconv} exhibits a local saddle-node
bifurcation as $\lb\downarrow\xi_-$, with $\xi_-\le \lb_-$.
\end{itemize}
\end{teor}
\begin{proof}
(i) We define
\[
\mg^*(\w,x):=\left\{\!\!\begin{array}{ll}
\mg(\w,\mr_{\lb_1}^-(\w))+\mg_x(\w,\mr_{\lb_1}^-(\w))\,(x-\mr_{\lb_1}^-(\w))&\\[1ex]
\qquad+\,(1/2)\,\mg_{xx}(\w,\mr_{\lb_1}^-(\w))\,(x-\mr_{\lb_1}^-(\w))^2\\[1ex]
\qquad-(x-\mr_{\lb_1}^-(\w))^3&\quad\text{if }x>\mr_{\lb_1}^-(\w)\,,\\[1ex]
\mg(\w,x)
&\quad\text{if }x\in[\ml_{\lb_2}(\w),\mr_{\lb_1}^-(\w)]\,,\\[1ex]
\mg(\w,\ml_{\lb_2}(\w))+\mg_x(\w,\ml_{\lb_2}(\w))\,(x-\ml_{\lb_2}(\w))&\\[1ex]
\qquad+\,(1/2)\,\mg_{xx}(\w,\ml_{\lb_2}(\w))\,(x-\ml_{\lb_2}(\w))^2\\[1ex]
\qquad-(x-\ml_{\lb_2}(\w))^3&\quad\text{if }x<\ml_{\lb_2}(\w)\,.
\end{array}\right.
\]
Note that $\mg^*\in C^{0,2}(\WR,\R)$, and that $\lim_{x\to\pm\infty}(\pm\mg^*(\w,x))=-\infty$
uniformly on $\W$. Since $\{\w\in\W\mid\, [\ml_{\lb_2}(\w),\mr_{\lb_1}^-(\w)]\subseteq(\mD_\mg)_\w\}
\subseteq\{\w\in\W\mid\,\mg_x(\w,\cdot)$ is concave on $[\ml_{\lb_2}(\w),\mr_{\lb_1}^-(\w)]\}=
\{\w\in\W\mid\,\mg_x^*(\w,\cdot)$ is concave$\}$, Definition \ref{def:5includedinband}(1)
ensures that $m(\{\w\in\W\mid\,x\mapsto\mg_x^*(\w,x)$ is concave$\})=1$ for all $m\in\merg$.
Analogously, Definition \ref{def:5includedinband}(2) yields
$m(\{\w\in\W\mid\,x\mapsto\mg_x^*(\w,x)$ is strictly concave$\})>0$ for all $m\in\merg$.

Proposition \ref{prop:4tapas}(i) ensures that $\{\ml_{\lb_2}\}$ is an attractive $\tau_{\lb_2}$-hyperbolic copy of the base.
In addition, the definition of $\mg^*$ ensures that, for each $m\in\merg$, the corresponding Lypaunov exponent
for $\tau^*_{\lb_2}$ coincides with that for $\tau_{\lb_2}$, and hence it is negative. According to Theorem \ref{th:2copia},
this ensures that $\{\ml_{\lb_2}\}$ is also an attractive $\tau^*_{\lb_2}$-hyperbolic copy of the base.

On the other hand, since $\mg^-\le\mg$, we have $(\mr_{\lb_1}^-)'(\w)\le
\mc(\w)+\md(\w)\,\mr_{\lb_1}^-(\w)+\mg(\w,\mr_{\lb_1}^-(\w))+\lb_1<\mc(\w)+\md(\w)\,\mr_{\lb_1}^-(\w)+
\mg^*(\w,\mr_{\lb_1}^-(\w))+\lb_2$. So, $\mr_{\lb_1}^-$ is a continuous
strict global lower solution for $\tau_{\lb_2}^*$. Consequently,
the maps $\mb^-_{\lb_2}(\w):=\lim_{s\to\infty}v_{\lb_2}^*(-s,\ws,\mr_{\lb_1}^-(\ws))$ and
$\mb^+_{\lb_2}(\w):=\lim_{s\to\infty}v_{\lb_2}^*(s,\w{\cdot}(-s),\mr_{\lb_1}^-(\w{\cdot}(-s)))$
are two bounded $\tau^*_{\lb_2}$-equilibria, upper and lower semicontinuous respectively, and with
$\mb^-_{\lb_2}(\w)<\mr_{\lb_1}^-(\w)<\mb^+_{\lb_2}(\w)$ for all $\w\in\W$ (see, e.g., \cite[Secion 2.2]{dno1} and
references therein). In addition, $v_{\lb_2}^*(-s,\ws,\mr_{\lb_1}^-(\ws))>v_{\lb_2}^*(-s,\ws,\ml_{\lb_2}(\ws))
=\ml_{\lb_2}(\w)$ for all $s$, and hence $\mb^-_{\lb_2}\ge\ml_{\lb_2}$. Let us
check that $\mb^-_{\lb_2}$ and $\ml_{\lb_2}$ are uniformly separated. Let $(k,\gamma)$ and $\rho$
be a dichotomy constant pair and a radius of uniform convergence for $\{\ml_{\lb_2}\}$
(with $k\ge 1$, $\gamma>0$ and $\rho>0$),
and assume for contradiction the existence of $\bar\w\in\W$ such that $\mb^-_{\lb_2}(\bar\w)-\ml_{\lb_2}(\bar\w)<\rho/2$.
Then, for large enough $s>0$, $|v^*_{\lb_2}(-s,\bws,\mr_{\lb_1}^-(\bws))-\ml_{\lb_2}(\bar\w)|<\rho$, and hence
$|\mr_{\lb_1}^-(\bar\w)-\ml_{\lb_2}(\bar\w)|<k\,e^{-\gamma s}$, which provides the contradiction
$\mr_{\lb_1}^-(\bar\w)=\ml_{\lb_2}(\bar\w)$. This uniform separation and the semicontinuity properties
imply that the closures of the graphs of $\ml_{\lb_2}$, $\mb_{\lb_2}^-$ and
$\mb_{\lb_2}^+$ define three compact $\tau_{\lb_2}^*$-invariant sets which project onto $\W$.
Therefore, \cite[Theorem 5.3]{dno4} guarantees that they are three hyperbolic
$\tau_{\lb_2}^*$-copies of the base, as asserted.
For next purposes, we point out that it also ensures that
$\ml_{\lb_2}$ and $\mb_{\lb_2}^+$ are the lower and upper equilibria of the global attractor for $\tau^*_{\lb_2}$,
and, consequently, that $\ml_{\lb_2}<\ml^*_\lb$ and $\mb_{\lb_2}^+<\muk^*_\lb$
for any $\lb>\lb_2$, where $\ml^*_\lb$ and $\muk^*_\lb$ are the lower and upper
equilibria of the global attractor of $\tau^*_\lb$: this follows, for instance, from
\cite[Proposition 5.5(iv)]{dno4}.

In this situation, \cite[Theorem 4.4]{dno5} describes the complete bifurcation diagram for
$\tau_\lb^*$, ensuring that a nonautonomous local saddle-node bifurcation takes place at a point
$\xi^+>\lb_2$, with involved hyperbolic copies of the base lying within the region $\bigcup_{\w\in\W}\big(\{\w\}\times[\ml_{\lb_2}(\w),\mb_{\lb_2}^-(\w)]\big)\subset
\bigcup_{\w\in\W}\big(\{\w\}\times[\ml_{\lb_2}(\w),\mr_{\lb_1}^-(\w)]\big)$.
(See \cite[Figure 6]{dno5}: $\xi^+$ is the point $\lb^+$ of that depiction.)
Hence, it is also a local nonautonomous saddle-node bifurcation for $\tau_\lb$, as stated.
We point out here that
\cite[Theorem 4.4]{dno5} also establishes that $\inf_{\w\in\W}(\muk^*_\lb(\w)-\ml_\lb^*(\w))=0$
if $\lb>\xi^+$.

Our next goal is to show that $\xi^+\ge\lb^+$. We assume for contradiction the existence of $\lb\in(\xi^+,\lb^+)$.
In particular, $\lb^+>\lb>\lb_2>\lb_1$, which ensures the existence and hyperbolicity of the $\tau_\lb^+$-copy of the base
$\{\ma_\lb^+\}$. Then $\mr_{\lb_1}^->\ma^+_\lb>\ml_{\lb_2}$: see Proposition \ref{prop:4tapas}(iii).
In particular, $\mg^*(\w,\ma_\lb^+(\w))=\mg(\w,\ma_\lb^+(\w))$ for all $\w\in\W$, which combined
with the fact that $\{\ma_\lb^+\}$ is a $\tau^+_\lb$-copy of the base and with $\mg^+\ge\mg$ yield $(\ma^+_\lb)'(\w)\ge\mc(\w)+\md(\w)\,x+\mg(\w,\ma_\lb^+(\w))+\lb=
\mc(\w)+\md(\w)\,x+\mg^*(\w,\ma_\lb^+(\w))+\lb$. This inequality and
\cite[Proposition 5.5(iv)]{dno4} applied to $\tau_\lb^*$ ensure $\ml_\lb^*\le\ma^+_\lb$,
which can be extended to $\ml_\lb^*\le\ma^+_\lb<\mr^+_\lb<\mr_{\lb_1}^-<\mb_{\lb_2}^+<\muk^*_\lb$:
see Proposition \ref{prop:4tapas}(iii) and recall the already established properties of
$\mb_{\lb_2}^+$. So, the
uniform separation of $\ma^+_\lb$ and $\mr^+_\lb$ implies the same property for
$\ml_\lb^*$ and $\muk^*_\lb$. This contradicts the last sentence of the previous paragraph,
and so $\lb^+\le\xi^+$, as asserted. The proof of (i) is complete.\hspace{-1cm}~
\smallskip

(ii) The argument to find the local saddle-node bifurcation as $\lb\downarrow\xi_-$
is analogous to that in the proof of (i) (and now $\xi_-$ is the point $\lb_-$ of
\cite[Figure 6]{dno5}).
\end{proof}
We recall that $\ma_\lb^+=\ml_{\lb}$ (resp.~$\ma_\lb^-=\mu_{\lb}$) whenever $\ma^+_\lb\le 0$ (resp.~$\ma^-_\lb\ge 0$):
see Proposition \ref{prop:4tapas}(i) (resp.~(ii)). These properties allow us to
reformulate
Theorem \ref{th:4band} in terms of $\ma^+_\lb$ and $\ma^-_\lb$.
In addition, the previous proof can be easily adapted to the next result,
which does not require $\lb_-<\lb^+$, and which completes this section.
\begin{teor}\label{th:4band-ext}
Assume that \ref{cc1}-\ref{cc6} hold with $\mg\in C^{0,2}(\WR,\R)$.
\begin{itemize}[leftmargin=22pt]
\item[{\rm(i)}] If, for a $\lb_2<\lb^+$, there exists a strict global continuous
$\tau_{\lb_2}$-lower solution $\ms_{\lb_2}\colon\W\to\R$ with $\ms_{\lb_2}>\ml_{\lb_2}$
and such that $\ml_{\lb_2}$ and $\ms_{\lb_2}$ are jointly included in the strict
d-concavity region $\mD^*_\mg$, then \eqref{eq:4conconv}
exhibits a local saddle-node bifurcation as $\lb\uparrow\xi^+$, with $\xi^+>\lb_2$.
\item[{\rm(ii)}] If, for a $\lb_1>\lb_-$, there exists a strict global continuous $\tau_{\lb_1}$-upper solution $\ms_{\lb_1}\colon\W\to\R$
with $\ms_{\lb_1}<\muk_{\lb_1}$ and such that $\ms_{\lb_1}$ and $\muk_{\lb_1}$ are jointly included in the strict
d-concavity region $\mD^*_\mg$, then \eqref{eq:4conconv} exhibits a local saddle-node bifurcation as $\lb\downarrow\xi_-$,
with $\xi_-<\lb_1$.
\end{itemize}
\end{teor}
\begin{proof} To prove (i), we proceed as in the proof of Theorem \ref{th:4band}(i) with $\ms_{\lb_2}$ in the
role of $\mr_{\lb_1}^-$. The proof of (ii) is analogous to that of Theorem \ref{th:4band}(ii).
\end{proof}
%%%%%%%%%%%%%%%%%%%%%%%%%%%%%%%%%%%%%%%%%%%%%%%%%%%%%%%%%%%%%%%%%%%%
%%%%%%%%%%%%%%%%%%%%%%%%%%%%%%%%%%%%%%%%%%%%%%%%%%%%%%%%%%%%%%%%%%%%
\subsection{Different possibilities for a transitive flow} \label{subsec:4todos}
Let us assume now that $(\W,\sigma)$ is transitive, and let $\bar\w$
be a point with dense $\sigma$-orbit.
Then, the parameters $\lb_-$ and $\lb^+$ of Theorem
\ref{th:4CCbifur} satisfy $\lb_-=\lb_-^{\bar\w}$ and $\lb^+=\lb^+_{\bar\w}$ (as seen in
Theorem \ref{th:3CLbifuruna} and in its linear-convex analogue).
In addition, since $\mb_\lb$ is continuous,
$\mu_-=\mu_-^{\bar\w}$ and $\mu^+=\mu^+_{\bar\w}$, being these parameters defined in
\eqref{def:4mu}, in Theorem \ref{th:3CLbifuruna}(iv) and in Remark \ref{rm:3convex}.
In this case,
\[
 \mu^+\le\mu_-,\,\qquad \lb_*\le\lb_-<\mu_-\,,\qquad  \mu^+<\lb^+\le\lb^*\,,
\]
as deduced from Proposition \ref{prop:3CL-DE}, Theorems \ref{th:3CLbifur} and \ref{th:3CLbifuruna},
and their linear-convex analogues.
These inequalities leave only five possible orders for $\lb_-,\,\lb^+\!,\,\mu_-$ and $\mu^+\!$:
\begin{enumerate}[leftmargin=23pt,label=\rm{\bf{o\arabic*}}]
\item\label{o1} $\lb_-<\mu^+\le\mu_-<\lb^+$,
\item\label{o2} $\mu^+\le\lb_-<\mu_-\le\lb^+$,
\item\label{o3} $\lb_-\le\mu^+<\lb^+\le\mu_-$,
\item\label{o4} $\mu^+<\lb_-<\lb^+<\mu_-$,
\item\label{o5} $\mu^+<\lb^+\le\lb_-<\mu_-$.
\end{enumerate}
Note that
\begin{itemize}[leftmargin=15pt]
\item[-] if the linear part of the family \eqref{eq:4conconv} is autonomous
(that is, if $\mc(\w)\equiv c\in\R$ and $\md(\w)\equiv d>0$),
then $\bar\mb_\lb\equiv-(\lb+c)/d$. So, $\mu_-=\mu_+=-c$, and hence the order in this case is that of \ref{o1}.
\item[-] In the case of a transitive flow, the hypotheses of Theorem \ref{th:4CCbifur} involve case \ref{o1}.
\item[-] Let $\mM\subseteq\W$ be a minimal subset, and let $\bar\lb_-,\,\bar\lb^+,\bar\mu_-,\bar\mu^+$ be the corresponding values of the
parameter. Theorem \ref{th:3CLbifur}(vii)
shows that $\bar\lb_-\le\lb_-$, and it is easy to deduce from \eqref{def:4mu} that
$\bar\mu^+\ge\mu^+$. So, if $\lb_-<\mu^+$, then $\bar\lb_-<\bar\mu^+$. Analogously, if $\mu_-<\lb^+$, then $\bar\mu_-<\bar\lb^+$.
\end{itemize}

Our next purpose is to show that cases \ref{o2}, \ref{o3}, \ref{o4}
and \ref{o5} can indeed occur, for which we combine theoretical results with a basic
numerical analysis.
We will construct the examples assuming one of the simplest types of
time-dependence under which they can exist: a constant map $\md$, a time-independent
map $\mg$, and a map $\mc$ defined over the circle $\W=\s^1=\R/(2\pi\Z)$.
More precisely, our family will be obtained by the hull procedure
(see Section \ref{subsec:2process}) from the equation
\begin{equation}\label{eq:4casos}
 x'=d\,x+\cos(t)+g(x)+\lb\,,
\end{equation}
where $d\in\R$ satisfies $d>0$, and where $g\colon\R\to\R$ is $C^2(\R,\R)$ and satisfies $g(0)=g'(0)=g''(0)=0$, $g'(x)<0$ for $x\ne 0$, $g''(x)<0$ for $x>0$ and $g''(x)>0$ for $x<0$.
So, if $g^-(x):=\min(0,g(x))$ and $g^+(x):=\max(0,g(x))$, then $g(x)=g^\pm(x)$ if $\mp x\ge 0$.
In addition, the family on the hull is $x'=d\,x+\cos(t+s)+g(x)+\lb$ for $s\in[0,2\pi)$,
and it is easy to check that \ref{cc1}-\ref{cc6} hold. In particular, the base flow is minimal.
Therefore, there are only four parameters: $\lb_*=\lb_-<\lb^+=\lb^*$ (see Theorem \ref{th:3CLbifur}(vi))
and $\mu_-\ge\mu^+$. In this case, \eqref{eq:4mu-mu+} provides
\begin{equation}\label{ig:4mu}
\begin{split}
\mu_-&=-\min_{t\in[0,2\pi]}\frac{d}{d^2+1}\,\big(d\cos(t)-\sin(t)\big)=\frac{d}{\sqrt{d^2+1}}\,,\\
\mu^+&=-\max_{t\in[0,2\pi]}\frac{d}{d^2+1}\,\big(d\cos(t)-\sin(t)\big)=-\frac{d}{\sqrt{d^2+1}}\,.
\end{split}
\end{equation}
We consider the Poincar\'{e} maps $T_\lb^\pm\colon\R\to\R$ of $x'=d\,x+\cos(t)+g^\pm(x)+\lb$ to
obtain information about the relative position of the parameters $\lb_-$, $\lb^+$, $\mu_-$ and $\mu^+$. So,
$T_\lb^\pm(x):=x_\lb^\pm(2\pi,0,x)$, where $x_\lb^\pm$ is the solution of the equation with value $x$
at $t=0$. Since the equation for the map $g^-$ is in the concave-linear case, Theorem \ref{th:3CLbifur} ensures
that it has: no bounded solutions for $\lb<\lb_-$ (which means no fixed points for $T_\lb^-$);
one periodic solution for $\lb=\lb_-$  (one fixed point for $T_\lb^-$)
since, in the scalar periodic case,
any minimal set for the skewproduct flow is given by the closure of just one orbit, and hence
it is a copy of the base representing a periodic solution;
and two periodic solutions for $\lb>\lb^-$ (two fixed points for $T_\lb^-$),
since there are two copies of the base.
Analogously, $T_\lb^+$ has two fixed points for $\lb<\lb^+$, one for $\lb=\lb^+$, and none for $\lb>\lb^+$.
So, the number of fixed points of $T^-_{\mu_+}$ (resp.~$T^+_{\mu^-}$) determines the relative position
of $\mu^+$ and $\lb_-$ (resp.~$\mu_-$ and $\lb^+$).  Table~\ref{tabla:posibilidades} partly reflects this fact.

\begin{table}[h]
\begin{center}
\begin{tabular}{|c|c|c|c|}
\hline
$\#\mathrm{fix}(T^-_{\mu^+})$ & $\#\mathrm{fix}(T^+_{\mu_-})$ & parameter order & case \\ \hline \hline
2 & 2 & $\lb_-<\mu^+<\mu_-<\lb^+$ & \ref{o1} \\
0 & 2 & $\mu^+<\lb_-<\mu_-<\lb^+$ & \ref{o2} \\
2 & 0 & $\lb_-<\mu^+<\lb^+<\mu_-$ & \ref{o3} \\
0 & 0 & $\mu^+<\lb_-$ and $\lb^+<\mu_-$ & \ref{o4} or \ref{o5} \\ \hline
\end{tabular}
\end{center}\par
\caption{Relative orders of the parameters $\lb_-$, $\lb^+$, $\mu_-$ and $\mu^+$
determined by the number of fixed points of the Poincar\'{e} maps $T^-_{\mu^+}$ of
$x'=d\,x+\cos(t)+g^-(x)+\mu^+$ and $T^+_{\mu_-}$ of $x'=d\,x+\cos(t)+g^+(x)+\mu_-$,
with $\mu^+$ and $\mu_-$ given by \eqref{ig:4mu}.}
\label{tabla:posibilidades}
\end{table}

To find explicit examples, we take $d=0.1$, so that
$\mu_-\approx0.0995$ and $\mu^+\approx-0.0995$, and define
\begin{equation}\label{def:4g}
 g(x):=\left\{\begin{array}{ll}
 -g_+\,x^3&\text{for }x\le 0\,,\\[.1cm]
 -g_-\,x^3&\text{for }x\ge 0\,,
 \end{array}\right.
\end{equation}
for $g_->0$ and $g_+>0$.
\begin{figure}
     \centering
         \includegraphics[width=\textwidth]{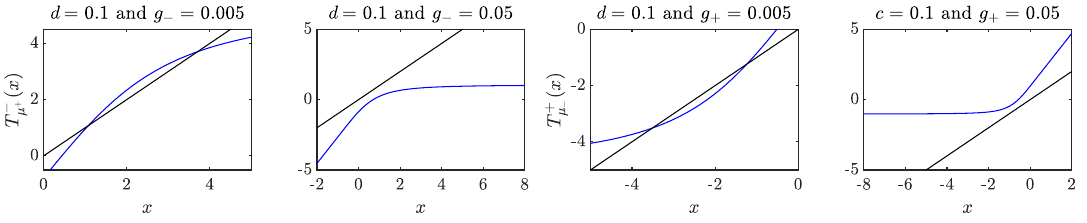}
         \caption{Numerical depiction (in blue) of the Poincar\'{e} maps $T^-_{\mu^+}$ of
         of $x'=d\,x+\cos(t)+\min(g(x),0)+\mu^+$ and
         $T^+_{\mu_-}$ of $x'=d\,x+\cos(t)+\max(g(x),0)+\mu_-$
         for $d=0.1$ and different values of $g_-$ and $g_+$ in \eqref{def:4g}.
         In black, the diagonal $y=x$. The two panels on the left
         (resp. on the right) represent two different situations for $T^-_{\mu^+}$
         (resp. $T^+_{\mu_-}$): two fixed points in the first panel with $g_-=0.005$
         (resp. in the third panel with $g_+=0.005$) and no fixed points in the
         second panel with $g_-=0.05$ (resp. in the fourth panel with $g_+=0.05$).
         The concavity of $g^-$ and the convexity of $g^+$ ensure that $T_\lb^-$ is
         concave and $T_\lb^+$ is convex for all $\lb$. So, clearly, there are no fixed points
         the representation windows.
         The Matlab2023a \texttt{ode45} algorithm was used
         with \texttt{AbsTol} and \texttt{RelTol} equal to \texttt{1e-12} for the numerical integration.
}
        \label{fig:poincare}
\end{figure}
\begin{figure}
     \centering
         \includegraphics[width=\textwidth]{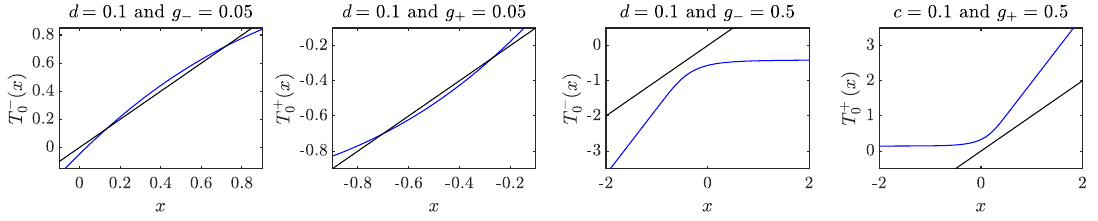}
         \caption{Numerical depiction (in blue) of the Poincar\'{e} maps $T^-_0$
         of $x'=d\,x+\cos(t)+\min(g(x),0)$ and $T^+_0$ of $x'=d\,x+\cos(t)+\max(g(x),0)$ for $d=0.1$ and different values
         of $g_-$ and $g_+$ in \eqref{def:4g}.
         The methods of calculation are those explained in the caption of Figure~\ref{fig:poincare}.}
        \label{fig:poincare2}
\end{figure}
\begin{figure}
     \centering
         \includegraphics[width=\textwidth]{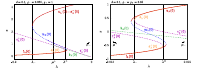}
         \caption{Main dynamical elements of \eqref{eq:4casos} for different values of the parameter.
         In the left-hand panel, $c=0.1$, and $g_-=0.005$ and $g_+=1$ (in \eqref{def:4g}). In Figure~\ref{fig:poincare},
         it has been checked that $T_{\mu^+}^-$ has two fixed points and in a similar way,
         it can be checked that $T_{\mu_-}^+$ has no fixed points.
         Therefore, as explained in Table~\ref{tabla:posibilidades}, \ref{o3} holds.
         In the right-hand panel, $d=0.1$, $g_-=g_+=0.05$.
         Figures~\ref{fig:poincare} and \ref{fig:poincare2} ensure that \ref{o4} holds,
         as explained in the last lines of Section \ref{subsec:4todos}.
         The lines depicted in the pictures represent the sections at $0\in\s^1$ of the
         $\tau_\lb$, $\tau_\lb^-$ and $\tau_\lb^+$ equilibria of the extension to the hull of
         \eqref{eq:4casos}$^\lb$ that have been described in this section and in
         the previous one.
         In particular, these sections correspond to the fixed points of the associated Poincar\'{e}
         maps. That is: $\muk_\lb(0)$, $\mm_\lb(0)$ and $\ml_\lb(0)$ are the fixed points of
         $x\mapsto x_\lb(2\pi,0,x)$; $\mr_\lb^+(0)$ and $\ma_\lb^+(0)$ are the fixed points of
         $x\mapsto x_\lb^+(2\pi,0,x)$; $\mr_\lb^-(0)$ and $\ma_\lb^-(0)$ are the fixed points of
         $x\mapsto x_\lb^-(2\pi,0,x)$; and $\mb_\lb^-(0)$ is the fixed point of
         $x\mapsto x_\lb^l(2\pi,0,x)$.}
         \label{fig:diagrams}
\end{figure}
Figure~\ref{fig:poincare} depicts several numerical approximations for the corresponding
Poincar\'{e} maps $T^-_{\mu^+}$ and $T^+_{\mu_-}$. It shows that different choices of $g_-$ and
$g_+$ give rise to the four possibilities described in Table~\ref{tabla:posibilidades}:
$g_-=g_+=0.005$ gives case \ref{o1}, $g_-=0.05$ and $g_+=0.005$ gives case \ref{o2},
$g_-=0.005$ and $g_+=0.05$ gives case \ref{o3} and $g_-=g_+=0.05$ gives either case \ref{o4} or \ref{o5}.
The numerical study in Figure~\ref{fig:poincare2} shows that, in the case of $g_-=g_+=0.05$,
both $T^-_0$ and $T^+_0$ have two fixed points.
This ensures that $\lb_-<0<\lb^+$, so $\mu^+<\lb_-<\lb^+<\mu_-$, that is, \ref{o4} takes place.
Figure~\ref{fig:poincare2} also shows that, if $g_-=g_+=0.5$, then both $T^-_0$ and $T^+_0$ have no fixed points,
ensuring that $\lb^+<0<\lb_-$ and hence $\mu^+<\lb^+<\lb_-<\mu_-$, that is, \ref{o5} holds. Finally,
Figure~\ref{fig:diagrams} represents a numerical approximation of the main dynamical elements of
\eqref{eq:4casos} in two different cases. The left-hand one corresponds to \ref{o3} while
the right-hand one corresponds to \ref{o4}.
%%%%%%%%%%%%%%%%%%%%%%%%%%%%%%%%%%%%%%%%%%%%%%%%%%%%%%%%%%%%%%%%%%%%
%%%%%%%%%%%%%%%%%%%%%%%%%%%%%%%%%%%%%%%%%%%%%%%%%%%%%%%%%%%%%%%%%%%%
\subsection{An example of case \ref{o5}}
Let us describe analytically an example of case \ref{o5}.
For each $\ep\in[0,1]$, we
%%\begin{equation}\label{4.gep}
%\[
% g\colon\R\to\R\,\quad x\mapsto\left\{\begin{array}{cl}
%  \quad x^2 &\text{ if }x< 0\,,\\
%  -x^2&\text{ if }x>0\,,
%  \end{array}\right.
%\]
%%\end{equation}
consider the one-parametric family of
$2\pi$-periodic ODEs
\[
 x'=\ep\,x+\cos t-x\,|x|+\lb
\]
for $\lb\in\R$. The hull procedure provides the compact set $\W=\s^1\equiv\R/(2\pi\Z)$,
the minimal and uniquely ergodic flow $\sigma_t(s)=s+t$, and, for each $\lb\in\R$, the family
\begin{equation}\label{eq:4lbep}
 x'=\ep\,x+\cos(s+t)-x\,|x|+\lb\,,\quad s\in\s^1\,,
\end{equation}
which determines the skewproduct flow $\tau_\ep^\lb$. In this case, the function $d$
is constant: $d\equiv\ep$, and hence $\ep=\int_\W d\,dm$
for the unique measure $m\in\merg$. It is easy to check that \ref{cc1}-\ref{cc6} hold if $\ep>0$,
and hence Proposition \ref{prop:4tapas}(iii) ensures that, in the situation described in the
next proposition, we are case \ref{o5} (i.e., $\lb_-(\ep)\ge\lb^+(\ep)$) whenever $\ep\in(0,\ep_0)$:
otherwise the upper and lower bounded equilibria would be different in
the interval $(\lb_-(\ep),\lb^+(\ep))$.
\begin{prop}
There exists $\ep_0\in(0,1]$ such that, if $\ep\in[0,\ep_0)$, then:
the attractor $\mA_\ep^\lb$ of \eqref{eq:4lbep}$^\lb_\ep$ is an attractive hyperbolic $\tau_\ep^\lb$-copy
of the base for all $\lb\in\R$, $\mA^\lb_\ep=\{\ml^\lb_\ep\}$; the map
$\R\to C(\s^1,\R),\;\lb\mapsto \ml^\lb_\ep$ is strictly increasing; and
$\lim_{\lb\to\pm\infty}\ml_\ep^\lb(s)=\pm\infty$ uniformly in $\s^1$.
\end{prop}
\begin{proof}
We begin by analyzing the situation for $d\equiv\ep=0$.
It is clear that $x\equiv 0$ does not solve $x'=\cos(s+t)-x\,|x|+\lb$ at any real interval.
Since $(-x|x|)'=-2\,|x|<0$ for $x\ne 0$, we conclude that all the Lyapunov exponents of any
minimal set are strictly negative. Hence, \cite[Theorem 3.4]{cano} ensures that
$\mA^0_\lb$ is an attractive hyperbolic $\tau_\lb$-copy of the base for any $\lb\in\R$:
$\mA^0_\lb=\{\ml^0_\lb\}$.
We will prove that the global dynamics is
similar to that of the previous paragraph if $\ep>0$ is small enough.
Once this is checked, the stated monotonicity and limiting properties with respect to $\lambda$ 
for each one of these values of $\ep$ follow from \cite[Proposition 4.2]{dno5}.
\par
Let us first check that there exists $\lb_1>0$ large enough to ensure that
$\mA_\ep^\lb$ is an attractive hyperbolic $\tau_\ep^\lb$-copy of the base for all
$\ep\in[0,1]$ whenever $\lb\ge \lb_1$. We call $f^\ep_\lb(s,x)=\ep\,x+\cos(s)-x\,|x|+\lb$,
so that $(f^\ep_\lb)_x(s,x)=\ep-2\,x$ if $x\ge 0$. If $x\ge 1$, then  $(f^\ep_\lb)_x(s,x)\le-1$.
Therefore, it is enough to show that $\mA_\ep^\lb\subset\s^1\times[1,\infty)$,
since in this case the upper Lyapunov exponent of every minimal set is
negative and we can apply again \cite[Theorem 3.4]{cano}.
It is easy to check that, if $x\le (-1-\sqrt{5})/2\le -1\le-\ep$ and $\lb>0$,
then $f^\ep_\lb(s,x)=\ep\,x+\cos(s)+x^2+\lb>x^2+x-1+\lb\ge \lb$. Therefore,
any solution is above $(-1-\sqrt{5})/2$ from a certain time onwards,
and this holds for any $\lb>0$.
Now we choose $\lb_1$ large enough to ensure that $f^\ep_{\lb_1}(s,x)>1$
whenever $x\in[(-1-\sqrt{5})/2,1]$, and note that this property
persists for $\lb\ge\lb_1$. The conclusion is that any forward
orbit enters $\s^1\times[1,\infty)$ if $\lb\ge\lb_1$, which implies our assertion about
$\mA_\ep^\lb$.
An analogous argument provides $\lb_2<0$ such that $\mA_\ep^\lb$
is an attractive hyperbolic $\tau_\ep^\lb$-copy of the base whenever $\lb\le\lb_2$ and for all $\ep\in[0,1]$.
\par
From now on, we work for $\lb\in[\lb_2,\lb_1]$. In order to subsequently apply
Arzela-Ascoli's theorem, we first prove that
the existence of $r_0>0$ large enough to ensure that
$\mA_\ep^\lb\subset\s^1\times[-r_0,r_0]$
whenever $\ep\in[0,1]$ and $\lb\in[\lb_2,\lb_1]$.
We take $r_0>0$ such that $\ep\,r+\cos(s)-r^2+\lb<-1$
whenever $\ep\in[0,1]$, $s\in[0,2\pi]$, $\lb\in[\lb_2,\lb_1]$ and $r\ge r_0$.
It is easy to check that the solution with initial datum
$r_0$ is unbounded as $t$ decreases. Similarly, if
$-\ep\,r+\cos(s)+r^2+\lb>1$ whenever $\ep\in[0,1]$,
$\lb\in[\lb_2,\lb_1]$ and $r\le-r_0$, then
the solution with initial datum $-r_0$ is unbounded as $t$ increases.
This proves the asserted property.
\par
It follows from the behavior for $\ep=0$ that, for any $\bar\lb\in[\lb_2,\lb_1]$,
there exists $\rho_{\bar\lb}>0$ such that the flow given by \eqref{eq:4lbep}$^{\bar\lb}$
has an attractive hyperbolic copy of the base whenever
$\lb\in[\bar\lb-\rho_{\bar\lb},\bar\lb+\rho_{\bar\lb}]$ and
$\ep\in[0,\rho_{\bar\lb}]$: see Theorem~\ref{th:2pers}.
We look for $\lb^1,\ldots,\lb^m$ such that $[\lb_2,\lb_1]\subset
[\lb^1-\rho_{\lb^1},\lb^1+\rho_{\lb^1}]\cup\cdots\cup
[\lb^m-\rho_{\lb^m},\lb^m+\rho_{\lb^m}]$ and define
$\wit\ep:=\min(\rho_{\lb^1},\ldots,\rho_{\lb^m})$. Therefore,
for any $(\lb,\ep)\in[\lb_2,\lb_1]\times[0,\wit\ep]$,
there exists an attractive hyperbolic $\tau^\lb_\ep$-copy of the base, say
$\{\mn^\lb_\ep\}$. We will prove that there exists $\ep_0\in(0,\wit\ep]$
such that $\mA^\lb_\ep=\{\mn^\lb_\ep\}$ for all $\ep\in[0,\ep_0]$
(which implies $\mn^\lb_\ep=\ml^\lb_\ep$). Assume for contradiction
the existence of sequences $(\ep_n)\downarrow 0$ and $(\lb_n)$ in $[\lb_2,\lb_1]$
such that $\mA^{\lb_n}_{\ep_n}\ne\{\mn^{\lb_n}_{\ep_n}\}$, and
observe that there is no restriction in assuming the existence of
$\lim_{n\to\infty}\lb_n=:\bar\lb\in[\lb_2,\lb_1]$. According to \cite[Theorem 3.4]{cano}, there exists
a $\tau_{\ep_n}^{\lb_n}$-minimal set $\mM^{\lb_n}_{\ep_n}$ with
nonnegative Lyapunov exponent. As said in Section \ref{subsec:4todos},
in the scalar periodic case, any minimal
set is a copy of the base: $\mM^{\lb_n}_{\ep_n}=\{\mm^{\lb_n}_{\ep_n}\}$.
Hence, $ 1/(2\pi)\int_0^{2\pi}\big(\ep_n+g'(\mm^{\lb_n}_{\ep_n}(s))\big)\,ds\ge 0$.
Arzel\'{a}-Ascoli's theorem ensures the existence of subsequences $(\ep_j)$ and
$(\lb_j)$ and of a bounded continuous function $\mm\colon\W\to\R$
such that $\lim_{j\to\infty}\mm^{\lb_j}_{\ep_j}=\mm$ uniformly in $\s^1$.
It is easy to check that $\mm$ is a bounded $\tau_0^{\bar\lb}$-equilibrium,
which as seen in the first paragraph of this proof means that
$\mm=\ml_0^{\bar\lb}$. We conclude from Lebesgue's Dominated Convergence Theorem
that $1/(2\pi)\int_0^{2\pi}g'(\ml_0^{\bar\lb}(s))\,ds\ge 0$,
which contradicts the previous analysis for $\ep=0$.
The proof is complete.
\end{proof}
%%%%%%%%%%%%%%%%%%%%%%%%%%%%%%%%%%%%%%%%%%%%%%%%%%%%%%%%%%%%%%%%%%%%%%%
%%%%%%%%%%%%%%%%%%%%%%%%%%%%%%%%%%%%%%%%%%%%%%%%%%%%%%%%%%%%%%%%%%%%%%%
%%%%%%%%%%%%%%%%%%%%%%%%%%%%%%%%%%%%%%%%%%%%%%%%%%%%%%%%%%%%%%%%%%%%%%%
%%%%%%%%%%%%%%%%%%%%%%%%%%%%%%%%%%%%%%%%%%%%%%%%%%%%%%%%%%%%%%%%%%%%%%%
\section{Critical transitions in concave-convex equations}\label{sec:5}
The results of this section show the existence of size-induced critical transitions
in parametric perturbations of nonautonomous concave-convex ODEs, of the type
\begin{equation}\label{eq:5ct}
 x'=h(t,x)-\rho\,k(t,x)\,,
\end{equation}
where $h$ and $k$ satisfy certain conditions, below described.
As explained in the Introduction,
a novel aspect of the approach to critical transitions in this paper is the lack of
assumptions about future dynamics, which were customary in our previous papers on critical transitions.
More precisely, the hypotheses on $h$ and $k$ will ensure the existence of: a {\em past equation}
for \eqref{eq:5ct}$^\rho$ for $\rho\ge 0$,
with three uniformly separated hyperbolic solutions determining its global dynamics; and
two locally pullback attractive solutions of \eqref{eq:5ct}$^\rho$, whose behavior
as time increases will determine the dynamical situation of {\em tracking\/} or {\em tipping}.
As usual, a {\em critical transition\/} will occur when a small variation in
$\rho$ means a jump from one of these two dynamical scenarios to the other.
The required conditions on $h$ and $k$ are:
\begin{enumerate}[leftmargin=23pt,label=\rm{\bf{ct\arabic*}}]
\item\label{ct1} $h\colon\RR\to\R$ is $C^2$-admissible and, if $\W_h$ is its hull with time-shift
flow $\sigma_h(t,\w_1)=\w_1{\cdot}t$, if $\mh\colon\W_h\times\R\to\R\,,(\w_1,x)\mapsto\w_1(0,x)$,
and if the skew-product flow given by the family $\{x'=\mh(\w_1{\cdot}t,x),|\;\w_1\in\W_h\}$
is $\tau(t,\w_1,x):=(\w_1{\cdot}t,v_0(t,\w_1,x))$,
then:
\begin{itemize}[leftmargin=12pt]
\item[-] the maps $\mc(\w_1):=\mh(\w_1,0)$, $\md(\w_1):=\mh_x(\w_1,0)$ and
$\mg(\w_1,x):=\mh(\w_1,x)-\mh(\w_1,0)-\mh_x(\w_1,0)\,x$ satisfy conditions \ref{cc1}-\ref{cc6};
\item[-] there are exactly three hyperbolic $\tau$-copies of the base $\{\ml_0\}$, $\{\mm_0\}$ and $\{\muk_0\}$,
with $\ml_0<\mm_0<\muk_0$, where $\ml_0$ and $\muk_0$ are attractive and the bounds of the global attractor,
and $\mm_0$ is repulsive;
\item[-] $\lim_{t\to\infty}(v_0(t,\w_1,x)-\ml_0(\w_1{\cdot}t))=0$ for all $x<\mm_0(\w_1)$,
$\lim_{t\to\infty}(v_0(t,\w_1,x)-\muk_0(\w_1{\cdot}t))=0$ for all $x>\mm_0(\w_1)$, and
$\lim_{t\to-\infty}(v_0(t,\w_1,x)-\mm_0(\w_1{\cdot}t))=0$ for all $x\in(\ml_0(\w_1),\muk_0(\w_1))$.
\end{itemize}
\item\label{ct2} $k\colon\RR\to\R$ is $C^1$-admissible, with $k(t,x)\ge 0$ for all
$(t,x)\in\RR$, $\lim_{x\to-\infty}k(t,x)=0$ uniformly on $\R$, and
$\lim_{t\to-\infty}k(t,x)=0$ uniformly on compact subsets of $\R$.
\item\label{ct3} $\ml_0(\w_1)\le 0$
for all $\w_1\in\W_h$, and $k_x(t,x)\ge 0$ whenever %$x\le\sup_{\w_1\in\W_h}\ml_0(\w_1)$.
whenever $x\leq\ml_0(h{\cdot}t)$.
\end{enumerate}
\begin{notas}
1.~Note that, if the family $x'=\mh(\w_1{\cdot}t,x)+\lb,\,\;\w_1\in\W_h$ is in case \ref{o1} (with
$\mc(\w_1)=\mh(\w_1,0)$, $\md(\w_1)=\mh_x(\w_1,0)$ and $\mg(\w_1,x)=\mh(\w_1,x)-\mh(\w_1,0)-\mh_x(\w_1,0)\,x$), and if
$\lb=0$ belongs to the set $(\lb_-,\mu^+)\cup(\mu_-,\lb^+)$ of Theorem \ref{th:4CCbifur}, then
\ref{ct1} holds. This is a
sufficient but not necessary condition. For instance, also under the conditions of
Theorem \ref{th:4CCbifur}, \ref{ct1} holds
whenever $\mm_0^+=\mm_0^-$ provides a repulsive hyperbolic copy of the base.
\par
2.~It is possible to establish sufficient conditions on $h$ guaranteing \ref{ct1}.
If $h$ is $C^1$, then \ref{cc1} holds. If $h_x(t,0)>\delta>0$ for all $t\in\R$,
then \ref{cc2} holds. If $\lim_{x\to\pm\infty} h(t,x)=\mp\infty$ uniformly on $\R$, then
\ref{cc3} holds. If $h$ is $C^2$-admissible and $g(t,x):=h(t,x)-h(t,0)-h_x(t,0)\,x$ is strictly
decreasing, and it satisfies $\sup_{t\in\R}\big(g_x(t,x_2)-g_x(t,x_1)\big)<0$
whenever $0<x_1<x_2$ and whenever $0>x_1>x_2$, then \ref{cc4}, \ref{cc5} and \ref{cc6} hold.
The proof of the sufficiency of these conditions, some of which are far away from necessary, 
can be an exercise for the reader.
\end{notas}

Let us assume \ref{ct1} and \ref{ct2}, and define
$\W:=\mathrm{closure}\{(h{\cdot}s,k{\cdot}s)\mid\,s\in\R\}$ in the compact-open topology
of $C(\R,\R^2)$ and $\sigma(t,\w):=\wt$ as the time-shift flow in $\W$. Note that any element of $\W$
is of the form $\w=(\w_1,\w_2)$, with $\w_1\in\W_h$. The maps $\bar\mh\colon\WR\to\R,\,(\w,x)
\mapsto\w_1(0,x)=\mh(\w_1,x)$ and $\bar\mk\colon\WR\to\R,\,(\w,x)\mapsto\w_2(0,x)$ are the extensions
of $h$ and $k$ to their common hull $\W$, and they respectively belong to $C^{0,2}(\WR,\R)$
and $C^{0,1}(\WR,\R)$ (with $\bar\mh_x(\w,x)=(\w_1)_x(0,x)=\mh_x(\w_1,x),\,\bar\mh_{xx}(\w,x)=
(\w_1)_{xx}(0,x)=\mh_{xx}(\w_1,x)$, and $\bar\mk_x(\w,x)=(\w_2)_x(0,x)$ ).
So, the $\rho$-parametric family \eqref{eq:5ct} is included in the $\rho$-parametric family
\begin{equation}\label{eq:5cthull}
 x'=\bar\mh(\wt,x)-\rho\,\bar\mk(\wt,x)\,,\qquad\w\in\W
\end{equation}
of families of equations over $\W$: more precisely, equation \eqref{eq:5ct}$^\rho$ coincides with equation
\eqref{eq:5cthull}$^\rho_{\w_0}$ for $\w_0:=(h,k)\in\W$.
Let $\bar\tau_\rho(t,\w,x):=(\wt,\bar v_\rho(t,\w,x))$ be the flow induced on $\WR$ by
\eqref{eq:5cthull}$^\rho$ (see Section \ref{subsec:2skew}). It is clear that
\begin{itemize}[leftmargin=12pt]
\item[1.] the maps
$\bar\ml_0(\w):=\ml_0(\w_1)$, $\bar\mm_0(\w):=\mm_0(\w_1)$ and $\bar\muk_0(\w):=\muk_0(\w_1)$
(for $\w=(\w_1,\w_2)$) define three hyperbolic $\bar\tau_0$-copies of the base, attractive, repulsive,
and attractive, since $\bar v_0(t,\w,x)=v_0(t,\w_1,x)$.
\item[2.] These equalities and one of the last hypotheses included in \ref{ct1} guarantee that
$\lim_{t\to\infty}(\bar v_0(t,\w,x)-\bar\ml_0(\wt))=0$
whenever $x<\bar\mm_0(\w)$.
\item[3.] Condition \ref{ct2} yields $\bar\mk(\w,x)\ge 0$ for all $(\w,x)\in\WR$,
$\lim_{x\to-\infty}\bar\mk(\w,x)=0$ uniformly on $\W$, and
$\lim_{t\to-\infty}\bar\mk(\bwt,x)=0$ for $\bar\w=(h,k)\in\W$ uniformly on compact subsets of $\R$.
If also \ref{ct3} holds, then $\bar\mk_x(\bwt,x)\ge 0$ whenever
$x<\bar\ml_0(\bwt)$, which follows from $\bar\mk_x(\bwt,x)=k_x(t,x)$ and from $\bar\ml_0(\bwt)=\ml_0(h{\cdot}t)$.
The definition of $\W$ and the continuity of $\bar\ml_0$ make it easy to deduce that
$\bar\mk_x(\w,x)\ge 0$ if $\w\in\W$ and $x<\bar\ml_0(\w)$, and
the continuity of $\bar\mk_x$ shows that $\bar\mk_x(\w,x)\ge 0$ if $\w\in\W$ and $x\le\bar\ml_0(\w)$.
\item[4.] If $\bar m\in\merg$, then its projection
$m(\mB):=\bar m(\{(\w_1,\w_2)\in\W\,|\;\w_1\in\mB\})$
for any Borel set $\mB\subset\W_h$ defines a measure $m\in\mathfrak{M}_\mathrm{erg}(\W_h,\sigma_h)$.
\end{itemize}
The map $\bar\mh_\rho:=\bar\mh-\rho\,\bar\mk$ belongs to $C^{0,1}(\WR,\R)$ for any $\rho\in\R$.
In addition, if $\rho\ge 0$, then $\limsup_{x\to\pm\infty}(\pm\bar\mh_\rho(\w,x))<0$ uniformly on $\W$:
this follows from condition \ref{cc3}, included in \ref{ct1}, which ensures that
$\limsup_{x\to\pm\infty}(\pm\bar\mh(\w,x))<0$ uniformly on $\W$, combined with $\bar\mk(\w,x)\ge 0$ for the
limit at $+\infty$ and with $\lim_{x\to-\infty}\bar\mk(\w,x)=0$ uniformly on $\W$ (see point 3 above)
for the limit at $-\infty$. As \cite[Proposition 5.5]{dno4} shows, these two properties ensure the
existence of the global attractor $\mA_\rho$ for $\bar\tau_\rho$,
\[
 \mA_\rho=\{(\w,x)\,|\;\w\in\W\text{ and }\bar\ml_\rho(\w)\le x\le\bar\muk_\rho(\w)\}\,,
\]
where $\bar\ml_\rho\colon\W\to\R$ and $\bar\muk_\rho\colon\W\to\R$ are $\bar\tau_{\rho}$-equilibria
respectively lower and upper semicontinuous, with $\rho\mapsto\bar\ml_\rho(\w)$ nonincreasing
and right-continuous and with $\rho\mapsto\bar\muk_\rho(\w)$ nonincreasing and left-continuous
on $\R$ for all $\w\in\W$ (see also Remark \ref{rm:4existeA}).

Let us fix $\bar\w=(h,k)\in\W$. The condition $\lim_{t\to-\infty}\bar\mk(\bwt,x)=0$ allows
us to understand $x'=\bar\mh(\bwt,x)$ (i.e., $x'=h(t,x)$)
as the {\em past equation} of $x'=\bar\mh_\rho(\bwt,x)$ (i.e.,
of \eqref{eq:5ct}$^\rho$). The lower and upper solutions
$t\mapsto\bar\ml_\rho(\bwt)$ and $t\mapsto\bar\muk_\rho(\bwt)$ of
\eqref{eq:5cthull}$^\rho_{\bar\w}=\;$\eqref{eq:5ct}$^\rho$
are characterized in a different way in our next result.
\begin{prop}\label{prop:5pull}
Assume that \ref{ct1} and \ref{ct2} hold, and fix $\bar\w=(h,k)\in\W$ and $\rho\ge 0$. Then,
\begin{itemize}[leftmargin=22pt]
\item[\rm(i)] the (global) solutions $t\mapsto\bar\ml_\rho(\bwt)$ and $t\mapsto\bar\muk_\rho(\bwt)$
of \eqref{eq:5cthull}$^\rho_{\bar\w}$ are the unique ones satisfying
$\lim_{t\to-\infty}(\bar\ml_\rho(\bwt)-\bar\ml_0(\bwt))=0$ and $\lim_{t\to-\infty}(\bar\muk_\rho(\bwt)-\bar\muk_0(\bwt))=0$,
and they are locally pullback attractive.
\item[\rm(ii)] If, in addition, \ref{ct3} holds, then
$\{\bar\ml_\rho\}$ is an attractive hyperbolic $\bar\tau_\rho$-copy of the base for all $\rho\ge0$,
$\lim_{t\to\infty}(\bar v_\rho(t,\w,x)-\bar\ml_\rho(\wt))=0$ whenever $x\le\bar\ml_0(\w)$, and
$[0,\infty)\to C(\WR,\R),\,\rho\mapsto\bar\ml_\rho$ is continuous.
\end{itemize}
\end{prop}
\begin{proof} To simplify the notation, we call
%$l_\rho(t):=\bar\ml_\rho(\bwt)$, $u_\rho(t):=\bar\muk_\rho(\bwt)$, and
$\bar h_\rho(t,x):=\bar\mh_\rho(\bwt,x)=
\mh(\bar\w_1{\cdot}t,x)-\rho\,\bar\mk(\bwt,x)=h(t,x)-\rho\,k(t,x)$.

(i)~This proof adapts several arguments from \cite{dno3,dno4}.
As a first step, reasoning as in the proofs of \cite[Proposition 6.2(i)]{dno4} and
\cite[Proposition 3.5(i)]{dno3}, we deduce the next property from the dynamics of
$x'=\bar h_0(t,x)$ (determined by \ref{ct1}): if $f$ is $C^2$-admissible and with $\limsup_{x\to\pm\infty}(\pm f(t,x))<0$,
if $\lim_{t\to-\infty}(f(t,x)-\bar h_0(t,x))=0$
uniformly in the compact sets of $\mK$, and if $x'=f(t,x)$ has three uniformly separated hyperbolic solutions $l_f<m_f<u_f$
with $l_f$ and $u_f$ attractive and $m_f$ repulsive, then
$l_f$ and $u_f$ are the lower and upper bounded solutions of $x'=f(t,x)$, and they
approach the attractive hyperbolic solutions of $x'=\bar h_0(t,x)$ as time decreases.
(A more self-contained proof is given in \cite[Proposition 4.12(ii)]{duenasphdthesis}.)

Now, we use those arguments of the proofs of \cite[Proposition 6.3(i)]{dno4} and \cite[Theorem 3.7]{dno3}
regarding the behavior as $t\to-\infty$. The main idea is to suitably approximate $x'=\bar h_0(t,x)$
by an equation $x'=f(t,x)$ with the same past equation (which, as explained above, is $x'=\bar h_0(t,x)$),
in such a way that $x'=f(t,x)$ fulfills the conditions of the previous paragraph. The fact that
the attractive hyperbolic solutions of $x'=f(t,x)$ are also the lower and upper bounded solutions
combined with the fact that $t\mapsto\bar\ml_\rho(\bwt)$ and $t\mapsto\bar\muk_\rho(\bwt)$
are the lower and upper bounded solutions of $x'=\bar h_\rho(t,x)$ allows us to repeat the arguments of
those proofs, in order to check all the assertions in (i).
(A more self-contained proof is given in part of the proof of \cite[Theorem 4.13(i)]{duenasphdthesis}.)
\smallskip

(ii) It follows from $\bar\mk(\w,\bar\ml_0(\w))\ge 0$ (see point 3 after \eqref{eq:5cthull})
that the upper semicontinuous map $\bar\ml_0$ is a global upper solution for
$\bar\tau_\rho$ for all $\rho>0$, with $\bar\ml_0\ge\bar\ml_\rho$. Hence, the map
$\bar\mb_\rho(\w):=\lim_{s\to\infty}\bar v_\rho(s,\w{\cdot}(-s),\bar\ml_0(\w{\cdot}(-s)))$
is an upper semicontinuous $\bar\tau_\rho$-equilibrium with $\bar\ml_\rho\le\bar\mb_\rho\le\bar\ml_0$
(see, e.g., \cite[Section 2.2]{dno1} and references therein).
Let us define $\mL_\rho:=\{(\w,x)\in\WR\,|\;\bar\ml_\rho(\w)\le x\le\bar\mb_\rho(\w)\}$, which is a compact
$\bar\tau_\rho$-invariant set due to the semicontinuity properties of the defining $\bar\tau_\rho$-equilibria.
Note that $x\le\bar\ml_0(\w)=\ml_0(\w_1)$ for all $(\w,x)\in\mL_\rho$ with $\w=(\w_1,\w_2)$. In addition,
$(\bar\mh_\rho)_x(\w,x)=\mh_x(\w_1,x)-\rho\,\bar\mk_x(\w,x)\le \mh_x(\w_1,x)$ for all $(\w,x)\in\WR$.
We take a measure $\bar m\in\merg$ and its projection $m\in\mathfrak{M}_{\text{\rm erg}}(\W_h,\sigma_h)$ (see property
4 of the list before the statement). Note that condition \ref{cc5}, included in \ref{ct1}, ensures that
$\bar m(\{\w\in\W\,|\;x\mapsto \bar\mh_x(\w,x)$ is nondecreasing on $(-\infty,0]\})=
m(\{\w_1\in\W_h\,|\;x\mapsto \mh_x(\w_1,x)$ is nondecreasing on $(-\infty,0]\})=1$. In addition,
the attractive hyperbolicity of $\ml_0$, also included in \ref{ct1}, ensures that
$\int_\W\mh_x(\w_1,\ml_0(\w_1))\,dm<0$ (see Theorem \ref{th:2copia}). Let
$\mb\colon\W\to\R$ be an $\bar m$-measurable $\bar\tau_\rho$-equilibrium
with graph contained in $\mL_\rho$. Then,
$\int_\W (\bar\mh_\rho)_x(\w,\mb(\w))\,d\bar m\le\int_\W\bar\mh_x(\w,\mb(\w))\,d\bar m=\int_\W\mh_x(\w_1,\mb(\w))\,dm\le
\int_\W\mh_x(\w_1,\ml_0(\w_1))\,dm<0$. This means that the upper
Lyapunov exponent of $\mL_\rho$ is negative. From here we can reason as in the proof of
Proposition \ref{prop:4tapas}(i) to check that
$\mL_\rho$ is an attractive hyperbolic copy of the base,
with $\mL_\rho=\{\bar\ml_\rho\}=\{\bar\mb_\rho\}$. This proves
the first assertion in (ii).

Since $\bar\ml_{\rho}$ is the lower equilibrium of the global attractor,
$\lim_{t\to\infty}(\bar v_\rho(t,\w,x)-\bar\ml_\rho(\wt))=0$
whenever $x\le\bar\ml_\rho(\w)$. Let us now check the same for $x=\bar\ml_0(\w)$
and hence, due to the fiber-monotonicity of the flow,
for all $x\in(\bar\mb_\rho(\w),\bar\ml_0(\w)]=(\bar\ml_\rho(\w),\bar\ml_0(\w)]$.
We take $\bar\w\in\W$, take a minimal set $\mM$ contained in
the omega-limit set of $(\bar\w,\bar\ml_0(\bar\w))$ for $\bar\tau_\rho$
and $\W_\mM$ to its projection on $\W$ (also minimal). Since $\bar\ml_0$ is a global upper
solution for $\bar\tau_\rho$, then $x\le\bar\ml_0(\w)$ for all $(\w,x)\in\mM$.
The upper boundary of $\mM$, $\bar\muk_\mM$, is an equilibrium for the
restriction of $\bar\tau_\rho$ to $\W_\mM\times\R$, and hence $\bar\mb_\rho(\w)=\bar\ml_\rho(\w)\le
\bar\muk_\mM(\w)\le\bar\ml_0(\w)$ for all $\w\in\W_\mM$. The definition of
$\bar\mb_\rho$ precludes the possibility $\bar\mb_\rho(\w_0)<\bar\muk_\mM(\w_0)$ for any $\w_0\in\W_\mM$:
if this were the case, then $\bar\muk_\mM(\w_0)=\bar v_\rho(s,\w_0{\cdot}(-s),
\bar\ml_0(\w_0{\cdot}(-s)))$ for an $s>0$; this and the minimality of $\mM$
would imply that the map $\bar\mb_s\colon\W_\mM\to\R$
given by $\bar\mb_s(\w):=\bar v_\rho(s,\w{\cdot}(-s),
\bar\ml_0(\w{\cdot}(-s)))$ defines a $\bar\tau_\rho|_{\W_\mM\times\R}$-equilibrium which
coincides with $\bar\mb_\rho$ and with
$\bar\muk_\mM$; and this is contradictory. Hence, $(\w,\bar\ml_\rho(\w))=
(\w,\bar\mb_\rho(\w))\in\mM$ for all $\w\in\W_\mM$.
The limiting behaviour stated in (ii) follows easily.
Finally, the persistence of the existence of a hyperbolic copy of the base
for small variations of $\rho$ (see Theorem \ref{th:2pers}), combined with the
monotonicity and right-continuity of $\rho\mapsto\bar\ml_{\rho}(\w)$ and with the
fact that $\bar v_\rho(t,\w,x)$ is unbounded if $x<\bar\ml_\rho(\w)$, provide an easy proof of the
continuity stated in (ii).
\end{proof}
\begin{defi}\label{def:5trackingtipping}
Under conditions \ref{ct1} and \ref{ct2}, we say that \eqref{eq:5ct}$^\rho$
exhibits {\em tracking} for $\rho\ge 0$ if $\inf_{t\in\R}(u_\rho(t)-l_\rho(t))>0$,
and {\em tipping} otherwise, where $l_\rho(t)$ and $u_\rho(t)$
are the lower and upper bounded solutions of the equation.
\end{defi}
In other words, if \ref{ct1} and \ref{ct2} hold, \eqref{eq:5ct}$^\rho$ exhibits tracking
if and only if it has uniformly separated bounded solutions.

The next result proves the existence of a unique tipping point for the parametric family
\eqref{eq:5ct} (i.e., a unique shift in $\rho$ from tracking to tipping) under some additional assumptions:
the two upper hyperbolic copies of the base of the past equation $x'=\mh(\wt,x)$
are strictly positive and bound an area
at which $\mk$ is strictly positive and convex with respect to $x$.
\begin{teor}\label{th:5ctconcave}
Assume that \ref{ct1}, \ref{ct2} and \ref{ct3} hold, as well as
\begin{itemize}[leftmargin=12pt]
\item[-] $\inf_{\w\in\W}\bar\mm_0(\w)>0$,
\item[-] $\mk(\w,x)$ is strictly positive on $\mB^*:=\bigcup_{\w\in\W}(\{\w\}\times[\bar\mm_0(\w), \bar\muk_0(\w)])$,
    and $x\mapsto\mk(\w,x)$ is convex in $[\bar\mm_0(\w),\bar\muk_0(\w)]$
    for all $\w\in\W$.
\end{itemize}
Then, there exists $\rho_0>0$ such that the equation \eqref{eq:5ct}$^\rho$
exhibits tracking for $\rho\in[0,\rho_0]$ and tipping for $\rho>\rho_0$.
\end{teor}
\begin{proof} Our first goal is to get part of the conclusions of \cite[Theorem 3.4]{dno5} 
for a suitable equation. Let us define $\mh^*\colon\WR\to\R$ by: $\mh^*(\w,x):=\bar\mh(\w,x)$ if $(\w,x)\in\mB^*$;
$\mh^*(\w,x):=\bar\mh(\w,\bar\muk_0(\w))+\bar\mh_x(\w,\bar\muk_0(\w))(x-\bar\muk_0(\w))-
(x-\bar \muk_0(\w))^2$ if $x>\bar\muk_0(\w)$; and $\mh^*(\w,x):=\bar\mh(\w,\bar\mm_0(\w))+
\bar\mh_x(\w,\bar\mm_0(\w))(x-\bar\mm_0(\w))-
(x-\bar\mm_0(\w))^2$ if $x<\bar\mm_0(\w)$.
And we define $\mk^*$ from $\bar\mk$ in a similar way, but changing
$-(x-\bar\muk_0(\w))^2$ and $-(x-\bar\mm_0(\w))^2$ by $(x-\bar\muk_0(\w))^2$ and $(x-\bar\mm_0(\w))^2$.
We take $\bar\w=(h,k)\in\W$. The mentioned equation is
\begin{equation}\label{eq:5ct*}
 x'=\mh^*(\bwt,x)-\rho\,\mk^*(\bwt,x)
\end{equation}
It is easy to check that $\mh^*$ and $\mk^*$ belong to $C^{0,1}(\WR,\R)$,
and hence the map $(\w,x)\mapsto\mf(\w,x,\rho):=\mh^*(\w,x)-\rho\,\mk^*(\w,x)$ belongs to 
$C^{0,1}(\WR,\R)$ for all $\rho\ge 0$.
The expressions of $\mh^*$ and $\mk^*$ for
$x\notin[\bar\mm_0(\w),\bar\muk_0(\w)]$ show that $\lim_{x\to\pm\infty}\mf(\w,x,\rho)=-\infty$
uniformly on $\W$ for all $\rho\ge 0$. It is also easy to deduce from the
definition of $\mh^*$ and $\mk^*$, from the concavity conditions \ref{cc5} and \ref{cc6} on $\W\times[0,\infty)$
(guaranteed by \ref{ct1}),
from the convexity condition on $\mk$ in the statement, and from the
relation between $\merg$ and $\mathfrak{M}_\mathrm{erg}(\W_h,\sigma_h)$
that, if $\rho\ge 0$, then $\bar m(\{\w\in\W\,|\;x\mapsto\mf(\w,x,\rho)$ is concave$\})=1$ and
$\bar m(\{\w\in\W\,|\;\,x\mapsto\mf_x(\w,x,\rho)$ is strictly decreasing on $\mJ\})>0$ for all compact
sets $\mJ\subset\R$ and $\bar m\in\merg$. All these properties show that $\mf$ fulfills the first hypothesis in
\cite[Theorems 3.3 and 3.4]{dno5}. The second one is the joint continuity of $\mf$ and $\mf_x$, obvious.
Finally, since $\mf_\rho(\w,x,\rho)=-\mk^*(\w,x)$,
the condition
$\sup_{(\w,x)\in\mB^*}(-\mk^*(\w,x))<0$ substitutes
the last hypothesis in \cite[Theorems 3.3 and 3.4]{dno5} (which, having in mind that the
perturbation parameter is $-\rho$, is $\sup_{(\w,x,\rho_1,\rho_2)\in\W\times\mJ\times\RR,\,\rho_1<\rho_2}
(\mf(\w,x,\rho_2)-\mf(\w,x,\rho_1))/(\rho_2-\rho_1)<0$ for each compact interval $\mJ$ in
$\R$). These properties allow us to repeat the proof of \cite[Theorem 3.3]{dno5} for
the family of equations $x'=\mh^*(\wt,x)-\rho\,\mk^*(\wt,x)$, $\w\in\W$ just focusing on
$\rho>0$, having in mind that its upper and lower bounded equilibria evolve monotonically
as $\rho$ increases as long as their graphs are contained in $\mB^*$ due to the condition 
$\mk^*(\w,x)=\mk(\w,x)\ge 0$ for $(\w,x)\in\mB^*$. With this done, we can also repeat part 
of the proof of \cite[Theorem 3.4]{dno5},
again for $\rho>0$. Its main conclusions allow us to continue with this proof.

Recall that $l_\rho(t):=\bar\ml_\rho(\bwt)$ and $u_\rho(t):=\bar\muk_\rho(\bwt)$,
are the lower and upper bounded solutions of \eqref{eq:5ct}$^\rho$, and that
$\rho\mapsto l_\rho,\,u_\rho$ are nonincreasing. The previous adaptation of
\cite[Theorem 3.4]{dno5} shows the existence of $\rho_0>0$
such that there are: two hyperbolic solutions $m_\rho$ and $v_\rho$ of
\eqref{eq:5ct*}
for all $\rho\in[0,\rho_0)$, with $l_\rho\le l_0<m_0<m_\rho<v_\rho<u_0$; at least a bounded
and nonhyperbolic solution $v_{\rho_0}$ for $\rho=\rho_0$, with $l_{\rho_0}\le l_0<m_0<v_{\rho_0}<u_0$;
and no bounded solutions with graph in  $\mB^*$ for $\rho>\rho_0$.
In particular, \eqref{eq:5ct}$^\rho$
exhibits tracking for $\rho\in[0,\rho_0]$, since any solution of \eqref{eq:5ct*}$^\rho$ is
also a solution of \eqref{eq:5ct}$^\rho$ as long as its graph lies in $\mB^*$,
and since $l_0$ and $m_0$ are uniformly separated (as $l_0$ and $m_0$): these properties ensure
that the bounded solutions $l_\rho$ and $v_\rho$ of \eqref{eq:5ct}$^\rho$ are uniformly separated.

Let us fix $\rho>\rho_0$. It remains to check that $\inf_{t\ge 0}(u_\rho(t)-l_\rho(t))=0$,
which ensures that \eqref{eq:5ct}$^\rho$ exhibits tipping. We know that $u_\rho(t)\le u_0(t)=\bar\muk_0(\bwt)$ for all
$t\in\R$, and hence the absence of bounded solutions for \eqref{eq:5ct*}$^\rho$
with graph in $\mB^*$ precludes the possibility
$u_\rho(t)\ge\bar\mm_0(\bwt)$ for all $t\in\R$. So, there exists $t_0\in\R$ with $u_\rho(t_0)<\bar\mm_0(\bwt_0)$,
which, as explained in point 2 after \eqref{eq:5cthull}, yields $\lim_{t\to\infty}(x_0(t,t_0,u_\rho(t_0))-l_0(t))=
\lim_{t\to\infty}(\bar v_0(t-t_0,\bwt_0,\bar\muk_\rho(\bwt_0))-\bar\ml_0(\bwt))=0$.
Since $u_\rho(t)=\bar v_\rho(t-t_0,\bwt_0,\bar\muk_\rho(\bwt_0))\le
\bar v_0(t-t_0,\bwt_0,\bar\muk_\rho(\bwt_0))$ for all $t\ge t_0$ (due to $\mk\ge 0$),
we conclude that $\limsup_{t\to\infty}(u_\rho(t)-l_0(t))\le 0$. Recall that
$l_0(t)=\bar\ml_0(\bwt)$. We take a sequence $(t_n)\uparrow\infty$ such that there exists
$(\w_0,x_0):=\lim_{n\to\infty}(\bwt_n,\bar\muk_\rho(\bwt_n))$. Since $\bar\ml_0$ is continuous,
$\bar\ml_0(\w_0)=\lim_{n\to\infty}\bar\ml_0(\bwt_n)$, and the previous property shows that $x_0\le\bar\ml_0(\w_0)$.
According to Proposition \ref{prop:5pull}(ii), $\lim_{t\to\infty}(v_\rho(t,\w_0,x_0)-\bar\ml_\rho(\w_0{\cdot}t))=0$,
from where $\inf_{t\ge 0}(\bar\muk_\rho(\bwt)-\bar\ml_\rho(\bwt))=0$ follows easily. This completes the proof.
\end{proof}
\begin{nota}\label{rm:5barrho}
In Section \ref{subsec:6modelo} we will consider an equation of the type \eqref{eq:5ct} for which all the hypotheses
of Theorem \ref{th:5ctconcave} are satisfied excepting one: the map $x\mapsto\mk(\w,x)$ is not
convex in $[\bar\mm_0(\w),\bar\muk_0(\w)]$. However, in that example,
there exists a value $\bar\rho>0$ such that
the dynamics on $\WR$ is like that required in the theorem for $\rho=0$ (ensured by \ref{ct1}), and such that
the map $x\mapsto\mk(\w,x)$ is convex in $[\bar\mm_{\bar\rho}(\w),\bar\muk_{\bar\rho}(\w)]$.
So, we can repeat the previous proof taking $\bar\mm_{\bar\rho}$ and $\bar\muk_{\bar\rho}$ to define $\mh^*$ and
$\mk^*$, and conclude the existence of a unique critical value $\rho_0>\bar\rho$.
\end{nota}
We conclude this section with a result on critical transitions related to the
d-concavity band of Section \ref{subsec:5band}. The assumptions added to those of Proposition
\ref{prop:5pull} are that the two upper hyperbolic copies of the base of $x'=\mh(\wt,x)$ provided by \ref{ct1} lie within the strict
d-concavity band, and that $-k$ preserves d-concavity in the region where these copies of the base are located.
\begin{teor}\label{th:5ctdconcave}
Assume that \ref{ct1}, \ref{ct2} and \ref{ct3} hold, and that $k$ is $C^2$-admissible, as well as
\begin{itemize}[leftmargin=22pt]
\item[-] $\mm_0$ and $\muk_0$ are jointly included in the strict
d-concavity region $\mD^*_\mg$ given by Definition {\em\ref{def:4strictband}} for the map
$\mg(\w_1,x):=\mh(\w_1,x)-\mh(\w_1,0)-\mh_x(\w_1,0)\,x$.
\item[-] $\mk(\w,x)$ is strictly positive on $\mB^*:=\bigcup_{\w\in\W}(\{\w\}\times[\bar\mm_0(\w),
    \bar\muk_0(\w)])$,
    and $x\mapsto\mk_x(\w,x)$ is convex in $[\bar\mm_0(\w),\bar\muk_0(\w)]$
    for all $\w\in\W$.
\end{itemize}
Then, there exists $\rho_0>0$ such that the equation \eqref{eq:5ct}$^\rho$
exhibits tracking for $\rho\in[0,\rho_0]$ and tipping for $\rho>\rho_0$.
\end{teor}
\begin{proof}
Let $\mh^*\colon\WR\to\R$ be defined by $\mh^*(\w,x):=\bar\mh(\w,x)$ if $(\w,x)\in\mB^*$, $\mh^*(\w,x):=\bar\mh(\w,\bar\muk_0(\w))+\bar\mh_x(\w,\bar\muk_0(\w))(x-\bar\muk_0(\w))+
(1/2)\,\bar\mh_{xx}(\w,\bar\muk_0(\w))(x-\bar\muk_0(\w))^2-(x-\bar\muk_0(\w))^3$ if $x>\bar\muk_0(\w)$, and $\mh^*(\w,x):=\bar\mh(\w,\bar\mm_0(\w))+\bar\mh_x(\w,\bar\mm_0(\w))(x-\bar\mm_0(\w))+
(1/2)\bar\mh_{xx}(\w,\bar\mm_0(\w))(x-\bar\mm_0(\w))^2-(x-\bar\mm_0(\w))^3$ if $x<\bar\mm_0(\w)$.
And we analogously define $\bar\mk$ from $\bar\mk$, but replacing $-(x-\bar\muk_0(\w))^3$ and $-(x-\bar\mm_0(\w))^3$ with $(x-\bar\muk_0(\w))^3$ and $(x-\bar\mm_0(\w))^3$.
The properties established for $\mg^*$ in the first paragraph of the proof of
Theorem \ref{th:4band}(i) are also satisfied by $\mh^*-\rho\,\mk^*$ for all $\rho\geq0$; and
reasoning as done in that proof for $\ml_{\lb_2}$, we check that $\{\bar\muk_0\}$
and $\{\bar\mm_0\}$ are hyperbolic copies of the base for $x'=\mh^*(\wt,x)$,
respectively attractive and repulsive.
Let $\mu>0$ be a radius of uniform exponential stability (at $-\infty$) of $\bar\mm_0$.
All the orbits for $x'=\mh^*(\wt,x)$ with initial data at 
$\bigcup_{\w\in\W}(\{\w\}\times[\bar\mm_0(\w)-\mu,\bar\mm_0(\w)+\mu])$ are globally
backward defined because of the stability of $\bar\mm_0$ and globally forward defined due to coerciveness.
So, the closure of the lower equilibria $\bar\ml_0^*\leq\bar\mm_0-\mu$ of the attractor
for $x'=\mh^*(\wt,x)$ is a compact invariant set projecting onto $\W$ and below $\{\bar\mm_0\}$, in turn
below $\{\bar\muk_0\}$, and hence \cite[Theorem 5.3]{dno4} proves the
existence of three hyperbolic copies of the base, always for $x'=\mh^*(\wt,x)$.

So, we can repeat some of the steps of the proofs of \cite[Theorems 4.4 and 4.5]{dno5} 
(and \cite[Theorem B.3(iii)]{dlo1}), working with the family
$x'=\mh^*(\wt,x)-\rho\,\mk^*(\wt,x)$. This requires to replace the perturbation 
$\lb\,\mg(\w)$ with $-\rho\,\mk^*(\w,x)$, which can be done with minor modifications 
focusing on $\rho>0$ and on the behavior
within $\mB^*$. The interval $\mI^+:=\{\mu>0\,|\;$
there are three hyperbolic copies of the base $\bar\ml^*_\rho<\bar\mm^*_\rho<\bar\muk^*_\rho$ for all $\rho\in[0,\mu)\}$
is finite and of the form $[0, \rho^+)$, which we prove: checking that $\bar\mm_0^*=\bar\mm_0\leq\bar\mm_\rho^*<\bar\muk_\rho^*\leq\bar\muk_0^*=\bar\muk_0$,
for $\rho\in\mI$; showing that,
for sufficiently large $\rho\geq0$, bounded solutions cannot exist within $\mB^*$, 
and applying Theorem \ref{th:2pers} to show that $\rho^+\notin\mI$.
We take $\bar\w=(h,k)$ and deduce from the previous inequalities
that $m_\rho(t):=\bar\mm_\rho^*(\bwt)$ and $u_\rho(t):=\bar\muk_\rho^*(\bwt)$
are bounded solutions of \eqref{eq:5ct}$^\rho$ for $\rho\in\mI$. Since they are uniformly 
separated, we have tracking. In addition, repeating the argument of
\cite[Theorem 4.5]{dno5}, we can check that $u_{\rho^+}(t)$ and
$m_{\rho^+}(t):=\lim_{\rho\to(\rho^+)^+}m_\rho(t)$ are not uniformly separated.
This fact
allows us to reason as in \cite[Theorem B.3(iii)]{dlo1} in order to prove that,
if $\rho>\rho^+$, then there exists $t_0\in\R$ such that $u_\rho(t_0)<\bar\mm_0(\bar\w{\cdot}t_0)$. 
From here, the last argument of the proof of Theorem \ref{th:5ctconcave} can be repeated 
to prove that we have tipping.
\end{proof}
%%%%%%%%%%%%%%%%%%%%%%%%%%%%%%%%%%%%%%%%%%%%%%%%%%%%%%%%%%%%%%%%%%%%%%%
%%%%%%%%%%%%%%%%%%%%%%%%%%%%%%%%%%%%%%%%%%%%%%%%%%%%%%%%%%%%%%%%%%%%%%%
%%%%%%%%%%%%%%%%%%%%%%%%%%%%%%%%%%%%%%%%%%%%%%%%%%%%%%%%%%%%%%%%%%%%%%%
%%%%%%%%%%%%%%%%%%%%%%%%%%%%%%%%%%%%%%%%%%%%%%%%%%%%%%%%%%%%%%%%%%%%%%%
\section{Numerical simulation and model fitting}\label{sec:6}
In this section, we show that concave-convex functions provide an adequate
fit to laboratory experimental data on the growth rates of certain mosquito
populations that show a strong Allee effect, thus illustrating the real-world
applicability of concave-convex modeling. This straightforward example leads
to autonomous concave-convex equations. More complex situations which are not
dealt with here, where species evolution is influenced by time-varying
external parameters, would result in nonautonomous concave-convex equations.
We also illustrate some of the paper's results with a numerical simulation of
a transition equation that is concave-convex in measure but not asymptotically
concave-convex, showing the wide range of models covered under the theory
developed in this paper.
%%%%%%%%%%%%%%%%%%%%%%%%%%%%%%%%%%%%%%%%%%%%%%%%%%%%%%%%%%%%%%%%%%%%%%%
%%%%%%%%%%%%%%%%%%%%%%%%%%%%%%%%%%%%%%%%%%%%%%%%%%%%%%%%%%%%%%%%%%%%%%%
\subsection{Modeling mosquito growth rate with concave-convex splines}
In this subsection, we utilize some of the data subsets provided by the Biological Resource 
Center EP-Coll in the database \cite{mosquitodata}, specifically some of those identified in 
\cite{vercken1} as exhibiting a strong Allee effect.
Our aim is to demonstrate, using these laboratory data, that certain sets of concave-convex splines 
can provide an adequate fit to the experimental data. The model fitted to these data, obtained in 
a controlled laboratory environment, will be an autonomous model.

The starting point is a dataset of the form $(x_j, y_j)$ for a single population, representing 
the population growth rate $y_j$ as a function of the population size $x_j$, where $j$ varies 
across the set of experiments. Knowing that the population dynamics are subject to the Allee effect, 
we seek a concave-convex cubic spline $x \mapsto \theta(x)$ that appropriately fits the point 
cloud $(x_j, y_j)$. In this way, we understand that the equation $x' = \theta(x)$ adequately 
describes the population dynamics.

As indicated in \cite{meyer1}, regression splines that impose shape constraints (such as monotonicity, 
concavity, or convexity) are highly robust with respect to the choice of spline knots, making them 
particularly suitable when the underlying law generating the data is presumed to possess such a shape constraint.
In this case, we will employ a similar approach to that of \cite{meyer1} to provide a set of concave-convex
splines, which we will fit to the data using least-squares regression with nonnegativeness constraints
on certain parameters.

\begin{figure}
     \centering
         \includegraphics[width=\textwidth]{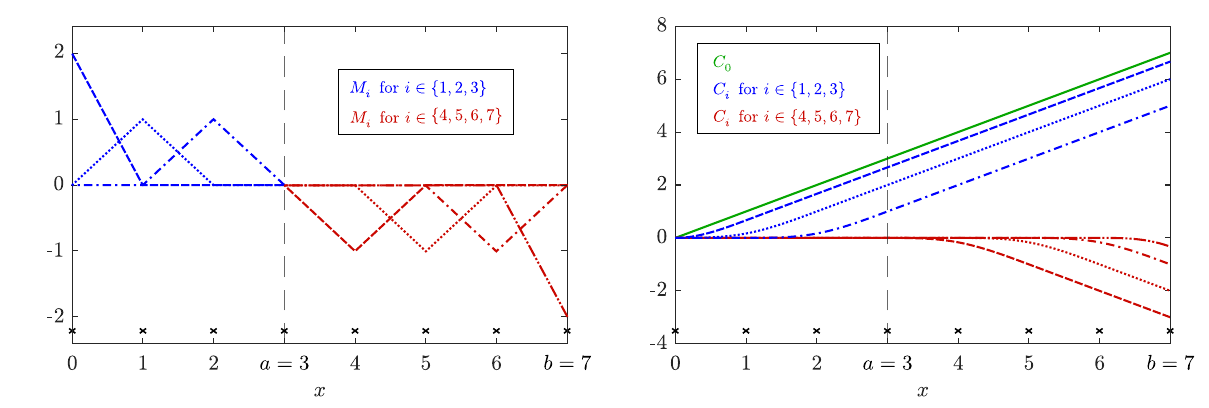}
         \caption{Let $m=2$, $n=3$, $a=3$, $b=7$.
         In the left panel, the linear splines $M_i$, in blue for $i\in\{1,2,3\}$ and in red for $i\in\{4,5,6,7\}$.
         In the right panel, the convex cubic splines $C_i$ for $i\in\{1,2,3\}$ (in blue),
         the concave cubic splines $C_i$ for $i\in\{4,5,6,7\}$ (in red) and the linear spline $C_0(x)=x$ (in green).
         In both panels, the knots are marked with crosses.
}
        \label{fig:splines}
\end{figure}

Let us explain the construction of $\theta$, which will be defined on an interval $[0,b]$,
concave in $[a, b]$ and convex in $[0, a]$
for a point $a\in(0,b)$. First, on $[0,a]$, we consider the $m+2$ equally spaced knots $t_i=i\,h$ for $i=0,\ldots,m+1$,
where $h:=a/(m+1)$ (so, $t_0=0$ and $t_{m+1}=a$), and define $M_i$ for $i=2,\ldots,m+1$
as the continuous piecewise linear maps on $[0,b]$
taking the value $1/h$ at $t_{i-1}$ and $0$ outside $(t_{i-2},t_i)$.
We define $M_1$ in a similar way, with
$M_1(t_0)=M_1(0)=2/h$ and $M_1(x)=0$ for all $x\in[t_1,b]$.
See the blue linear splines in the left panel of Figure \ref{fig:splines}.
Second, we define the convex cubic splines $C_i$ for $i\in\{1,\dots,m+1\}$ as
\begin{equation}\label{def:6spline}
 C_i(x)=\int_0^x\int_0^y M_i(t)\,dt\,dy\,\,\quad\text{for }x\in[0,b]\,.
\end{equation}
Note that $C_i''(x)=M_i(x)$, and hence $C_i$ is strictly convex in the interval in which $M_i>0$ and
linear otherwise, i.e., when $M_i=0$, as in $[a,b]$ (where all of them have slope $1=\int_0^a M_i(t)\,dt$).
See the blue splines in the right panel of Figure \ref{fig:splines}. We also define $C_0(x)=x$.

We proceed in a similar way to define concave cubic splines on
$[0,b]$. First, we divide $[a,b]$ in $n+1$ subintervals of length $\tilde h:=(b-a)/(n+1)$,
with end-points $a=t_{m+1}<t_{m+2}<\cdots<t_{m+n+2}=b$, define
$M_i$ for $i=m+2,\ldots,m+n+1$ as the continuous piecewise linear maps on $[0,b]$
taking the value $-1/\tilde h$ at $t_i$ and $0$ outside $(t_{i-1},t_{i+1})$,
and define $M_{m+n+2}$ similarly, with $M_{m+n+2}(t_{m+n+2})=M_{m+n+2}(b)=-2/\tilde h$
and $M_{m+n+2}(x)=0$ for all $x\in[0,t_{m+n+1}]$. See the red
splines in the left panel of Figure \ref{fig:splines}. Now, we define
$C_{m+2},\ldots,C_{m+n+2}$ by \eqref{def:6spline}, and observe that
they vanish (at least) in $[0,a]$, are globally concave, and are strictly concave when
$M_i$ is strictly negative. See the red splines in the right panel of Figure \ref{fig:splines}.
Note that $C_i(0)=0$ for $i\in\{0,1,\dots,m+n+2\}$.
This condition will eventually represent that the growth rate of the species under
study becomes null when the population is extinct.

Now, let $\{(x_j, y_j)\,|\;j\in\{1,2,\dots,N\}\}$ be the sccater plot that we
aim to fit using a suitable combination of the described cubic splines. The first step is observing
the data in order to choose the parameters $m$, $n$, $a$, and $b$. The second one is looking for
$\alpha_i$ for $i\in\{0,1,\ldots,m+n+2\}$ with the constrain $\alpha_i\ge 0$ for $i\ge 1$, such that, if
\[
 \theta(x):=\sum_{i=0}^{n+m+2}\alpha_i\,C_i(x)\,,
\]
then the squared error function $\sum_{j=1}^N(y_j-\theta(x_j))^2$ is minimum. That is,
we minimize the distance of the point cloud to the
{\em concave-convex regression spline\/} $\theta$ in the least squares sense.
If the vectors $(C_i(x_1),C_i(x_2),\dots,C_i(x_N))$ for
$i\in\{0,1,\ldots,m+n+2\}$ are linearly independent, then the possible solutions
for the parameters are at most finite. We carry out this process with
the \texttt{lsqlin} algorithm from MATLAB R2023a.

\begin{figure}
     \centering
         \includegraphics[width=\textwidth]{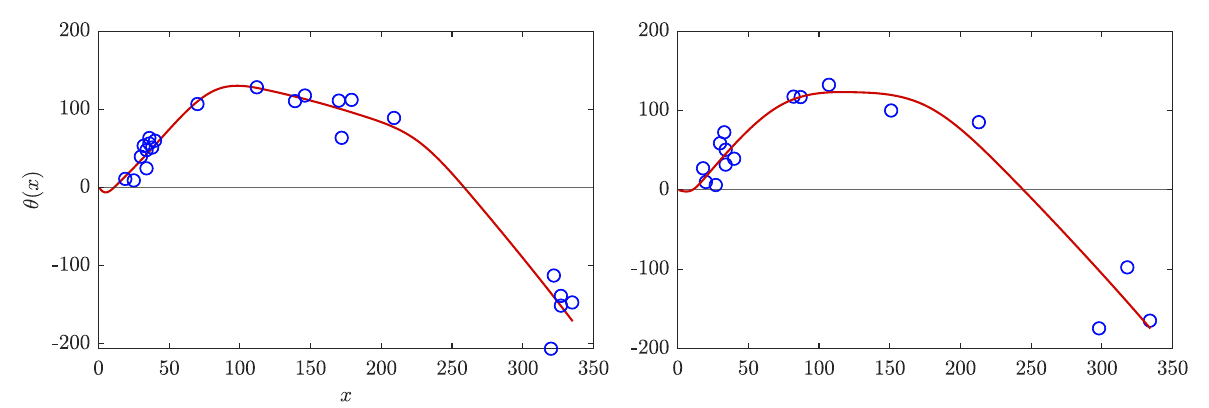}
         \caption{Concave-convex regression splines for two datasets of parasitoid wasps of the genus 
         \emph{Trichogramma Cacoeciae} growth rates.
         Datasets are taken from \cite{mosquitodata}.
         In the left panel (resp. right panel), the dataset \texttt{cac5} (resp. \texttt{cac7}) with $a=45$, 
         $b=335$, $m=3$ and $n=7$ (resp. $a=18.2$, $b=335$, $m=1$ and $n=5$).
         In \cite{vercken1}, both datasets are identified as exhibiting strong Allee effect, in correspondence 
         with the three real roots of each of the splines obtained.
         To ensure that the regression spline is strictly concave on $[a, b]$ and strictly convex on $[0, a]$, the
         additional condition $\alpha_i \ge 0.1$ for $i \in \{1,2,\dots, m+n+2\}$ has been imposed for the minimization.
         In both cases, it has been checked that the vectors $(C_i(x_1),C_i(x_2),\dots,C_i(x_N))$ for $i\in\{0,1,\dots,m+n+2\}$
         are linearly independent.
}
        \label{fig:splines2}
\end{figure}

Now, let us apply the previously described method to some laboratory data.
Each of the experiments of the datasets \texttt{cac5} and \texttt{cac7} from \cite{mosquitodata},
which are reported in \cite{vercken1} to exhibit a strong Allee effect, provide us with pairs of
population sizes $P_j(t)$ and $P_j(t+1)$, each of which corresponds to two successive generations of the insect under study.
Data corresponding to extinction in the generation $P_j(t+1)$ have been excluded.
We interpret that a fixed unit of time has elapsed between the measurement of population size
$P_j(t)$ and $P_j(t+1)$ in all the experiments. If we assume that the growth during this time interval is
governed by $x'=r(x)\,x$ and it is assumed that the per capita growth rate function $r(x)\equiv r_j$ remains
constant within the interval $[P_j(t),P_j(t+1)]$, then $r_j=\log(P_j(t+1)/P_j(t))$.
This approximation, commonly used in biology, is reliable as long as the time unit elapsed between the
measurement of $P_j(t)$ and $P_j(t+1)$ is sufficiently small. Thus, the dataset that we want to adjust through
concave-convex splines encompass data $(x_j,y_j)$ of the form $(P_j(t),P_j(t)\,\log(P_j(t+1)/P_j(t)))$
which approximate the pair (population size, population growth rate). The index $j$
varies on the set of different experiments,

By selecting appropriate values for $a$, $b$, $m$, and $n$, we obtain the concave-convex regression
splines depicted in Figure \ref{fig:splines2} for these datasets.
The strong Allee effect mentioned in \cite{vercken1} is reflected in the fact that the obtained
regression splines have three roots: the middle one provides the Allee cooperation threshold.
It should be noted that the resulting regression splines fit the data reasonably well and that the regression
splines for both datasets are similar to each other, as might be expected given that these datasets pertain
to the same species, albeit collected from different plants or locations.
Moreover, the slight variation between the data and results of \texttt{cac5} and \texttt{cac7}
shows the influence of external factors on the growth rate (climate and other external factors
vary from one location to another), while also motivating the interest of a more complex and
not addressed here nonautonomous modeling approach (since these factors also change over time
within the same plant and location).
Finally, it is worth noting that further research is needed on how to select
the inflection point $a$ from concavity to convexity based on the dataset, and that
minimizing the error by allowing unrestricted variation of this
parameter $a$ is not a good option, since it may compromise the concave-convex shape
constraint of the objective function.
%%%%%%%%%%%%%%%%%%%%%%%%%%%%%%%%%%%%%%%%%%%%%%%%%%%%%%%%%%%%%%%%%%%%%%%
%%%%%%%%%%%%%%%%%%%%%%%%%%%%%%%%%%%%%%%%%%%%%%%%%%%%%%%%%%%%%%%%%%%%%%%
\subsection{Numerical simulation of a concave-convex critical transition.}\label{subsec:6modelo}
In this section, we present a numerical example of a critical transition within the
parametric framework described by Theorem~\ref{th:5ctconcave}. At the same time, this example is a
sample of a concave-convex which is not d-concave model for the subsequent choice of
the involved coefficients.
As starting point, we consider the single-species population model
\begin{equation}\label{eq:6}
y'=r\,y\,\left(1-\frac{y}{K(t)}\right)\frac{y-S}{K(t)}-\Delta(t)\,\frac{y^2}{b+y^2}\,.
\end{equation}
The cubic term represents a population subject to the Allee effect, that is,
it incorporates some cooperative mechanism within the model. Here, $r$ is the intrinsic per capita
growth rate, $K(t)$ is the amount over which the growth rate becomes negative in absence of predation
for each $t\in\R$ (closely related to the limit of resources of the environment)
and $S$ is the Allee threshold, i.e., the critical population size below
which the growth rate becomes negative, also in the absence of predation.
On the other hand, $-\Delta(t)\,y^2/(b+y^2)$ is a Holling type III functional response
term representing the contribution of predation to the population net growth rate:
$b$ is related to the amount of time that elapses between attacks, and $\Delta(t)$ is
proportional to the density of predators present in the patch at
each instant of time $t\in\R$.

\begin{figure}
     \centering
         \includegraphics[width=\textwidth]{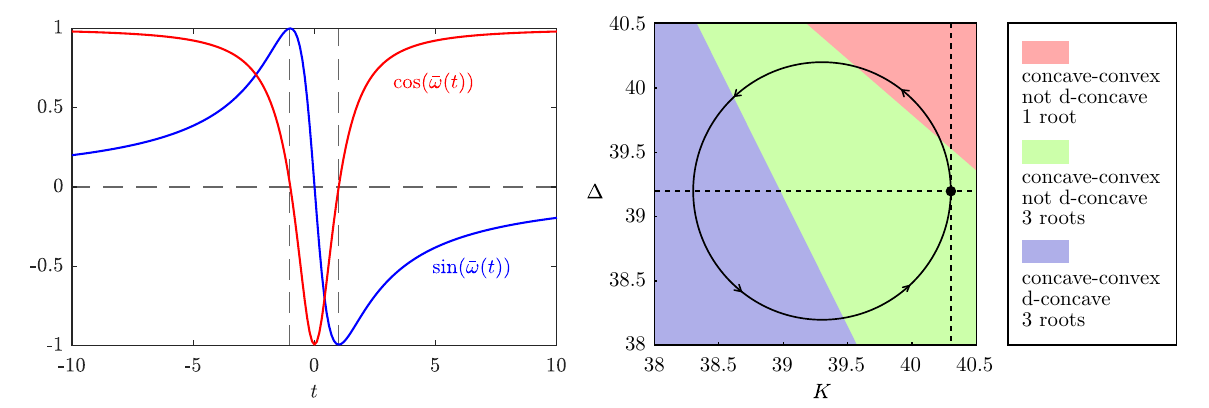}
         \caption{In the left panel, in blue $\sin(\bar\w(t))$ and in red $\cos(\bar\w(t))$, where
         $\bar \w(t)$ is
         the (strictly increasing) solution of $\w'=1-\cos(\w)$ satisfying $\bar\w(-2{\cdot}10^5)=10^{-5}$,
         with limits $0$ at $-\infty$ and $2\pi$ at $+\infty$.
         The map $t\mapsto\sin(\bar\w(t))$ is strictly decreasing between the two dashed vertical
         lines, within the same time interval on which $t\mapsto\cos(\bar\w(t))$ is below 0,
         which is the average value between its maximum 1 and its minimum -1.\\[.1cm]
         In the right panel, a numerical approximation to the map of the properties
         of concavity-convexity, d-concavity and the number of roots of the map
         $y\mapsto r\,y\,(1-y/K)(y-S)/K-\Delta\,y^2/(b+y^2)$ depending on the choices
         of $K$ and $\Delta$ for $r=3$, $S=0.3$ and $b=620$.
         The black circle of radius 1 centered at $(K_0,\Delta_0)=(39.3,39.2)$ represents
         the values that $(\tilde K(\w),\tilde\Delta(\w))$ takes as $\w$ goes
         from $0$ to $2\pi$, with $K_1=\Delta_1=1$. The verification of each property was
         conducted numerically on a $100\times100$ grid, with the edges between regions
         smoothed in the rendering.
         }
        \label{fig:figuremap}
\end{figure}

We assume that the two involved time-dependent parameter functions $K$ and $\Delta$ are drifted
by a common external effect, quantified through the variable $\w$, related to phenomena
such as climate variability or the Earth's rotation, and hence depending on time.
This allows a wide variety of biologically meaningful functions to be modelled in a simple
way. In order to propose a simple model for this influence, we define
$\tilde K(\w):=K_0+K_1\,\cos(\w)$ and $\tilde\Delta(\w):=
\Delta_0+\Delta_1\sin(\w)$, with $K_0>K_1>0$ and $\Delta_0>\Delta_1>0$ to get $\tilde K>0$
and $\tilde\Delta>0$.
Then, $K(t):=\tilde K(\bar\w(t))$ and $\Delta(t):=\tilde\Delta(\bar\w(t))$, where $\bar\w(t)$
stands for the time evolution of the driving variable. For this variable,
we take a simple asymptotically constant variation $\bar\w(t)$ between the
values $0$ and $2\pi$ through the ODE $\w'=1-\cos(\w)$, so that
$\bar\w(t)$ is  monotonically increasing. In this way,
the density of predators (and hence, proportionally, $\Delta(t)$)
increases when $K(t)$ exceeds the average value between its maximum and minimum.
This can be understood as the resource limitation for the species under study also
representing a resource limitation for the predator species.
The left panel of Figure \ref{fig:figuremap} shows
the shape of the components $\cos(\bar\w(t))$ and $\sin(\bar\w(t))$ of $K(t)$ and $\Delta(t)$
for a particular choice of $\bar\w(t)$.
Alternative forms of variation for $\w$, which governs the external effect,
may be considered to examine more realistic situations, as well
as situations with different biological implications.

With our choices, the evolution of the population modeled by \eqref{eq:6} is hence given by
the second component of an orbit of the flow induced by the system of ODEs
\begin{equation}\label{eq:6sistapl}
\left\{\begin{array}{lll}
\w'&=&1-\cos(\w)\,,\\[1ex]
y'&=&\displaystyle{r\,y\,\left(1-\frac{y}{\tilde K(\w)}\right)\frac{y-S}{\tilde K(\w)}-\tilde\Delta(\w)\,\frac{y^2}{b+y^2}}
\end{array}
\right.
\end{equation}
on $\s^1\times\R$, identifying $\mathbb{S}^1$ with the quotient space $\R/2\pi\Z$.
Thus, $\s^1$ is invariant by the action
of the flow $\sigma_{\s^1}$ induced by the first equation $\w'=1-\cos(\w)$, and
we can read the second equation as a family of equations over a skewproduct base $(\s^1,\sigma_{\s^1})$.
We represent the right-hand of the second equation by $\mh(\w,y)$, and note that it can be understood as a
map defined either on $\R\times\R$ or on $\s^1\times\R$.

The first step is to examine the concavity-convexity properties of the second equation in \eqref{eq:6sistapl}.
To do so, note first that $\s^1$ consists of the fixed point $0$ and the heteroclinic orbit $ (0, 2\pi)$,
which corresponds to all the non constant (and strictly increasing) solutions
of $\w'=1-\cos(\w)$.
Thus, $\merg=\{\delta_0\}$, where $\delta_0$ stands for the Dirac measure concentrated at $0$.
Hence, it suffices to study the concavity-convexity properties of $\mh(0,y)$ to deduce
concavity-convexity in full measure.

\begin{figure}
     \centering
         \includegraphics[width=\textwidth]{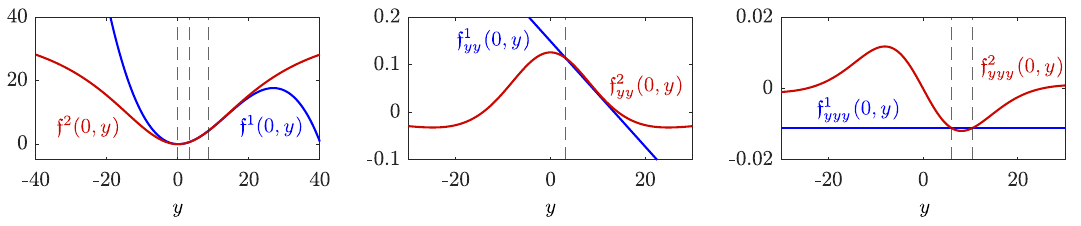}
         \caption{Numerical depiction of the study of the roots of $\mh(0,y)$, $\mh_{yy}(0,y)$ and $\mh_{yyy}(0,y)$.
         We study separately its two terms $\mf^1(0,y)=r\,y\,(1-y/\tilde K(0))(y-S)/\tilde K(0)$ 
         and $\mf^2(0,y)=\tilde\Delta(0)\,y^2/(b+y^2)$.
         In the left panel,
         the vertical dashed lines stand for the intersections of $\mf^1(0,y)$ (in blue) and $\mf^2(0,y)$ (in red);
         i.e., to the three roots of $\mh(0,y)$.
         In the central panel, %the blue curve represents $\mf^1_{yy}(0,y)$ and the red one $\mf^2_{yy}(0,y)$. The
         the vertical dashed line stands for the intersection of $\mf^1_{yy}(0,y)$ (in blue) and $\mf^2_{yy}(0,y)$ (in red);
         i.e., to the unique root $\mb(0)$ of $\mh_{yy}(0,y)$.
         Hence, $\mh(0,y)$ is concave on $[\mb(0),\infty)$ and convex on $(-\infty,\mb(0)]$, that is, it is concave-convex.
         In the right panel, %the blue curve represents $\mf^1_{yyy}(0,y)$ and the red one $\mf^2_{yyy}(0,y)$.
         The vertical dashed lines stand for the intersections of $\mf^1_{yyy}(0,y)$ (in blue) and $\mf^2_{yyy}(0,y)$ (in red);
         i.e., to the two roots of $\mh_{yyy}(0,x)$.
}
        \label{fig:intuition}
\end{figure}

\begin{figure}
     \centering
         \includegraphics[width=\textwidth]{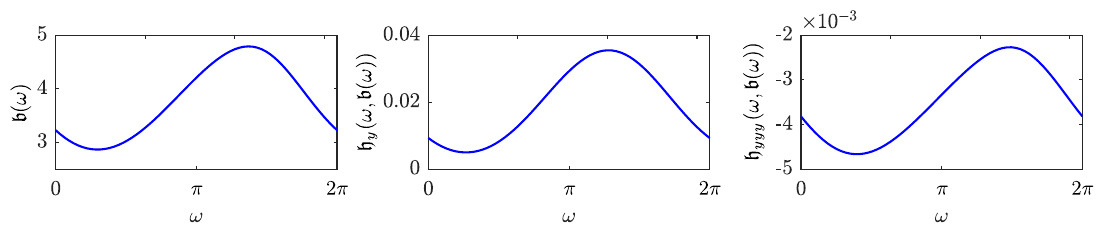}
         \caption{In the left panel, numerical approximation to the curve $\w\mapsto\mb(\w)$ which satisfies
         $\mh_{yy}(\w,\mb(\w))=0$ for all $\w\in\s^1$, for the choices of the parameters in Figure \ref{fig:figuremap}.
         The approximation has been carried out using algorithm \texttt{fzero} from Matlab R2023a on a grid over $[0,2\pi]$.
         In the central and right panels, evaluation of $\mh_y(\w,\mb(\w))$ and $\mh_{yyy}(\w,\mb(\w))$ 
         using the numerical approximation of the left panel, respectively.
}
        \label{fig:comega}
\end{figure}

The next selection of parameters provide a model as classical as \eqref{eq:6sistapl}
driven by a concave-convex ODE that is not d-concave.
Let $r=3$, $S=0.3$, $b=620$, $K_0=39.3$, $K_1=1$, $\Delta_0=39.2$, $\Delta_1=1$.
Note that fixing any $\w\in\R$ means fixing a value of
$(K,\Delta)=(\tilde K(\w),\tilde\Delta(\w))$ which lies in the square
$\mR:=[38.3,\,40,3]\times[38.2,\,40.2]$.
(And that, if the solution $\bar\w$ is not constant,
$\{(\tilde K(\bar\w(t)),\tilde \Delta(\bar\w(t)))\,|\;t\in\R\}$ is the
circumference of radius $1$ centered at $(39.3,\,39.2)$.)
With this choices, we numerically check that: $y\mapsto\mh(0,y)$ has three real roots;
$y\mapsto\mh_{yy}(0,y)$ has a unique real root; and $y\mapsto\mh_{yyy}(0,y)$ does not
have constant sign (see Figure \ref{fig:intuition}). That is, $y\mapsto\mh(0,y)$ has three real
roots, is concave-convex and is not d-concave.
A verification analogous to the one made for $\w=0$ has been conducted on a grid of
values of $(K,\Delta)$ on a square containing $\mR$.
The result is partly depicted in the right panel of Figure \ref{fig:figuremap} (which in particular shows that
our equation \eqref{eq:6} is not globally d-concave). Figure \ref{fig:figuremap} does not depict two fundamental
facts for the next step, also observed numerically at all the points of the grid:
that, for any $\w\in\R$, there exists a unique solution
$y=\mb(\w)>0$ of the algebraic equation $\mh_{yy}(\w,y)=0$ (which hence varies continuously
with respect to the coefficients $r,\,b,\,S,\,\tilde K(\w)$ and $\tilde\Delta(\w)$)
with $\mh_{yy}(\w,y)>0$ for $y<\mb(\w)$ and $\mh_{yy}(\w,y)<0$ for $y>\mb(\w)$;
and that $\mh_{yyy}(\w,\mb(\w))<0$ for all $\w\in\R$. Clearly
$\mb$ is $2\pi$-periodic, and hence the maps $\s^1\to\R,\,\w\mapsto\mb(\w)$
and $\s^1\to\R,\,\w\mapsto \mh_{yyy}(\w,\mb(\w))$ are well-defined and continuous.
These curves are numerically approximated in the left and right panels of Figure \ref{fig:comega}.
The Implicit Function Theorem can be applied to $\RR\to\R,\,(\w,y)\mapsto\mh_{yy}(\w,y)$
to deduce the existence and continuity of the derivative $\dot{\mb}(\w)$, which is also $2\pi$-periodic and hence
well defined on $\s^1$.
As usual, we denote by $\mb'(\w)=(d/dt)\,\mb(\sigma_{\s^1}(t,\w))|_{t=0}$. Is it easy to check that
$\mb'(\w)=\dot{\mb}(\w)\,(1-\cos(\w))$, and hence $\mb'\colon\s^1\to\R$ is continuous.

For purposes which will be clear some paragraphs below,
we make the change of variables $(\w,y)\mapsto(\w,x)$ with $x=y-\mb(\w)$,
which takes the second equation in \eqref{eq:6sistapl} to
$x'=\tilde\mh(\w,x)$, with $\tilde\mh(\w,x):=\mh(\w,x+\mb(\w))-\mb'(\w)=\mc(\w)+\md(\w)\,x+\mg(\w,x)$ for
$\mc(\w):=\mh(\w,\mb(\w))-\mb'(\w)$, $\md(\w):=\mh_y(\w,\mb(\w))$
and $\mg(\w,x):=\mh(\w,x+\mb(\w))-\mh(\w,\mb(\w))-\mh_y(\w,\mb(\w))\,x$.
It is easy to check that
\ref{cc1} is satisfied. The central panel of Figure \ref{fig:comega} numerically
shows that \ref{cc2} also holds. Condition \ref{cc3} follows from $\lim_{y\to\pm\infty}\mh(\w,y)=\mp\infty$
uniformly on $\s^1$. It is clear that, for all $\w\in\s^1$, $\mg(\w,0)=\mg_x(\w,0)=0$, and the choice of
$\mb$ ensures that $\mg_{xx}(\w,x)>0$ whenever $x<0$ and $\mg_{xx}(\w,x)<0$ if $x>0$,
from where \ref{cc4}, \ref{cc5} and \ref{cc6} follow easily.


\vspace{.1cm}

\textbf{The transition equation.}
We consider that the previously described system undergoes a transition induced by an upsurge
in predation,
which subsequently decreases to a level slightly above the initial one.
We study the influence of the size $\rho$ of the increase into the transition.
Thus, we consider
\begin{equation}\label{eq:6trans-sim}
y'=r\,y\,\left(1-\frac{y}{K(t)}\right)\frac{y-S}{K(t)}-\big(\Delta(t)+\rho\,\Gamma(t)\big)\frac{y^2}{b+y^2}\:,
\end{equation}
where $\G(t):=15\arctan(t+10)/\pi-14\arctan(t-10)/\pi+1/2$ is a strictly positive
nonmonotonic asymptotically constant function which reaches its maximum value
(which is slightly above 14 and represents the upsurge) for $t\approx0.174$,
and such that $\lim_{t\to-\infty}\G(t)=0<1=\lim_{t\to\infty}\G(t)$
(so that the subsequent level is slightly above the initial one).
From now on, we fix a (typical) solution $\bar\w(t)$ of $\w'=1-\cos(\w)$ satisfying
$\lim_{t\to-\infty}\bar\w(t)=0$ and $\lim_{t\to\infty}\bar\w(t)=2\pi$.
In particular, that with $\bar\w(-2{\cdot}10^5)=10^{-5}$.
Recall that, in our model \eqref{eq:6trans-sim},
$K(t):=\tilde K(\bar\w(t))$ and $\Delta(t):=\tilde\Delta(\bar\w(t))$.

Our plan is to apply Theorem \ref{th:5ctconcave} (or, more precisely, the information
provided by Remark \ref{rm:5barrho}) to deduce the existence of a
critical transition for \eqref{eq:6trans-sim}. This cannot be done directly:
although, later, we will conduct the simulations using the variable $y$, we must use the previously
introduced variable $x=y-\mb(\w)$ shifting
the change-from-concavity-to-convexity curve $\mb$ to $0$ in order to get an equation
satisfying the required hypotheses of Theorem \ref{th:5ctconcave}: otherwise, the condition
on $\mg$ included in \ref{ct1} would not hold.
We also change the perturbation term $-\rho\,\Gamma(t)\,y^2/(b+y^2)$ by $-\rho\,\Gamma(t)\,f(y)$,
with $f(y):=y^2/(b+y^2)$ for $y\ge 0$ and $f(y):=0$ for $y<0$ (in order to get \ref{ct2} and \ref{ct3},
as explained below).
Note that this change does not modify the structure and properties of the set of nonnegative solutions,
which are the biologically significant ones. So, the change of variables
takes the (modified) equation \eqref{eq:6trans-sim}$^\rho$ to an ODE of the form \eqref{eq:5ct}$^\rho$ for
$h(t,x):=\mh(\bar\w(t),x+\mb(\bar\w(t)))-\mb'(\bar\w(t))$ and $k(t,x):=\G(t)\,f(x+\mb(\bar\w(t)))$.

\begin{figure}
     \centering
         \includegraphics[width=\textwidth]{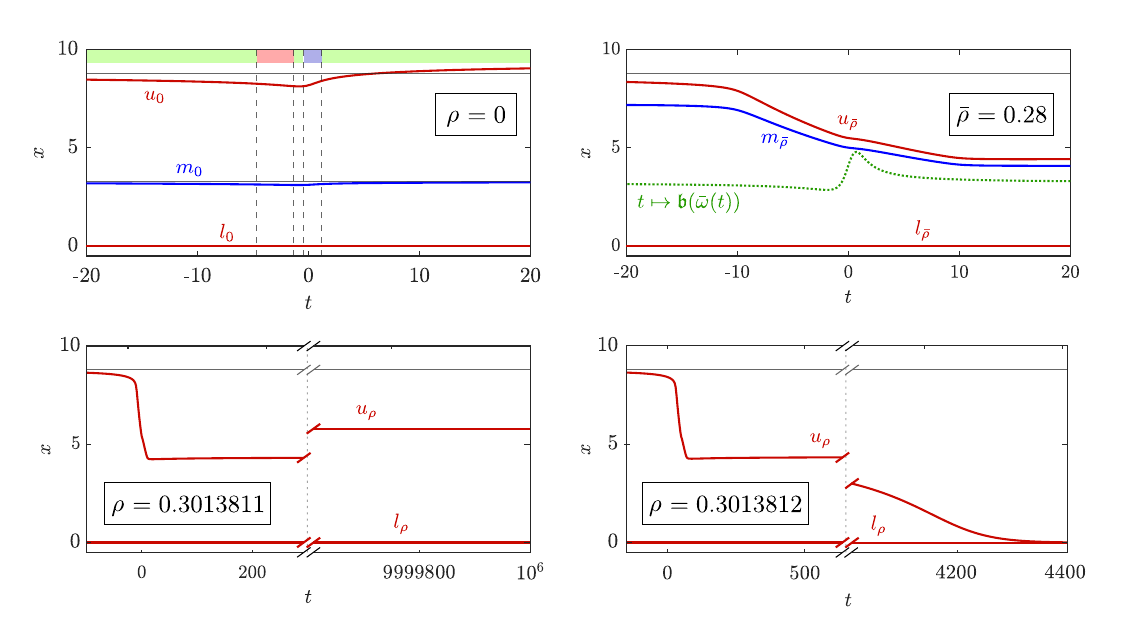}
         \caption{Numerical approximations (in red) to the locally pullback attractive solutions
         $l_{\rho}$ and $u_{\rho}$ of
         \eqref{eq:6trans-sim}$^\rho$ for different values of the parameter.
         In all the simulations, the map $\bar\w(t)$ use to define
         $K(t):=\tilde K(\bar\w(t))$ and $\Delta(t):=\tilde\Delta(\bar\w(t))$
         satisfies the initial condition $\bar\w(-2{\cdot}10^5)=10^{-5}$.\\[.1cm]
         In the top left panel, the solid black horizontal lines approximate the three hyperbolic fixed points of
         $y'=\mh(0,y)$. In the rest of the panels, the black horizontal lines point out the attractive
         hyperbolic fixed points of $y'=\mh(0,y)$.
         The two locally pullback attractive solutions (in red) are hyperbolic attractive solutions of
         \eqref{eq:6trans-sim}$^0$, as they asymptotically connect stable fixed points of $y'=\mh(0,y)$.
         The third hyperbolic solution,
         which is repulsive, is shown in blue. The dashed vertical lines indicate the points where $\bar\w(t)$
         transitions between the different regions in the right panel of Figure \ref{fig:figuremap}: the color band
         at the top shows the region where $\bar\w(t)$ it is located, using the same color code as in the right panel of
         Figure \ref{fig:figuremap}.
         In the top right panel, the three hyperbolic solutions for $\bar\rho=0.28$ with the same
         color code of the previous panel. As explained in the text, they approach the fixed points of
         $y'=\mh(0,y)$ as $t\to-\infty$ and those of $y'=\mh(0,y)-\bar\rho\,y^2/(b+y^2)$ as $t\to\infty$.
         \\[.1cm]
         In the two bottom panels, the two locally pullback attractive solutions (in red) of
         \eqref{eq:6trans-sim}$^\rho$ are shown for two values of $\rho$ very close to the bifurcation point
         $\rho_0$. Note that the time axis has been broken to capture the most representative parts
         for our study. The numerical approximation of the bifurcation point was obtained by taking a very
         large finite time, $10^6$, with $\varepsilon = 10^{-3}$, and examining whether the locally pullback
         attractive solutions at $10^6$ are $\epsilon$-close or not. This criterion was used to approximate
         the tipping point (at finite time), found between $0.3342651$ and $0.3342652$, via a bisection-like method.
         The left panel shows a tracking scenario at finite time $10^6$, while the right panel
         shows a tipping scenario.
         \\[.1cm]
         Numerical integration has been carried out with \texttt{ode89} from Matlab R2023a.
}
        \label{fig:ahacer}
\end{figure}

Let us check a ``suitably modified" \ref{ct1}.
We already know that $\mc$, $\md$ and $\mg$ satisfy \ref{cc1}-\ref{cc6}.
The second panel of Figure \ref{fig:ahacer} shows numerical evidence of the existence of
three hyperbolic solutions $l_{\bar\rho}=0<m_{\bar\rho}<u_{\bar\rho}$ of \eqref{eq:6trans-sim}$^{\bar\rho}$
for $\bar\rho=0.28$. Their asymptotic limits as $t\to-\infty$ are the three roots $c^-_l=0<c^-_m<c^-_u$ of
$\mh(0,y)$ already found in the left panel of Figure \ref{fig:intuition}
(which define three constant hyperbolic solutions of $y'=\mh(0,y)$, with $0$ and $c^-_u$ bounding
the set of bounded solutions). Their asymptotic limits as $t\to\infty$ are the three roots
$c^+_l=0<c^+_m<c^+_u$ of $\mh(0,y)-\bar\rho\,y^2/(b+y^2)$ (whose existence can be
numerically checked, and which play the same role as those for $\rho=0$).
The hyperbolicity of these solutions is guaranteed by this
limiting behavior: see \cite[Lemma 4.15]{duenasphdthesis}, which extends \cite[Lemma 3.8]{dno3}.
The arguments of Proposition \ref{prop:5pull} can be repeated to check another observed fact: that
$0$ and $u_{\bar\rho}$ are the lower and upper bounded solutions.
Let us check that these three maps provide hyperbolic copies of a suitable base.
We define $\W:=\mathrm{closure}\{(\bar\w(s),\G(s))\,|\;s\in\R\}$ in $\s^1\times\R$. It is easy to check that
$\W=\{(0,0)\}\cup\{(\bar\w(s),\G(s))\,|\;s\in\R\}\cup\{(0,1)\}$, and that the map $\sigma$ given by
$\sigma(t,0,0):=(0,0)$, $\sigma(t,\bar\w(s),\G(s)):=(\bar\w(t+s),\G(t+s))$ and $\sigma(t,0,1):=(0,1)$
is a continuous flow. This set will play the role played by the common hull in Section \ref{sec:5}.
The maps $\mh$ (resp.~$k$) is extended to a map on $\WR$ by $(\w,\gamma,y)\mapsto\mh(\w,y)$
(resp.~by $(0,0,y)\mapsto 0,\,(0,1,y)\mapsto f(t)$ and $(\bar\w(s),\G(s),y)\mapsto\G(s)\,f(y)$).
We focus on $u_{\bar\rho}$ and define $\mK_{\bar\rho}^u:=
\mathrm{closure}\{(\bar\w(t),\G(t),u_{\bar\rho}(t))\,|\;t\in\R\}\subset\W\times\R$,
which is given by the orbit plus the points $(0,0,c_u^-)$ and $(0,1,c_u^+)$.
So, $\mK_{\bar\rho}^u$ is a copy of the base, say $\mK_{\bar\rho}^u=\{\tilde\muk_{\bar\rho}\}$.
It is also clear that it is hyperbolic attractive (see, e.g., \cite[Proposition 1.54]{duenasphdthesis}).
We define $\{\bar\mm_{\tilde\rho}\}$ (which is repulsive) and $\{\tilde\ml_{\bar\rho}\}$ (again attractive) in
an analogous way. Then, $\tilde\ml_{\bar\rho}=0$ and $\tilde\muk_{\bar\rho}$ are the lower and upper equilibria of the
global attractor for the equation extended to $\W$ corresponding to $\bar\rho$.
A numerical computation with a large number of initial data is
conducted to check that $\tilde\mm_{\bar\rho}$ plays the role of the boundary between the basins of attraction of the other
copies of the base, which in particular ensures that there are no more hyperbolic copies of the base.
Observe also that the change of variables $\WR\to\WR,\,(\w,\gamma,y)\mapsto(\w,\gamma,x)$ with $x=y-\mb(\w)$
takes this hyperbolic structure to the same one for the transformed family of equations over $\W$:
the maps $\bar\ml_{\bar\rho}(\w,\gamma):=\tilde\ml_{\bar\rho}(\w,\gamma)-\mb(\w)=-\mb(\w)$,
$\bar\mm_{\bar\rho}(\w,\gamma):=\tilde\mm_{\bar\rho}(\w,\gamma)-\mb(\w)$ and
$\bar\muk_{\bar\rho}(\w,\gamma):=\tilde\muk_{\bar\rho}(\w,\gamma)-\mb(\w)$ define the adequate
hyperbolic copies of the base. Altogether, we have checked
the condition equivalent to \ref{ct1} required to reach the conclusion of Theorem \ref{th:5ctconcave}:
see Remark \ref{rm:5barrho}.

Let us now check \ref{ct2} and \ref{ct3}.
We have $k(t,x)\geq0$ for all $(t,x)\in\RR$. The definition of $\G$ guarantees
that $\lim_{t\to-\infty}k(t,x)=0$ uniformly on compact sets.
And it is clear that $\lim_{x\to-\infty}k(t,x)=0$ uniformly on $\R$,
that $\tilde\ml_0(\w,\gamma)=-\mb(\w)<0$ for all $(\w,\gamma)\in\W$,
and hence that $k_x(t,x)=0$ whenever $x\le\tilde\ml_0(\sigma(t,\w,\gamma))<0$.

It remains to check the additional hypotheses of the statement of Theorem
\ref{th:5ctconcave}, always for $\bar\rho=0.28$. The first one reduces to
$\tilde\mm_{\bar\rho}(\w,\gamma)>\mb(\w)$ for all $(\w,\gamma)\in\W$; i.e., to
$m_{\bar\rho}(t)>\mb(\bar\w(t))$ for $t\in\R$, $c_m^->\mb(0)$ and $c_m^+>\mb(0)$. This
can be numerically checked: see the second panel in Figure \ref{fig:ahacer}.
Finally, we take into account that $\Gamma(t)\,f(y)>0$ for all $t\in\R$ and $y>0$, and
that $(\partial^2/\partial y^2) f(y)>0$ for $y\in(0,2\sqrt{155/3})\approx(0,14.3759)$.
The second panel of Figure \ref{fig:ahacer}, combined with the asymptotic limits of the represented
solutions (which allow us to infer the bounds
outside a finite window), shows that this convexity area comprises the region $\{(t,y)\in\RR\,|\;
m_{\bar\rho}(t)\le y\le u_{\bar\rho}(t)\}$ as well as the asymptotic limits at $-\infty$ and $+\infty$,
which (after the change of variable) ensures the last conditions required in Theorem \ref{th:5ctconcave}.

Theorem \ref{th:5ctconcave} ensures the existence of $\rho_0>0$ such that
\eqref{eq:6trans-sim}$_\rho$ exhibits tracking for $\rho\in[\bar\rho,\rho_0]$ and
tipping for $\rho>\rho_0$ (recall Definition \ref{def:5trackingtipping}). It is easy
to check that $\rho\to u_\rho$ decreases, and clearly $l_\rho\le 0=l_{\rho_0}$
for all $\rho$. Hence, there is tracking for $\rho\in[0,\rho_0]$. This change of dynamics
can be observed in the panels of Figure \ref{fig:ahacer} for different parameter values.
In this figure, the bifurcation point has been approximated using a bisection-like method, verifying
whether the numerical approximations of the local pullback attractors remain separated by a
tolerance $\varepsilon$ after a finite but sufficiently long integration time $10^6$ set a priori.
In this example, understanding that the upper locally pullback attractive solution represents
the population size under study throughout the transition, the tipping situation means
extinction, with the population size approaching zero,
while tracking signifies persistence.
So, the application of the theorem, illustrated in the numerical simulation, demonstrates that
there exists a specific threshold size $\rho_0$ for the upsurge in predation: an increase beyond
this threshold will lead to the extinction of the species, whereas an increase below the threshold
allows for its survival.
%%%%%%%%%%%%%%%%%%%%%%%%%%%%%%%%%%%%%%%%%%%%%%%%%%%%%%%%%%%%%%%%%%%%%%%%%%%%%
%%%%%%%%%%%%%%%%%%%%%%%%%%%%%%%%%%%%%%%%%%%%%%%%%%%%%%%%%%%%%%%%%%%%%%%%%%%%%
%%%%%%%%%%%%%%%%%%%%%%%%%%%%%%%%%%%%%%%%%%%%%%%%%%%%%%%%%%%%%%%%%%%%%%%%%%%%%
%%%%%%%%%%%%%%%%%%%%%%%%%%%%%%%%%%%%%%%%%%%%%%%%%%%%%%%%%%%%%%%%%%%%%%%%%%%%%

\end{document}